\documentclass[a4paper]{article}

\usepackage[margin=1in]{geometry}

\usepackage[english]{babel}
\usepackage{csquotes}
\usepackage[url=false]{biblatex}
\usepackage[
    pdfauthor={Lukas Pieper},
    pdftitle={Discrete Temporal CSPs et al.},
]{hyperref}
\usepackage{microtype}

\usepackage{mathtools}
\usepackage{amssymb}
\usepackage{amsthm}

\usepackage{xcolor}

\usepackage{tikz}
\usetikzlibrary{cd}

\theoremstyle{plain}
\newtheorem{theorem}{Theorem}
\newtheorem{proposition}[theorem]{Proposition}
\newtheorem{lemma}[theorem]{Lemma}
\newtheorem{corollary}[theorem]{Corollary}
\newtheorem{observation}[theorem]{Observation}
\newtheorem{fact}[theorem]{Fact}
\theoremstyle{remark}
\newtheorem{remark}[theorem]{Remark}
\newtheorem{example}[theorem]{Example}
\theoremstyle{definition}
\newtheorem{definition}[theorem]{Definition}

\DeclareMathOperator{\Csp}{CSP}
\DeclareMathOperator{\Pcsp}{PCSP}
\DeclareMathOperator{\Pol}{Pol}
\DeclareMathOperator{\Succ}{succ}
\DeclareMathOperator{\Min}{min}
\DeclareMathOperator{\Max}{max}
\DeclareMathOperator{\Pow}{\mathcal{P}}
\DeclareMathOperator{\Proj}{Proj}
\DeclareMathOperator{\Diff}{diff}
\DeclareMathOperator{\Dist}{dist}
\DeclareMathOperator{\Th}{Th}
\DeclareMathOperator{\Diag}{Diag}

\DeclareMathOperator{\Sep}{Sep}
\DeclareMathOperator{\Cycl}{Cycl}

\DeclareMathOperator{\Hom}{Hom}
\DeclareMathOperator{\Aut}{Aut}
\DeclareMathOperator{\End}{End}

\DeclareMathOperator{\Id}{id}
\DeclareMathOperator{\Acl}{acl}

\newcommand*{\coloniff}{\;\vcentcolon\Longleftrightarrow\;}

\newcommand*{\Ptime}{\mathsf{P}}
\newcommand*{\Np}{\mathsf{NP}}
\newcommand*{\Npc}{\Np\textnormal{-}\mathrm{c}}

\newcommand*{\Sac}{{\mdseries\rmfamily\upshape\scshape SinglAC}}
\newcommand*{\Csac}{{\mdseries\rmfamily\upshape\scshape CSinglAC}}
\newcommand*{\Saip}{{\mdseries\rmfamily\upshape\scshape SinglAIP}}

\title{Discrete Temporal CSPs et al.}
\author{
    Lukas Pieper\thanks{
        Funded by the European Union (ERC, POCOCOP, 101071674).
        Views and opinions expressed are however those of the author(s) only and do not necessarily reflect those of the European Union or the European Research Council Executive Agency.
        Neither the European Union nor the granting authority can be held responsible for them.
    }
}
\date{2026}

\begin{document}
\maketitle
\begin{abstract}
    We revisit the complexity classification of constraint satisfaction problems (CSPs) for first-order reducts of $(\mathbb{Z},<)$, established by \citeauthor{bodirsky2018} in \citeyear{bodirsky2018}.
    Following their approach, we give a self-contained proof of this classification, correcting several errors in the argument.
    We further show that, with the exception of cases arising from finite-domain CSPs, all tractable cases are solved by a well-known uniform algorithm, eliminating the need for the collection of ad-hoc algorithms used in the original paper.

    A central component of the proof allows the analysis to be restricted to first-order reducts of $(\mathbb{Z},<)$ that contain the successor relation.
    We substantially simplify this reduction and generalize the underlying construction.
    This generalization allows us to apply the same technique to first-order reducts of the triangular lattice graph and of infinite regular trees, yielding a complete complexity classification for these classes of structures, in particular a $\Ptime$ versus $\Np$-complete dichotomy.
\end{abstract}

\section{Introduction}
The constraint satisfaction problem $\Csp(\mathcal{M})$ of a relational structure $\mathcal{M}$ asks, given a primitive positive sentence $\phi$ over the signature of $\mathcal{M}$, that is, an existentially quantified conjunction of atomic formulas, whether $\phi$ holds in $\mathcal{M}$.
Equivalently, viewing $\phi$ as encoding, in the natural way, a structure $\mathcal{N}$ over the same signature, this is the question of whether there exists a homomorphism from $\mathcal{N}$ to $\mathcal{M}$.
When $\mathcal{M}$ is finite, the complexity of $\Csp(\mathcal{M})$ is by now completely understood: the celebrated proof of the Feder--Vardi dichotomy conjecture shows that every such problem is solvable in polynomial time or $\Np$-complete.
For infinite $\mathcal{M}$, no such general classification can be expected, but $\Csp(\mathcal{N})$ has been studied extensively for \emph{first-order reducts} $\mathcal{N}$ of a fixed infinite structure $\mathcal{M}$.
A first-order reduct of $\mathcal{M}$ is a structure with the same domain whose relations are first-order definable in $\mathcal{M}$.
Complete complexity classifications, and in several cases full $\Ptime$ versus $\Np$-complete dichotomies, have been obtained for the CSPs of first-order reducts of the order of the rationals, the random graph, the random poset, and various other $\omega$-categorical structures.

Among these structures, first-order reducts of $(\mathbb{Z},<)$ occupy a distinguished position: $(\mathbb{Z},<)$ is not $\omega$-categorical, a property that in the other cases provides powerful tools, most notably a correspondence between primitive positive definability and preservation by polymorphisms, for classifying the complexity of the corresponding CSPs.
Essentially, the relations that first-order reducts of $(\mathbb{Z},<)$ can express beyond those already definable over $(\mathbb{Q},<)$ are relations that constrain the \emph{distance} between elements.
\begin{example}
    Let $\Diff_2\coloneq\{(a,a+2)\mid a\in\mathbb{Z}\}$.
    Then $\exists u\forall v(x<v<y\leftrightarrow u=v)$ defines $\Diff_2$ in $(\mathbb{Z},<)$ and hence $(\mathbb{Z},\Diff_2)$ is a first-order reduct of $(\mathbb{Z},<)$.
\end{example}
At the same time, every CSP of a first-order reduct of $(\mathbb{Q},<)$ equals $\Csp(\Gamma)$ for some first-order reduct $\Gamma$ of $(\mathbb{Z},<)$.
So the latter class strictly generalizes the former.
The complexity of $\Csp(\Gamma)$ for every finite-signature, first-order reduct $\Gamma$ of $(\mathbb{Z},<)$ was classified by \textcite{bodirsky2018}, building on an earlier classification, due to \textcite{bodirsky2016}, of the CSPs of \emph{locally finite}, finite-signature, first-order reducts of $(\mathbb{Z},\Succ)$, the integers equipped with the successor relation.
More recently, the classification of \textcite{bodirsky2018} has found application beyond its original setting: \textcite{dabrowski2026} use it as a key ingredient in the parametrized complexity analysis of disjunctive temporal problems.

In this paper we revisit the classification of \textcite{bodirsky2018} and make the following contributions.
We correct several errors in the original proof.
Most importantly, the argument in \cite{bodirsky2018} relies on the claim that every first-order reduct of $(\mathbb{Z},<)$ that contains the successor relation either primitively positively defines $<$ or is itself a first-order reduct of $(\mathbb{Z},\Succ)$.
We show that this claim is false, and we close the resulting gap in the classification.
We also show that, with the exception of those cases whose CSP equals the CSP of a finite structure, every tractable $\Csp(\Gamma)$ for a first-order reduct $\Gamma$ of $(\mathbb{Z},<)$ can be solved by a single, well-known uniform algorithm, which removes the need for the collection of ad-hoc algorithms designed separately for individual cases in \cite{bodirsky2018} and in \cite{bodirsky2016}.
Our main technical contribution, however, is a significant simplification and generalization of the technique underlying the reduction, in \cite{bodirsky2018}, of the full classification to the case of first-order reducts of $(\mathbb{Z},<)$ that contain the successor relation.
As a minor additional point, whereas the original proof relies on quantifier elimination in a way that only works for finite signatures, our proof does not rely on quantifier elimination in the same way and therefore also applies to reducts with infinite signatures.
Finally, we demonstrate the reach of this generalized technique by applying it to two further classes of structures: first-order reducts of the triangular lattice graph, and first-order reducts of infinite regular trees.
In both cases we obtain a full complexity classification, and in particular a $\Ptime$ versus $\Np$-complete dichotomy, for the corresponding CSPs.

\paragraph*{Tool and computational resource disclosure}
Large language models were used to assist with the writing and editing of this manuscript.
All mathematical content was developed without AI assistance.
GPT-6 Astra by OpenAI was used to review an earlier version of the manuscript and identified mathematical errors.
All identified errors were corrected without AI assistance.

\paragraph*{Notation}
$\mathbb{N}$ denotes the non-negative integers, including zero, and for $n\in\mathbb{N}$, $[n]\coloneq\{1,\dots,n\}$.
We write $\bar{a}$ for a tuple $(a_1,\dots,a_k)$; if the tuple itself carries an index, as in $\bar{a}_i$, we instead write its entries as $a_i^{(1)},\dots,a_i^{(k)}$.
The projection onto the $i$-th coordinate is written $\pi_i$, i.e., $\pi_i(\bar{a})\coloneq a_i$, and the projection of a relation $R$ of arity $k$ onto coordinates $I\subseteq[k]$ is written $\Proj_I(R)\coloneq\{(\pi_i(\bar{a}))_{i\in I}\mid\bar{a}\in R\}$.
Given an equivalence relation $\sim$ and an element $a$, $[a]_\sim$ denotes the $\sim$-equivalence class of $a$.
$\mathcal{A}$, $\mathcal{M}$, $\mathcal{N}$, and occasionally $\mathbb{M}$, denote structures with respective domains $A$, $M$, and $N$, while $\Gamma$ and $\Delta$ are reserved for first-order reducts of a fixed structure.
The $\sigma$-reduct of a $\tau$-structure $\mathcal{M}$, for $\sigma\subseteq\tau$, is written $\mathcal{M}\restriction\sigma$; the substructure induced by $A\subseteq M$ is written $\mathcal{M}|_A$, and the same bar notation restricts a function or relation to a subset of its domain.
We write $\mathcal{M}\preceq\mathcal{N}$ to mean that $\mathcal{M}$ is an elementary substructure of $\mathcal{N}$.
Let $\Hom(\mathcal{M},\mathcal{N})$ denote the set of homomorphisms from $\mathcal{M}$ to $\mathcal{N}$, $\End(\mathcal{M})\coloneq\Hom(\mathcal{M},\mathcal{M})$, and $\Aut(\mathcal{M})$ the automorphism group of $\mathcal{M}$.
When $\mathcal{M}$ and $\mathcal{N}$ are relational, we further let $\Pol(\mathcal{M},\mathcal{N})\coloneq\bigcup_{k\in\mathbb{N}}\Hom(\mathcal{M}^k,\mathcal{N})$.
For a binary relation symbol $R$, $n\in\mathbb{Z}$, and variables $x,y$, define
\[R^nxy\coloneq
\begin{cases}
    x=y & \text{if }n=0\\
    R^{-n}yx & \text{if }n<0\\
    \exists z_n(R^{n-1}xz_n\land Rz_ny) & \text{if }n>0.
\end{cases}\]
Finally, for $S,T\subseteq\mathbb{Z}$ and $\circ\in\{+,-,\cdot\}$, we write $S\circ T\coloneq\{s\circ t\mid s\in S,\,t\in T\}$, and $n\circ T\coloneq\{n\}\circ T$, $S\circ n\coloneq S\circ\{n\}$ for $n\in\mathbb{Z}$.

\section{Saturation}
A recurring theme in our proofs, already used by \textcite{bodirsky2018}, is that in order to study $\Csp(\Gamma)$ for a first-order reduct $\Gamma$ of $(\mathbb{Z},<)$, it is useful to look not at $\Gamma$ itself, but at a suitably enriched, \emph{saturated} version of it.
Recall that a countable structure is saturated if, informally, it realizes every consistent pattern of properties expressible with finitely many parameters, that is, every consistent type over a finite parameter set is realized in the structure.
Passing to a saturated elementary extension does not change which first-order sentences are true, and in particular preserves the CSP.

The countable saturated models of the theory of $(\mathbb{Z},<)$ are simple: they look like a dense linear order in which every point has been replaced by a copy of $(\mathbb{Z},<)$, or equivalently, $\mathbb{Q}\times\mathbb{Z}$ ordered lexicographically.
\begin{definition}
    For a first-order reduct $\Gamma$ of $(\mathbb{Z},<)$, let $\mathbb{Q}.\Gamma$ be the first-order reduct of $(\mathbb{Q}\times\mathbb{Z},<)$ over the same signature whose relations are defined by the same formulas.
\end{definition}
To work with such structures $\mathbb{Q}.\Gamma$ conveniently, we extend the difference operator on $\mathbb{Z}$ to $\mathbb{Q}\times\mathbb{Z}$: for $(p,a),(q,b)\in\mathbb{Q}\times\mathbb{Z}$, we define
\[(p,a)-(q,b)\coloneq
\begin{cases}
    a-b & \text{if }p=q\\
    \infty & \text{if }p>q\\
    -\infty & \text{if }p<q.
\end{cases}\]

What changes when passing from $\Gamma$ to $\mathbb{Q}.\Gamma$, and what makes this passage useful, is the set of endomorphisms.
Recall that a relation is primitive positive, or pp-, definable in a structure if it can be defined by an existentially quantified conjunction of atomic formulas over that structure, and that every pp-definable relation is preserved by all endomorphisms of the structure.
This matters for the complexity of $\Csp(\Gamma)$: if a relation is pp-definable in $\Gamma$, a constraint on it can be simulated using the relations of $\Gamma$ that are already present, so adding it to the signature does not change the complexity of the problem.
A saturated extension typically admits many endomorphisms that $\Gamma$ does not, and this richer supply of endomorphisms lets us rule out pp-definability for relations where the endomorphisms of $\Gamma$ alone would not suffice.
\begin{example}
    Let $\Gamma\coloneq(\mathbb{Z},\Diff_2,\lnot\Succ,\neq)$.
    Then $\End(\Gamma)=\{\mathbb{Z}\to\mathbb{Z},a\mapsto a+n\mid n\in\mathbb{Z}\}$, while $\mathbb{Q}.\Gamma$ has an endomorphism $e$ with $|e(a)-e(b)|=\infty$ for all $a,b\in\mathbb{Q}\times\mathbb{Z}$ with $a-b\notin 2\mathbb{Z}$.
    It follows that $\Gamma$ does not pp-define $\Succ$, a fact that is not visible just from the endomorphisms of $\Gamma$ itself.
\end{example}

Even more useful, and central to our proofs, is a partial converse of this fact.
By compactness, if every endomorphism of $\mathbb{Q}.\Gamma$ exhibits some fixed first-order definable behavior, then this behavior is already forced by a finite substructure, in the sense that the way endomorphisms act on finitely many elements already entails how they act everywhere else.
By encoding this finite witness as a primitive positive formula in the natural way, we can turn a purely algebraic hypothesis, namely that some property is preserved by all endomorphisms of $\mathbb{Q}.\Gamma$, into a pp-definability statement about $\Gamma$.

\subsection{Local and global}
The finite-witness phenomenon described above is, at heart, a compactness argument, so it is natural to make it precise in the framework of first-order logic.
The basic tool for reasoning about homomorphisms in this framework is the \emph{positive diagram} $\Diag^+(\mathcal{M})$ of a structure $\mathcal{M}$, defined as the set of atomic sentences valid in $(\mathcal{M},(a)_{a\in M})$.
Its usefulness stems from the easily verified fact that $(\mathcal{N},(b_a)_{a\in M})\models\Diag^+(\mathcal{M})$ if and only if $a\mapsto b_a$ defines a homomorphism $\mathcal{M}\to\mathcal{N}$.

The full strength of saturation is not needed for our arguments; a weaker property suffices.
Recall that a structure $\mathcal{M}$ is \emph{$\kappa$-universal}, for a cardinal $\kappa$, if every structure $\mathcal{N}\equiv\mathcal{M}$ with $|N|<\kappa$ embeds elementarily into $\mathcal{M}$.
Every saturated structure $\mathcal{M}$ is $|M|^+$-universal \cite[see][Theorem~10.1.6]{hodges1993}, so everything we prove for universal structures applies to saturated ones.

As a warm-up that already makes use of universality, we prove the following folklore result, which we will need later.
\begin{lemma}\label{lem:finite_csp_means_finite_range_endomorphism}
    Let $\mathcal{M}$ be an infinite relational $\kappa$-universal $\tau$-structure with $\aleph_0+|\tau|<\kappa$.
    Then $\mathcal{M}$ has a finite-range endomorphism if and only if there exists a finite structure $\mathcal{N}$ with $\Csp(\mathcal{M})=\Csp(\mathcal{N})$.
\end{lemma}
\begin{proof}
    Since the other direction is trivial, we assume that there is some finite $\tau$-structure $\mathcal{N}$ with $\Csp(\mathcal{M})=\Csp(\mathcal{N})$.
    Now, for all finite $A\subseteq M$ and $\sigma\subseteq\tau$ there exists a homomorphism $h:(\mathcal{M}|_A)\restriction\sigma\to\mathcal{N}\restriction\sigma$ and thus $\mathcal{N},(h(a))_{a\in A}\models\Diag^+((\mathcal{M}|_A)\restriction\sigma)$.
    By compactness it follows that $\Th(\mathcal{N})\cup\Diag^+(\mathcal{M})$ is consistent and hence has a model $\mathcal{N}'$.
    Since $\mathcal{N}$ is finite we have $\mathcal{N}\cong\mathcal{N}'\restriction\tau$.
    Moreover, $\mathcal{N}'\models\Diag^+(\mathcal{M})$ implies that there exists a homomorphism $f:\mathcal{M}\to\mathcal{N}$.

    Similarly, for all finite $\sigma\subseteq\tau$ there is a homomorphism $\mathcal{N}\restriction\sigma\to\mathcal{M}\restriction\sigma$ and by compactness the set $\Th(\mathcal{M})\cup\Diag^+(\mathcal{N})$ is consistent.
    Hence, by the Löwenheim-Skolem theorem \cite[see][Corollary~3.1.5]{hodges1993} it has a model $\mathcal{M}'$ of size less than $\kappa$, whose $\tau$-reduct admits a homomorphism from $\mathcal{N}$.
    Since $\mathcal{M}$ is $\kappa$-universal and $\mathcal{M}\equiv\mathcal{M}'\restriction\tau$, there exists an elementary embedding $\mathcal{M}'\restriction\tau\to\mathcal{M}$.
    Thus, we have a homomorphism $g:\mathcal{N}\to\mathcal{M}$ and $g\circ f$ is a finite-range endomorphism of $\mathcal{M}$.
\end{proof}
With the technique rehearsed, we now turn to Proposition~\ref{prop:compactness}, which we prove in a more general form than the preceding discussion requires.
First, rather than restricting attention to endomorphisms, we consider homomorphisms from an arbitrary structure to $\mathcal{M}$.
Second, we assume only universality instead of full saturation.
\begin{proposition}\label{prop:compactness}
    Let $\mathcal{M}$ be an infinite relational $\kappa$-universal $\tau$-structure with $\aleph_0+|\tau|<\kappa$ and $\mathcal{N}$ a structure with signature $\sigma\subseteq\tau$ and $|N|<\kappa$.
    Suppose $\Phi$ is a set of $\tau$-formulas in variables $\{x_a\mid a\in N\}$ such that
    \[\text{for all }h\in\Hom(\mathcal{N},\mathcal{M}\restriction\sigma)\text{ there is }\phi\in\Phi\text{ with }\mathcal{M}\models\phi((h(a))_{a\in N}).\]
    Then there are finite $A\subseteq N$, $\rho\subseteq\sigma$ and $\Psi\subseteq\Phi$ in variables $\{x_a\mid a\in A\}$ such that
    \[\text{for all }h\in\Hom(\mathcal{A},\mathcal{M}\restriction\rho)\text{ there is }\phi\in\Psi\text{ with }\mathcal{M}\models\phi((h(a))_{a\in A})\]
    where $\mathcal{A}$ is the $\rho$-reduct of $\mathcal{N}|_A$.
\end{proposition}
\begin{proof}
    Suppose for contradiction that the set
    \[\Th(\mathcal{M})\cup\Diag^+(\mathcal{N})\cup\{\lnot\phi((a)_{a\in N})\mid\phi\in\Phi\}\]
    is consistent.
    Then by the Löwenheim-Skolem theorem \cite[see][Corollary~3.1.5]{hodges1993} it has a model $\mathcal{B}$ of size less than $\kappa$.
    Since $\mathcal{B}\models\Diag^+(\mathcal{N})\cup\{\lnot\phi((a)_{a\in N})\mid\phi\in\Phi\}$, there is a homomorphism $f:\mathcal{N}\to\mathcal{B}\restriction\sigma$ with $\mathcal{B}\models\lnot\phi((f(a))_{a\in N})$ for all $\phi\in\Phi$.
    Because $\mathcal{B}\models\Th(\mathcal{M})$ and $\mathcal{M}$ is $\kappa$-universal there exists an elementary embedding $g:\mathcal{B}\restriction\tau\to\mathcal{M}$.
    However, now $h\coloneq g\circ f$ is a homomorphism $\mathcal{N}\to\mathcal{M}\restriction\sigma$ with $\mathcal{M}\models\lnot\phi((h(a))_{a\in N})$ for all $\phi\in\Phi$ contradicting the assumption.

    Hence, the set is inconsistent and by compactness there are finite $A\subseteq N$, $\rho\subseteq\sigma$ and $\Psi\subseteq\Phi$ in variables $\{x_a\mid a\in A\}$ such that for $\mathcal{A}\coloneq(\mathcal{N}|_A)\restriction\rho$ the subset
    \[\Th(\mathcal{M})\cup\Diag^+(\mathcal{A})\cup\{\lnot\phi((a)_{a\in A})\mid\phi\in\Psi\}\]
    is inconsistent and thus
    \[\Th(\mathcal{M})\cup\Diag^+(\mathcal{A})\vdash\bigvee_{\phi\in\Psi}\phi((a)_{a\in A}).\]
    Now, if $h:\mathcal{A}\to\mathcal{M}\restriction\rho$ is a homomorphism we have $\mathcal{M},(h(a))_{a\in A}\models\Diag^+(\mathcal{A})$, and hence there is some $\phi\in\Psi$ with $\mathcal{M}\models\phi((h(a))_{a\in A})$.
\end{proof}

\section{A general framework}\label{sec:model_theory}
Apart from saturation, $(\mathbb{Q}\times\mathbb{Z},<)$ has further properties that we exploit in our simplified proof of the reduction, mentioned in the introduction, of the classification of $\Csp(\Gamma)$, for first-order reducts $\Gamma$ of $(\mathbb{Z},<)$, to the case where $\Gamma$ contains the successor relation.
First, its automorphism group is \emph{transitive}, i.e., for all $a,b\in\mathbb{Q}\times\mathbb{Z}$ there is $f\in\Aut(\mathbb{Q}\times\mathbb{Z},<)$ with $f(a)=b$.
Second, every automorphism of a single copy of $(\mathbb{Z},<)$ inside $(\mathbb{Q}\times\mathbb{Z},<)$ can be extended to an automorphism of the full structure by leaving the elements of all other copies unchanged.
Third, the copy of $(\mathbb{Z},<)$ containing a given element of $\mathbb{Q}\times\mathbb{Z}$ coincides with the \emph{algebraic closure} of that element in $(\mathbb{Q}\times\mathbb{Z},<)$; moreover, this copy clearly embeds elementarily into $(\mathbb{Q}\times\mathbb{Z},<)$.
Formally, for a $\tau$-structure $\mathcal{M}$ and $A\subseteq M$, the algebraic closure $\Acl_\mathcal{M}(A)$ consists of all $a\in M$ for which there is a $\tau$-formula $\phi(x,y_1,\dots,y_k)$ and $\bar{b}\in A^k$ such that $\mathcal{M}\models\phi(a,\bar{b})$ and the set $\{a'\in M\mid\mathcal{M}\models\phi(a',\bar{b})\}$ is finite.
However, in a saturated structure algebraic closure admits an easier characterization: if $\mathcal{M}$ is saturated and $|A|<|M|$, then $a\in\Acl_\mathcal{M}(A)$ if and only if the set $\{f(a)\mid f\in\Aut(\mathcal{M}),f|_A=\Id_A\}$ is finite \cite[see][Proposition~4.3.27]{marker2002}.
We introduce a notion for structures satisfying suitable abstractions of the three properties above.
\begin{definition}
    A relational, infinite, saturated structure $\mathbb{M}$ is \emph{nice} if
    \begin{gather}
        \Aut(\mathbb{M})\text{ is transitive},\tag{trans}\label{eq:nice_transitive}\\
        \mathbb{M}|_{\Acl(a)}\preceq\mathbb{M}\text{ for some (equivalently all) }a\in M,\tag{emb}\label{eq:nice_elementary_embedding}\\
        a\in\Acl(b)\implies b\in\Acl(a),\tag{sym}\label{eq:nice_symmetry}\\
        f\in\Aut(\mathbb{M}|_{\Acl(a)})\implies f\cup\Id_{M\setminus\Acl(a)}\in\Aut(\mathbb{M}).\tag{ext}\label{eq:nice_extension}
    \end{gather}
\end{definition}
For elements $a,b$ of a nice structure $\mathbb{M}$, we write $a\asymp b$ if $\Acl(a)=\Acl(b)$.
By \eqref{eq:nice_symmetry}, the $\asymp$-class of $a$ is exactly $\Acl(a)$, i.e., $\Acl(a)=[a]_\asymp$.
\begin{lemma}\label{lem:equivalence_of_finite_range_endomorphisms}
    Let $\mathcal{M}$ be a first-order reduct of a nice structure $\mathbb{M}$.
    Then the following are equivalent.
    \begin{enumerate}
        \item $\mathcal{M}|_{\Acl_\mathbb{M}(a)}$ has a finite-range endomorphism for some (equivalently all) $a\in M$.
        \item $\Csp(\mathcal{M})=\Csp(\mathcal{N})$ for a finite structure $\mathcal{N}$.
        \item $\mathcal{M}$ has a finite-range endomorphism.
    \end{enumerate}
\end{lemma}
\begin{proof}
    If $e\in\End(\mathcal{M}|_{\Acl_\mathbb{M}(a)})$ has finite range, then $\mathcal{M}|_{e(\Acl_\mathbb{M}(a))}$ is finite with the same CSP as $\mathcal{M}|_{\Acl_\mathbb{M}(a)}$.
    By \eqref{eq:nice_elementary_embedding} we have $\mathcal{M}\equiv\mathcal{M}|_{\Acl_\mathbb{M}(a)}$ and hence $\Csp(\mathcal{M})=\Csp(\mathcal{M}|_{e(\Acl_\mathbb{M}(a))})$.

    If $\Csp(\mathcal{M})=\Csp(\mathcal{N})$ for a finite structure $\mathcal{N}$, then $\mathcal{M}$ has a finite-range endomorphism by Lemma~\ref{lem:finite_csp_means_finite_range_endomorphism}.

    Suppose $e\in\End(\mathcal{M})$ such that $e(M)\eqcolon A$ is finite and has minimal size, i.e., $|e(M)|\leq|f(M)|$ for all $f\in\End(\mathcal{M})$.
    Then $e(A)=A$, which means that $e|_A$ is a permutation of $A$.
    Since $e^{|A|!}|_A=\Id_A$ let us assume w.l.o.g.\ that already $e|_A=\Id_A$.
    Fix any $a\in A$ and suppose for contradiction that there exists $b'\in\Acl_\mathbb{M}(a)$ with $b\coloneq e(b')\notin\Acl_\mathbb{M}(a)$.
    By \eqref{eq:nice_extension} and \eqref{eq:nice_transitive} there is $f\in\Aut(\mathbb{M})$ with $f(a)=b'$ and $f(c)=c$ for all $c\in M\setminus\Acl_\mathbb{M}(a)$.
    Now, for $g\coloneq e\circ f$ the restriction $g|_A$ is a permutation of $A$ with $g|_A(a)=b=g|_A(b)$.
    Hence, $g(A)\subsetneq A$, which means that $(g\circ e)(M)\subsetneq A$ contradicts our choice of $e$.
    It follows that $e(\Acl_\mathbb{M}(a))\subseteq\Acl_\mathbb{M}(a)$ and thus $e|_{\Acl_\mathbb{M}(a)}$ is a finite-range endomorphism of $\mathcal{M}|_{\Acl_\mathbb{M}(a)}$.
\end{proof}
We also fix some general notation used throughout.
For a set $A$, $F\subseteq A^A$, and $\bar{a}\in A^k$, we let $\bar{a}^F\coloneq\{f(\bar{a})\mid f\in F\}$.
Typically, $A$ is the domain of a structure $\mathcal{A}$ and $F$ is $\Aut(\mathcal{A})$ or $\End(\mathcal{A})$; in the former case we call $\bar{a}^F$ the \emph{orbit} of $\bar{a}$ in $\mathcal{A}$.
An orbit of a pair is called an \emph{orbital}.
\begin{lemma}\label{lem:isolated_orbits}
    Let $\mathbb{M}$ be nice.
    If $a_1\asymp\dots\asymp a_k$, then $\bar{a}^{\Aut(\mathbb{M})}$ is first-order definable in $\mathbb{M}$.
\end{lemma}
\begin{proof}
    If $k=1$ we are done because $\Aut(\mathbb{M})$ is transitive and hence $a_1^{\Aut(\mathbb{M})}$ is defined by the formula $x=x$.
    Suppose $k>1$ and $(a_1,\dots,a_{k-1})^{\Aut(\mathbb{M})}$ is defined by a formula $\phi(\bar{x})$.
    Let $A\subseteq M$ be minimal such that it contains $a_k$ and is definable in $\mathbb{M}$ by a formula $\psi(a_1,\dots,a_{k-1},x)$.
    Since $a_k\in\Acl_\mathbb{M}(a_1,\dots,a_{k-1})$ we know that such a set exists and that it is finite.
    Clearly, $a_k^{\Aut(\mathbb{M},a_1,\dots,a_{k-1})}\subseteq A$.
    Now, suppose for contradiction that there exists $b\in A\setminus a_k^{\Aut(\mathbb{M},a_1,\dots,a_{k-1})}$.
    Since $\mathbb{M}$ is saturated there is a formula $\chi(a_1,\dots,a_{k-1},x)$ such that $\mathbb{M}\models\chi(a_1,\dots,a_{k-1},a_k)$ but $\mathbb{M}\not\models\chi(a_1,\dots,a_{k-1},b)$ \cite[see][Propositions~4.2.13 and~4.3.3]{marker2002}.
    However, now the formula $\psi(a_1,\dots,a_{k-1},x)\land\chi(a_1,\dots,a_{k-1},x)$ defines a proper subset of $A$ that contains $a_k$, which contradicts the minimality of $A$.
    Therefore, $\psi(a_1,\dots,a_{k-1},x)$ defines $a_k^{\Aut(\mathbb{M},a_1,\dots,a_{k-1})}$.
    It follows that $\phi(x_1,\dots,x_{k-1})\land\psi(x_1,\dots,x_k)$ defines $(a_1,\dots,a_k)^{\Aut(\mathbb{M})}$.
\end{proof}
\begin{definition}\label{def:rank_relation}
    For a first-order reduct $\mathcal{M}$ of a nice structure $\mathbb{M}$ we define ${\sim_\mathcal{M}}\subseteq M^2$ such that
    \[a\sim_\mathcal{M}b\iff e(a)\asymp e(b)\text{ for all }e\in\End(\mathcal{M}).\]
\end{definition}
Note, that $\sim_\mathcal{M}$ is an equivalence relation and preserved by $\End(\mathcal{M})$.
\begin{lemma}\label{lem:pp_define_image_under_endomorphisms}
    Let $\mathcal{M}$ be a first-order reduct of a nice structure $\mathbb{M}$.
    If $a_1\sim_\mathcal{M}\dots\sim_\mathcal{M}a_k$, then $\mathcal{M}$ pp-defines $\bar{a}^{\End(\mathcal{M})}$, which is a finite union of $k$-orbits of $\Aut(\mathbb{M})$.
\end{lemma}
\begin{proof}
    Suppose w.l.o.g.\ that $\mathcal{M}$ is a $\sigma$-reduct of $\mathbb{M}$.
    For any $e\in\End(\mathcal{M})$ we have $e(a_1)\asymp\dots\asymp e(a_k)$ and thus by Lemma~\ref{lem:isolated_orbits} a formula $\phi_e((x_a)_{a\in\bar{a}})$ that defines $e(\bar{a})^{\Aut(\mathbb{M})}$ in $\mathbb{M}$.
    Let $\Phi\coloneq\{\phi_e\mid e\in\End(\mathcal{M})\}$.
    Then by Proposition~\ref{prop:compactness} there are finite $A\subseteq M$, $\rho\subseteq\sigma$ and $\Psi\subseteq\Phi$ such that for all $h\in\Hom(\mathcal{A},\mathcal{M}\restriction\rho)$ there is $\phi\in\Psi$ with $\mathbb{M}\models\phi((h(a))_{a\in\bar{a}})$, where $\mathcal{A}$ is the $\rho$-reduct of $\mathcal{M}|_A$.
    Let $\psi((x_a)_{a\in\bar{a}})$ be the pp-formula obtained by existentially quantifying all $x_b$ with $b\notin\bar{a}$ in $\bigwedge\{Rx_{b_1}\dots x_{b_\ell}\mid R\in\rho,(b_1,\dots,b_\ell)\in R^\mathcal{A}\}$.
    Now, $\mathbb{M}\models\psi((b_a)_{a\in\bar{a}})$ iff there is $h\in\Hom(\mathcal{A},\mathcal{M}\restriction\rho)$ with $h(a)=b_a$ for all $a\in\bar{a}$.
    Hence, if $\mathbb{M}\models\psi((b_a)_{a\in\bar{a}})$ then $\mathbb{M}\models\phi((b_a)_{a\in\bar{a}})$ for some $\phi\in\Psi$ and thus $(b_a)_{a\in\bar{a}}\in\bar{a}^{\End(\mathcal{M})}$.
    Because every $\phi\in\Psi$ defines a $k$-orbit of $\Aut(\mathbb{M})$ it follows that the relation defined by $\psi$ is a finite union of $k$-orbits.
    On the other hand, every $e\in\End(\mathcal{M})$ induces a homomorphism $\mathcal{A}\to\mathcal{M}\restriction\rho$ and thus we have $\mathbb{M}\models\psi(e(\bar{a}))$ for all $e\in\End(\mathcal{M})$.
\end{proof}
\begin{proposition}\label{prop:endomorphism_with_image_in_relation}
    Let $\mathcal{M}$ be a first-order reduct of a nice structure $\mathbb{M}$ and let ${\approx}\subseteq{\sim_\mathcal{M}}$ be an equivalence relation that is preserved by $\End(\mathcal{M})$ such that for all $a\sim_\mathcal{M}b$ there is $e\in\End(\mathcal{M})$ with $e(a)\approx e(b)$.
    Then there exists $e\in\End(\mathcal{M})$ such that
    \begin{enumerate}
        \item $e(M)$ is finite or
        \item\label{item:nice_structures_proposition_infinite} for all $a,b\in M$ either $e(a)\approx e(b)$ or $e(a)\not\asymp e(b)$.
    \end{enumerate}
\end{proposition}
\begin{proof}
    First, we claim that whenever $B\subseteq M$ and for all $e\in\End(\mathcal{M})$ there are $a,b\in B$ with $e(a)\not\approx e(b)$ and $e(a)\asymp e(b)$, then there exists a finite $A\subseteq B$ such that
    \begin{multline}\tag{$*$}\label{eq:endomorphism_with_image_in_relation}
        \text{there is a finite union }U\subseteq{\asymp}\text{ of orbitals of }{\Aut(\mathbb{M})}\text{ so that}\\
        \text{for all }e\in\End(\mathcal{M})\text{ there are }a,b\in e(A)\text{ with }a\not\approx b\text{ and }(a,b)\in U.
    \end{multline}
    For every $e\in\End(\mathcal{M})$ pick some $a,b\in B$ with $e(a)\not\approx e(b)$ and $e(a)\asymp e(b)$ and let $\phi_e(x_a,x_b)$ be a formula defining $(e(a),e(b))^{\Aut(\mathbb{M})}$, which exists by Lemma~\ref{lem:isolated_orbits}.
    Then by Proposition~\ref{prop:compactness}, using the set $\Phi\coloneq\{\phi_e\mid e\in\End(\mathcal{M})\}$), there are finite $A\subseteq B$ and $\Psi\subseteq\Phi$ such that for all $e\in\End(\mathcal{M})$ there are $\phi\in\Psi$ and $a,b\in A$ with $\mathbb{M}\models\phi(e(a),e(b))$.
    Since each $\phi\in\Psi$ defines an orbital of $\Aut(\mathbb{M})$, the set $U\coloneq\{(a,b)\in M^2\mid\phi\in\Psi,\mathbb{M}\models\phi(a,b)\}$ is a finite union of orbitals.
    Because $\approx$ is preserved by $\Aut(\mathbb{M})\subseteq\End(\mathcal{M})$ we have $U\cap{\approx}=\emptyset$, which proves the claim.

    Suppose we are not in Case~\ref{item:nice_structures_proposition_infinite} and let $A\subseteq M$ be finite of minimal size satisfying \eqref{eq:endomorphism_with_image_in_relation}.
    If $A=\{a,b\}$ then for all $e\in\End(\mathcal{M})$ we have $e(a)\asymp e(b)$ and thus $a\sim_\mathcal{M}b$.
    Hence, by assumption there exists some $f\in\End(\mathcal{M})$ such that $f(a)\approx f(b)$, which contradicts $A$ satisfying \eqref{eq:endomorphism_with_image_in_relation}.
    It follows that $A$ contains at least three elements.

    Now, suppose for contradiction that there are $e\in\End(\mathcal{M})$ and distinct $a,b\in e(A)$ with $a\approx b$.
    Clearly, $e(A)$ satisfies \eqref{eq:endomorphism_with_image_in_relation}.
    By our minimality assumption on the size of $A$ the set $e(A)\setminus\{a\}$ does not satisfy \eqref{eq:endomorphism_with_image_in_relation}.
    From what we proved earlier it follows that there exists some $f\in\End(\mathcal{M})$ such that for all $c,d\in e(A)\setminus\{a\}$ we have $f(c)\approx f(d)$ or $f(c)\not\asymp f(d)$.
    Because $e(A)$ satisfies \eqref{eq:endomorphism_with_image_in_relation} there is some $c\in e(A)\setminus\{a\}$ such that $f(a)\not\approx f(c)$ and $f(a)\asymp f(c)$.
    Since $\approx$ is preserved by $\End(\mathcal{M})$ we have $f(a)\approx f(b)$ and in particular $f(a)\asymp f(b)$.
    It follows that $f(b)\not\approx f(c)$ and $f(b)\asymp f(c)$ contradicting the choice of $f$.

    Therefore, for every $e\in\End(\mathcal{M})$ and all distinct $a,b\in A$ we have $e(a)\not\approx e(b)$.
    Now, for all $e\in\End(\mathcal{M})$ and $B\subsetneq e(A)$, since $B$ does not satisfy \eqref{eq:endomorphism_with_image_in_relation}, there is some $f_B\in\End(\mathcal{M})$ with $f_B(a)\not\asymp f_B(b)$ for all distinct $a,b\in B$.
    Because $A$ contains at least three elements it follows that there are $e\in\End(\mathcal{M})$ and $a\in e(A)$ such that $a\not\asymp b$ for all $b\in e(A)\setminus\{a\}\eqcolon B$.
    Fix any $a'\in\Acl_\mathbb{M}(a)$.
    By \eqref{eq:nice_transitive} and \eqref{eq:nice_extension}, there exists $g\in\Aut(\mathbb{M})$ with $g(a)=a'$ and $g(b)=b$ for all $b\in B$.
    Now, as $e(A)$ satisfies \eqref{eq:endomorphism_with_image_in_relation}, there is some $b\in B$ such that $\big(f_B(a'),f_B(b)\big)=\big((f_B\circ g)(a),(f_B\circ g)(b)\big)\in U$.
    Since we picked $a'\in\Acl_\mathbb{M}(a)$ arbitrarily we get that $f_B(\Acl_\mathbb{M}(a))\subseteq\{c\mid b\in f_B(B),(c,b)\in U\}$.
    It follows that $f_B(\Acl_\mathbb{M}(a))$ is finite.
    Now, $f_B|_{\Acl_\mathbb{M}(a)}$ is a homomorphism $\mathcal{M}|_{\Acl_\mathbb{M}(a)}\to\mathcal{M}$ with finite range.
    By \eqref{eq:nice_elementary_embedding} we have $\mathcal{M}\equiv\mathcal{M}|_{\Acl_\mathbb{M}(a)}$ and thus by Lemma~\ref{lem:finite_csp_means_finite_range_endomorphism} there exists an endomorphism of $\mathcal{M}$ with finite range.
\end{proof}
\begin{lemma}\label{lem:avoiding_orbitals}
    Let $\mathcal{M}$ be a first-order reduct of a nice structure $\mathbb{M}$ and ${\approx}\subseteq{\sim_\Gamma}$ an equivalence relation that is preserved by $\Aut(\mathbb{M})$ and such that ${\approx}\cap\Acl_\mathbb{M}(a)^2$ has finitely many equivalence classes for all $a\in M$.
    Then the following are equivalent.
    \begin{enumerate}
        \item\label{item:avoiding_orbitals_finite} $\mathcal{M}$ has a finite-range endomorphism.
        \item\label{item:avoiding_orbitals_union} There exists a finite union $U$ of orbitals of $\Aut(\mathbb{M})$ such that for all $a\approx b$ there exists $e\in\End(\mathcal{M})$ with $(e(a),e(b))\in U$.
    \end{enumerate}
\end{lemma}
\begin{proof}
    Since the other direction is trivial, we prove that Item~\ref{item:avoiding_orbitals_union} implies Item~\ref{item:avoiding_orbitals_finite}.
    Because ${\sim_\Gamma}\supseteq{\approx}$ is preserved by $\End(\mathcal{M})$, we can assume w.l.o.g.\ that $U\subseteq{\sim_\Gamma}$.
    Let $V\coloneq U^{\End(\mathcal{M})}$.
    Then by Lemma~\ref{lem:pp_define_image_under_endomorphisms} the set $V$ is a finite union of orbitals of $\Aut(\mathbb{M})$.

    Suppose for contradiction that for all $e\in\End(\mathcal{M})$ there are $a,b\in M$ such that $e(a)\approx e(b)$ and $(e(a),e(b))\notin V$.
    For each $e\in\End(\mathcal{M})$ we fix such $a,b\in M$ and let $\phi_e(x_a,x_b)$ be a formula defining the orbital $(e(a),e(b))^{\Aut(\mathbb{M})}$ in $\mathbb{M}$.
    Such formulas exist by Lemma~\ref{lem:isolated_orbits}.
    Setting $\Phi\coloneq\{\phi_e\mid e\in\End(\mathcal{M})\}$, by Proposition~\ref{prop:compactness} there exists a finite $A\subseteq M$ such that
    \begin{equation}\tag{\dag}\label{eq:avoiding_orbitals}
        \text{for all }e\in\End(\mathcal{M})\text{ there are }a,b\in A\text{ with }e(a)\approx e(b)\text{ and }(e(a),e(b))\notin V.
    \end{equation}
    For all finite $A\subseteq M$ we write
    \[W(A)\coloneq\{(a,b)\in A^2\mid\text{exists }e\in\End(\mathcal{M})\text{ with }e(a)\approx e(b)\text{ and }(e(a),e(b))\notin V\}\]
    and fix a finite $A\subseteq M$ satisfying \eqref{eq:avoiding_orbitals} such that $W(A)$ has minimal size.
    Now, pick any $(a,b)\in W(A)$ and let $e\in\End(\mathcal{M})$ with $e(a)\approx e(b)$ and $(e(a),e(b))\notin V$.
    By assumption there exists $f\in\End(\mathcal{M})$ such that for $g\coloneq f\circ e$ we have $(g(a),g(b))\in U$.
    Now, $g(A)$ satisfies \eqref{eq:avoiding_orbitals} with $W(g(A))\subseteq g(W(A))\setminus\{(g(a),g(b))\}$ contradicting our choice of $A$.

    Hence, there exists $e\in\End(\mathcal{M})$ such that for all $a,b\in M$ we have $e(a)\not\approx e(b)$ or $(e(a),e(b))\in V$.
    Fix $a\in M$ and let $B\subseteq\Acl_\mathbb{M}(a)$ be minimal such that $e(B)$ contains exactly one representative of each equivalence class of ${\approx}\cap e(\Acl_\mathbb{M}(a))$.
    Since ${\sim_\Gamma}\cap\Acl_\mathbb{M}(a)^2\supseteq{\approx}\cap\Acl_\mathbb{M}(a)^2$ has finitely many equivalence classes, there is a finite $C\subseteq M$ such that $e(B)\subseteq\bigcup_{c\in C}\Acl_\mathbb{M}(c)$.
    Because ${\approx}\cap\Acl_\mathbb{M}(c)^2$ has finitely many equivalence classes for each $c\in C$, it follows that $B$ is finite.
    Now, for each $c\in\Acl_\mathbb{M}(a)$ there exists $b\in B$ such that $(e(b),e(c))\in V$.
    Since $V\subseteq{\asymp}$ is a finite union of orbitals of $\Aut(\mathbb{M})$ it follows that $e(\Acl_\mathbb{M}(a))\subseteq\{c\mid b\in B,(e(b),c)\in V\}$ is finite.
    Thus, $e|_{\Acl_\mathbb{M}(a)}:\mathcal{M}|_{\Acl_\mathbb{M}(a)}\to\mathcal{M}$ is a homomorphism with finite range and by Lemma~\ref{lem:finite_csp_means_finite_range_endomorphism} with \eqref{eq:nice_elementary_embedding} it follows that $\mathcal{M}$ has a finite-range endomorphism.
\end{proof}

\section{PP-Defining successor}\label{sec:defining_successor}
We now apply the results from Section~\ref{sec:model_theory} to first-order reducts of $(\mathbb{Z},<)$ in order to prove Theorem~\ref{thm:defining_successor}, which shows that every such reduct falls into one of four categories.
For two of these categories, the complexity of the corresponding CSP is already fully classified, and for a third we can establish $\Np$-hardness by a simple argument (Lemma~\ref{lem:hardness_from_distances}).
This leaves only the case in which the reduct pp-defines the successor relation, so for the remainder of the paper we may restrict attention to reducts containing $\Succ$.

This was already shown by \textcite{bodirsky2018}, but their proof is long, spanning pages~10 to~27, and only applies to reducts with finite signature, owing to the way it relies on quantifier elimination.
In contrast, our proof amounts to applying the results of Section~\ref{sec:model_theory} to an appropriately chosen relation (Definition~\ref{def:tight_rank_relation}); the main task of this section is to verify that this relation satisfies the required assumptions in the case relevant to our classification.

Two basic facts about $(\mathbb{Q}\times\mathbb{Z},<)$ will be useful throughout: for all $a,b\in\mathbb{Q}\times\mathbb{Z}$ the orbital $(a,b)^{\Aut(\mathbb{Q}\times\mathbb{Z},<)}$ is $\{(c,d)\in(\mathbb{Q}\times\mathbb{Z})^2\mid a-b=c-d\}$, and it is not hard to see that $(\mathbb{Q}\times\mathbb{Z},<)$ is nice, with $a\asymp b\iff|a-b|<\infty$.

It is also convenient to introduce some notation for working with first-order reducts of $(\mathbb{Z},<)$ and their saturated counterparts.
Given $S\subseteq\mathbb{Z}$, let $\Diff_S\coloneq\{(a,b)\in\mathbb{Z}^2\mid b-a\in S\}$ and $\Dist_S\coloneq\Diff_S\cup\Diff_{-S}$, and use the same notation for the corresponding relations on $\mathbb{Q}\times\mathbb{Z}$ whenever $S$ is finite, with the intended domain always clear from context.
Singletons are abbreviated by $\Diff_n\coloneq\Diff_{\{n\}}$ and $\Dist_n\coloneq\Dist_{\{n\}}$ for $n\in\mathbb{Z}$, and we write $(p,a)+n\coloneq(p,a+n)$ for $(p,a)\in\mathbb{Q}\times\mathbb{Z}$ and $n\in\mathbb{Z}$.
\begin{lemma}\label{lem:equivalence_relations_are_periodic}
    If ${\sim}\subseteq(\mathbb{Q}\times\mathbb{Z})^2$ is a non-total equivalence relation that is preserved by $\Aut(\mathbb{Q}\times\mathbb{Z},<)$ then there exists $n\in\mathbb{N}$ such that $a\sim b\iff a-b\in n\mathbb{Z}$.
\end{lemma}
\begin{proof}
    If there are $a\sim b$ with $|a-b|=\infty$, it follows from \eqref{eq:nice_transitive}, \eqref{eq:nice_elementary_embedding} and \eqref{eq:nice_extension} that $\sim$ is total against the assumption.
    Since $\sim$ is preserved by $\Aut(\mathbb{Q}\times\mathbb{Z},<)$, it is a union of orbitals and thus $a\sim b\iff a-b\in S$ for some non-empty $S\subseteq\mathbb{Z}$.
    If $S=\{0\}$ we are done setting $n\coloneq 0$.
    Otherwise, let $n\coloneq\gcd S$.
    By B\'{e}zout's identity there are $s_1,\dots,s_k\in S$ and $m_1,\dots,m_k\in\mathbb{Z}$ such that $m_1s_1+\dots+m_ks_k=n$.
    Hence, for any $a\in\mathbb{Q}\times\mathbb{Z}$ there is a $\sim$-path from $a$ to $a+n$.
    It follows that $n\in S$ and thus $a\sim b\iff a-b\in n\mathbb{Z}$.
\end{proof}
To simplify notation, for a first-order reduct $\Gamma$ of $(\mathbb{Z},<)$ we will write $\sim_\Gamma$ instead of $\sim_{\mathbb{Q}.\Gamma}$ (see Definition~\ref{def:rank_relation}).
\begin{definition}\label{def:tight_rank_relation}
    For a first-order reduct $\Gamma$ of $(\mathbb{Z},<)$ let ${\approx_\Gamma}\subseteq(\mathbb{Q}\times\mathbb{Z})^2$ such that
    \[a\approx_\Gamma b\iff|e(a)-e(b)|\leq|a-b|<\infty\text{ for all }e\in\End(\mathbb{Q}.\Gamma).\]
\end{definition}
For Proposition~\ref{prop:endomorphism_with_image_in_relation} to apply, we would like $\approx_\Gamma$ to be an equivalence relation that is preserved by $\End(\mathbb{Q}.\Gamma)$.
In general, however, neither property holds, as the following example shows.
\begin{example}
    Let $\Gamma\coloneq(\mathbb{Z},\Dist_{\{0,1,2\}},\Dist_{\{0,1,2,3\}})$.
    Then $(0,0)\approx_\Gamma(0,3)\approx_\Gamma(0,1)$.
    However,
    \[e:\mathbb{Q}.\Gamma\to\mathbb{Q}.\Gamma,(p,a)\mapsto
    \begin{cases}
        (p,2) & \text{if }a=1\\
        (p,0) & \text{if }a<1\\
        (p,a) & \text{if }a>1
    \end{cases}\]
    is an endomorphism of $\mathbb{Q}.\Gamma$ with $|e(0,0)-e(0,1)|=2>1=|(0,0)-(0,1)|$, so $(0,0)\not\approx_\Gamma(0,1)$; hence $\approx_\Gamma$ is not transitive.
    Moreover, $e(0,1)=(0,2)\not\approx_\Gamma(0,3)=e(0,3)$, so $\approx_\Gamma$ is not preserved by $e$.
\end{example}
Fortunately, in the case relevant to our classification, namely, whenever $\Gamma$ has no finite-range endomorphism, $\approx_\Gamma$ is both an equivalence relation and preserved by $\End(\mathbb{Q}.\Gamma)$, as we show below.
The proof of the next lemma is partly based on that of \cite[Lemma~5.23]{bodirsky2018}.
\begin{lemma}\label{lem:pp_define_distances}
    Let $\Gamma$ be a first-order reduct of $(\mathbb{Z},<)$ without a finite-range endomorphism.
    If $a\approx_\Gamma b$ then $\mathbb{Q}.\Gamma$ pp-defines $\Diff_{b-a}$ or $\Dist_{b-a}$.
\end{lemma}
\begin{proof}
    By Lemma~\ref{lem:equivalence_relations_are_periodic} there exists $r\in\mathbb{N}$ such that $c\sim_\Gamma d\iff c-d\in r\mathbb{Z}$ for all $c,d\in\mathbb{Q}\times\mathbb{Z}$.
    If $r=0$, since ${\approx_\Gamma}\subseteq{\sim_\Gamma}$ we have $a=b$ and $\Diff_{a-b}$ is pp-defined by $x=y$.
    Suppose $r>0$.
    By Lemma~\ref{lem:pp_define_image_under_endomorphisms} there exists $m\coloneq\max\{|e(c)-e(c+r)|\mid c\in\mathbb{Q}\times\mathbb{Z},e\in\End(\mathbb{Q}.\Gamma)\}\in r\mathbb{N}_{>0}$.
    From Lemma~\ref{lem:equivalence_of_finite_range_endomorphisms} it follows that $\mathbb{Q}.\Gamma$ does not have a finite-range endomorphism and hence by Lemma~\ref{lem:avoiding_orbitals} there is some $n\in r\mathbb{N}$ such that $|e(c)-e(c+n)|>m$ for all $e\in\End(\mathbb{Q}.\Gamma)$ and $c\in\mathbb{Q}\times\mathbb{Z}$.
    Taking $n\in r\mathbb{N}$ minimal with this property yields $|e(c)-e(c+n)|\geq n>m$ for all $e\in\End(\mathbb{Q}.\Gamma)$ and $c\in\mathbb{Q}\times\mathbb{Z}$.
    Now, suppose for contradiction that there are $e\in\End(\mathbb{Q}.\Gamma)$ and $c,d\in\mathbb{Q}\times\mathbb{Z}$ with $c-d\in r\mathbb{Z}$ such that $e(c)-e(c+n)\geq n\leq e(d+n)-e(d)$.
    We assume w.l.o.g.\ that $|c-d|=r$.
    Then
    \begin{align*}
        e(c)-e(c+n)&\geq n\\
        e(d+n)-e(d)&\geq n\\
        e(d)-e(c)&\geq -m\\
        e(c+n)-e(d+n)&\geq -m.
    \end{align*}
    Summing over these inequalities yields $0\geq 2n-2m$ contradicting $n>m$.
    Hence, for all $\ell\in\mathbb{Z}$, $e\in\End(\mathbb{Q}.\Gamma)$ and $c\in\mathbb{Q}\times\mathbb{Z}$ we have $|e(c)-e(c+\ell n)|\geq|\ell n|$.
    Choosing $\ell\coloneq b-a$, by Lemma~\ref{lem:pp_define_image_under_endomorphisms} $\mathbb{Q}.\Gamma$ pp-defines $Q\coloneq(a,a+n(b-a))^{\End(\mathbb{Q}.\Gamma)}$.
    By Lemma~\ref{lem:pp_define_image_under_endomorphisms} $\mathbb{Q}.\Gamma$ also pp-defines $R\coloneq(a,b)^{\End(\mathbb{Q}.\Gamma)}$.
    It is now easy to see that $\exists z(Rxy\land R^{n-1}yz\land Qxz)$ pp-defines $\Diff_{b-a}$ or $\Dist_{b-a}$.
    \begin{center}
        \begin{tikzcd}[sep=huge]
            x \arrow[r, dash, "\leq|a-b|", "R" swap] \arrow[rr, dash, bend right, "Q", "\geq n|a-b|" swap] & y \arrow[r, dash, "\leq(n-1)|a-b|", "R^{k-1}" swap] & \exists z
        \end{tikzcd}\
    \end{center}
\end{proof}
\begin{lemma}\label{lem:tight_rank_relation_is_equivalence}
    Let $\Gamma$ be a first-order reduct of $(\mathbb{Z},<)$ without a finite-range endomorphism.
    Then $\approx_\Gamma$ is an equivalence relation.
\end{lemma}
\begin{proof}
    It is clear that $\approx_\Gamma$ is reflexive and symmetric.
    Moreover, if $a\approx_\Gamma b\approx_\Gamma c$ with $a\leq b\leq c$ then also $a\approx_\Gamma c$.
    It remains to show that $a\approx_\Gamma b\approx_\Gamma c$ with $a<c<b$ implies $a\approx_\Gamma c$.
    Since $\approx_\Gamma$ is preserved by $\Aut(\mathbb{Q}.\Gamma)$, for $d\coloneq a-(b-c)$ we have $d\approx_\Gamma a$ and $d\approx_\Gamma c$.
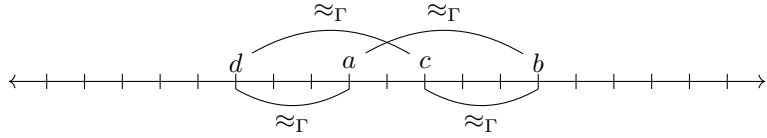
    \begin{center}
        \begin{tikzpicture}
            \draw[<->] (-5,0) -- (5,0);
            \foreach \x in {-4.5,-4,-3.5,-3,-2.5,-2,-1.5,-1,-0.5,0,0.5,1,1.5,2,2.5,3,3.5,4,4.5} {
                \draw (\x,0.1) -- (\x,-0.1);
            }
            \node (a) at (-0.5,0.25) {$a$};
            \node (b) at (2,0.25) {$b$};
            \node (c) at (0.5,0.25) {$c$};
            \node (d) at (-2,0.25) {$d$};
            \draw[bend left] (a) to node[above] {$\approx_\Gamma$} (b);
            \draw[bend right] (0.5,-0.1) to node[below] {$\approx_\Gamma$} (2,-0.1);
            \draw[bend right] (c) to node[above] {$\approx_\Gamma$} (d);
            \draw[bend left] (-0.5,-0.1) to node[below] {$\approx_\Gamma$} (-2,-0.1);
        \end{tikzpicture}
    \end{center}
    Hence, $b\approx_\Gamma d$ and by Lemma~\ref{lem:pp_define_distances} we get $|e(b)-e(d)|=|b-d|$ for all $e\in\End(\mathbb{Q}.\Gamma)$.
    It follows that $|e(a)-e(c)|=|a-c|$ for all $e\in\End(\mathbb{Q}.\Gamma)$ and thus $a\approx_\Gamma c$.
\end{proof}
\begin{fact}\label{fact:quantifier_elimination_order}
    The structures $(\mathbb{Z},<)$ and $(\mathbb{Z},(x+n<y)_{n\in\mathbb{Z}})$ are interdefinable and $\Th(\mathbb{Z},(x+n<y)_{n\in\mathbb{Z}})$ has quantifier elimination.
\end{fact}
\begin{lemma}[{\cite[see][Lemma~5.28]{bodirsky2018}}]\label{lem:substructure_isomorphic_to_fo_reduct}
    Let $\Gamma$ be a first-order reduct of $(\mathbb{Z},<)$ and $n\in\mathbb{N}_{>0}$.
    Then there exists a first-order reduct $\Delta$ of $(\mathbb{Z},<)$ such that $\Gamma|_{n\mathbb{Z}}\to\Delta,a\mapsto a/n$ and $\mathbb{Q}.\Gamma|_{\mathbb{Q}\times n\mathbb{Z}}\to\mathbb{Q}.\Delta,(p,a)\mapsto(p,a/n)$ are isomorphisms.
\end{lemma}
\begin{proof}
    By Fact~\ref{fact:quantifier_elimination_order}, for every $R$ in the signature $\tau$ of $\Gamma$ there exists a quantifier-free formula $\phi_R(\bar{x})$ in the language $\sigma\coloneq\{x+m<y\mid m\in\mathbb{Z}\}$ such that $(\mathbb{Z},\sigma)\models\phi_R(\bar{a})\iff\bar{a}\in R^\Gamma$.
    For each $R\in\tau$ let $\psi_R$ be the formula $\phi_R$ with every atomic formula of the form $x+m<y$ replaced by $x+\lfloor m/n\rfloor<y$.
    We pick $\Delta$ to be the first-order reduct of $(\mathbb{Z},<)$ with signature $\tau$ and $\bar{a}\in R^\Delta\coloniff(\mathbb{Z},\sigma)\models\psi_R(\bar{a})$ for all $R\in\tau$.
    Then for all $a,b\in\mathbb{Z}$ we get $na+m<nb\iff a+\lfloor m/n\rfloor<b$ and thus
    \begin{align*}
        (na_1,\dots,na_k)\in R^\Gamma&\iff(\mathbb{Z},\sigma)\models\phi_R(na_1,\dots,na_k)\\
        &\iff(\mathbb{Z},\sigma)\models\psi_R(a_1,\dots,a_k)\\
        &\iff(a_1,\dots,a_k)\in R^\Delta
    \end{align*}
    for all $R\in\tau$ (of arity $k$) and $\bar{a}\in\mathbb{Z}^k$.

    The proof for $\mathbb{Q}.\Gamma$ is analogous.
\end{proof}
\begin{theorem}\label{thm:defining_successor}
    Let $\Gamma$ be a first-order reduct of $(\mathbb{Z},<)$.
    Then $\mathbb{Q}.\Gamma$ is homomorphically equivalent to one of the following:
    \begin{enumerate}
        \item\label{item:defining_successor_finite} a finite structure;
        \item\label{item:defining_successor_rationals} a first-order reduct of $(\mathbb{Q},<)$;
        \item\label{item:defining_successor_succ} a first-order reduct of $(\mathbb{Q}\times\mathbb{Z},<)$ that pp-defines $\Succ$;
        \item\label{item:defining_successor_dist} a first-order reduct of $(\mathbb{Q}\times\mathbb{Z},<)$ that pp-defines $\Dist_n$ for every $n\in\mathbb{N}$.
    \end{enumerate}
\end{theorem}
\begin{proof}
    If $\Gamma$ has an endomorphism $e$ with finite range, then it has the same CSP as the finite structure $\Gamma|_{e(\mathbb{Z})}$ and we are in Case~\ref{item:defining_successor_finite} by Lemma~\ref{lem:finite_csp_means_finite_range_endomorphism}.
    Suppose from now on that $\Gamma$ does not have a finite-range endomorphism and note that by Lemma~\ref{lem:equivalence_of_finite_range_endomorphisms} neither does $\mathbb{Q}.\Gamma$ have a finite-range endomorphism.
    Then by Lemma~\ref{lem:tight_rank_relation_is_equivalence}, $\approx_\Gamma$ is an equivalence relation.
    Since ${\approx_\Gamma}\subseteq{\sim_\Gamma}$, by Lemma~\ref{lem:pp_define_image_under_endomorphisms} for all $a\sim_\Gamma b$ there exists $e\in\End(\mathbb{Q}.\Gamma)$, maximizing $|e(a)-e(b)|$, with $e(a)\approx_\Gamma e(b)$.
    Now, by Proposition~\ref{prop:endomorphism_with_image_in_relation} there exists $e\in\End(\mathbb{Q}.\Gamma)$ such that for all $a,b\in\mathbb{Q}\times\mathbb{Z}$ either $e(a)\approx_\Gamma e(b)$ or $|e(a)-e(b)|=\infty$.
    By Lemma~\ref{lem:equivalence_relations_are_periodic} there exists $n\in\mathbb{N}$ such that $a\approx_\Gamma b\iff a-b\in n\mathbb{Z}$.
    Composing $e$ with an appropriate automorphism of $(\mathbb{Q}\times\mathbb{Z},<)$ we get an endomorphism of $\mathbb{Q}.\Gamma$ whose image is contained in $\mathbb{Q}\times n\mathbb{Z}$.
    Let us assume that already $e(\mathbb{Q}\times\mathbb{Z})\subseteq\mathbb{Q}\times n\mathbb{Z}$.
    Then $\mathbb{Q}.\Gamma$ is homomorphically equivalent to $\mathbb{Q}.\Gamma|_{\mathbb{Q}\times n\mathbb{Z}}$.
    If $n=0$, it is easy to see that $\mathbb{Q}.\Gamma|_{\mathbb{Q}\times n\mathbb{Z}}$ is isomorphic to a first-order reduct of $(\mathbb{Q},<)$ and we are in Case~\ref{item:defining_successor_rationals}.
    Suppose that $n>0$.
    Then by Lemma~\ref{lem:substructure_isomorphic_to_fo_reduct} there exists a first-order reduct $\Delta$ of $(\mathbb{Z},<)$ such that $\iota:\mathbb{Q}.\Gamma|_{\mathbb{Q}\times n\mathbb{Z}}\to\mathbb{Q}.\Delta,(p,a)\mapsto(p,a/n)$ is an isomorphism.

    We claim that $a\approx_\Delta b\iff|a-b|<\infty$.
    Suppose for contradiction that there are $a,b\in\mathbb{Q}\times\mathbb{Z}$ with $|a-b|<\infty$ and $a\not\approx_\Delta b$.
    Then there exists $f\in\End(\mathbb{Q}.\Delta)$ with $|f(a)-f(b)|>|a-b|$.
    \begin{center}
        \begin{tikzcd}
            \mathbb{Q}.\Gamma\arrow[r, shift left, "e"] & \mathbb{Q}.\Gamma|_{\mathbb{Q}\times n\mathbb{Z}}\arrow[loop above, "h"]\arrow[l, shift left, "\operatorname{id}"]\arrow[r, shift left, "\iota"] & \mathbb{Q}.\Delta\arrow[l, shift left, "\iota^{-1}"]\arrow[loop right, "f"]
        \end{tikzcd}\
    \end{center}
    Hence, writing $na\coloneq\iota^{-1}(a)$ and $nb\coloneq\iota^{-1}(b)$, the function $g\coloneq\iota^{-1}\circ f\circ\iota$ is an endomorphism of $\mathbb{Q}.\Gamma|_{\mathbb{Q}\times n\mathbb{Z}}$ with $|g(na)-g(nb)|>|na-nb|$.
    Let $c,d\in\mathbb{Q}\times\mathbb{Z}$ with $c-d=na-nb\in n\mathbb{Z}$ and thus $c\approx_\Gamma d$.
    Then by Lemma~\ref{lem:pp_define_distances}, $\Diff_{c-d}$ or $\Dist_{c-d}$ is pp-definable in $\mathbb{Q}.\Gamma$, which means that $|e(c)-e(d)|=|c-d|$.
    Let $h\in\Aut(\mathbb{Q}\times\mathbb{Z},<)$ with $h(\{e(c),e(d)\})=\{na,nb\}$.
    Then $|(g\circ h\circ e)(c)-(g\circ h\circ e)(d)|=|g(na)-g(nb)|>|na-nb|=|c-d|$, which contradicts $c\approx_\Gamma d$.

    By Lemma~\ref{lem:pp_define_distances} for each $n\in\mathbb{N}$ the structure $\mathbb{Q}.\Delta$ pp-defines $\Diff_n$ or $\Dist_n$.
    If it pp-defines $\Diff_1$ we are in Case~\ref{item:defining_successor_succ}.
    Suppose it pp-defines $\Dist_1$.
    Then it pp-defines $\Dist_n$ for all $n\in\mathbb{N}$ and we are in Case~\ref{item:defining_successor_dist} or there is some $n\in\mathbb{N}_{>0}$ such that $\mathbb{Q}.\Delta$ pp-defines $\Diff_n$ and $\exists z(\Dist_1(x,y)\land\Dist_1^{n-1}(y,z)\land\Diff_n(x,z))$ is a pp-definiton of $\Succ$, which again means that we are in Case~\ref{item:defining_successor_succ}.
    \begin{center}
        \begin{tikzcd}[sep=huge]
            x \arrow[r, dash, "1" swap, "\Dist_1"] \arrow[rr, bend right, "\Diff_n", "n" swap] & y \arrow[r, dash, "n-1" swap, "\Dist_1^{n-1}"] & \exists z
        \end{tikzcd}\
    \end{center}
\end{proof}
\begin{example}
    Let $\Gamma\coloneq(\mathbb{Z},\Diff_{\{4,10,16\}})$.
    Then for any injective $f:\{0,1\}\times\mathbb{Q}\to\mathbb{Q}$ the map
    \[e:\mathbb{Q}\times\mathbb{Z}\to\mathbb{Q}\times\mathbb{Z},(p,a)\mapsto(f(a\bmod 2,p),a)\]
    is an endomorphism of $\mathbb{Q}.\Gamma$ with $|e(p,a)-e(p,a+1)|=\infty$.
    Since $(a,a+10),(a+6,a+10),(a+2,a+6)\in\Diff_{\{4,10,16\}}$ we get $a\sim_\Gamma a+2$ and hence $a\sim_\Gamma b\iff a-b\in 2\mathbb{Z}$.
    The function
    \[g:\mathbb{Q}\times\mathbb{Z}\to\mathbb{Q}\times\mathbb{Z},(p,a)\mapsto
    \begin{cases}
        e(p,a) & \text{if }a\equiv 0\mod{4}\\
        e(p,a)-1 & \text{if }a\equiv 1\mod{4}\\
        e(p,a)+6 & \text{if }a\equiv 2\mod{4}\\
        e(p,a)+5 & \text{if }a\equiv 3\mod{4}
    \end{cases}\]
    is an endomorphism of $\mathbb{Q}.\Gamma$ with $|g(a)-g(a+2)|=8$ and thus $a\not\approx_\Gamma a+2$.
    Since
    \[\exists z(\Diff_{\{\textcolor{red}{4},10,16\}}(x,y)\land\Diff_{\{4,10,\textcolor{red}{16}\}}(x,z)\land\Diff_{\{\textcolor{red}{4},10,16\}}^3(y,z))\]
    pp-defines $\Diff_4$ and ${\approx_\Gamma}\subseteq{\sim_\Gamma}$ it follows that $a\approx_\Gamma b\iff a-b\in 4\mathbb{Z}$.
    We have $g(\mathbb{Q}\times\mathbb{Z})\subseteq\mathbb{Q}\times 4\mathbb{Z}$ and hence $\mathbb{Q}.\Gamma$ is homomorphically equivalent to $\mathbb{Q}.\Gamma|_{\mathbb{Q}\times 4\mathbb{Z}}$.
    Let $\Delta\coloneq(\mathbb{Z},\Diff_{\{1,4\}})$.
    Then $\mathbb{Q}.\Delta$ is isomorphic to $\mathbb{Q}.\Gamma|_{\mathbb{Q}\times 4\mathbb{Z}}$ and pp-defines $\Succ$ via the formula
    \[\exists z(\Diff_{\{\textcolor{red}{1},4\}}(x,y)\land\Diff_{\{1,\textcolor{red}{4}\}}(x,z)\land\Diff_{\{\textcolor{red}{1},4\}}^3(y,z)).\]
\end{example}
\begin{lemma}[{\cite[see][Proposition~27]{bodirsky2016}}]\label{lem:hardness_from_distances}
    $\Csp(\mathbb{Z},\Dist_1,\Dist_3)$ is $\Np$-hard.
\end{lemma}
\begin{proof}
    The formula $\exists z(\Dist_3(x,z)\land\Dist_1(y,z))$ pp-defines $\Dist_{\{2,4\}}$.
    The map
    \[e:\mathbb{Z}\to\mathbb{Z},a\mapsto
    \begin{cases}
        a\bmod 6 & \text{if }a\in 2\mathbb{Z}\\
        (a-1)\bmod 6 & \text{otherwise}
    \end{cases}\]
    is an endomorphism of $(\mathbb{Z},\Dist_{\{2,4\}})$ with $e(\mathbb{Z})=\{0,2,4\}$.
    Hence, $\Csp(\{0,2,4\},\Dist_{\{2,4\}})$ reduces in polynomial time to $\Csp(\mathbb{Z},\Dist_1,\Dist_3)$.
    Since $(\{0,2,4\},\Dist_{\{2,4\}})$ is isomorphic to the complete graph $\mathbb{K}_3$ on three elements, its CSP is $\Np$-hard by \cite{hell1990}.
\end{proof}

\section{Tractable cases}\label{sec:tractability}
Our next goal is to show that the CSP of every finite-signature first-order reduct of $(\mathbb{Z},<)$ that is preserved by a so-called modular min or modular max polymorphism, or whose relations are all definable by certain Horn formulas, is tractable.
Rather than reproducing the ad-hoc algorithms of \cite{bodirsky2018}, some of which originate in \cite{bodirsky2016}, we show that these CSPs are already solved by well-known uniform algorithms.
Like the original algorithms, ours rely on \emph{sampling}: given an instance of the CSP, we can efficiently compute a finite substructure of the reduct, called the \emph{sample}, such that the instance is a yes-instance of the original CSP if and only if it is a yes-instance of the CSP of the sample.
If, in addition, the CSP of the sample can be solved by an algorithm running in time polynomial in the size of the input and in the size of the sample, then combining the two steps yields an efficient algorithm for the original CSP.

A convenient framework for reasoning about sampling is a generalization of CSPs, the \emph{promise constraint satisfaction problem}.
For relational $\tau$-structures $\mathcal{M}$ and $\mathcal{N}$ with a homomorphism $\mathcal{M}\to\mathcal{N}$, the problem $\Pcsp(\mathcal{M},\mathcal{N})$ takes as input a finite structure $\mathcal{A}$ of finite signature $\sigma\subseteq\tau$, and asks to distinguish the case that $\mathcal{A}$ maps homomorphically to $\mathcal{M}\restriction\sigma$ from the case that $\mathcal{A}$ does not map homomorphically to $\mathcal{N}\restriction\sigma$.
Since $\Pcsp(\mathcal{M},\mathcal{M})$ coincides with $\Csp(\mathcal{M})$, promise constraint satisfaction problems generalize constraint satisfaction problems.
\begin{definition}
    A \emph{uniform} algorithm is a polynomial time algorithm that accepts a pair of finite structures $\mathcal{A}$ and $\mathcal{M}$ with finite relational signatures if $\mathcal{A}$ maps homomorphically to $\mathcal{M}$.
    It is allowed to accept even when such a homomorphism does not exist.
\end{definition}
Assuming $\mathcal{M}$ is finite with signature $\tau$, we say a uniform algorithm \emph{solves} $\Pcsp(\mathcal{M},\mathcal{N})$ if, whenever it accepts $(\mathcal{A},\mathcal{M}\restriction\sigma)$ for a finite $\sigma\subseteq\tau$ and a finite $\sigma$-structure $\mathcal{A}$, then $\mathcal{A}$ maps homomorphically to $\mathcal{N}\restriction\sigma$.

The following lemma isolates what makes sampling work for first-order reducts of $(\mathbb{Z},<)$.
\begin{lemma}[{\cite[see][Lemma~5.4]{bodirsky2018}}]\label{lem:substitution}
    Let $\Gamma$ be a first-order reduct of $(\mathbb{Z},<)$ with a finite signature $\tau$.
    Then there exists $n\in\mathbb{N}$ such that for every $\tau$-structure $\mathcal{M}$ and all functions $f,g:M\to\mathbb{Q}\times\mathbb{Z}$ with
    \begin{enumerate}
        \item $f(a)<f(b)\iff g(a)<g(b)$ and
        \item $f(a)-f(b)=m\iff g(a)-g(b)=m$ for all $m\in\{-n,\dots,n\}$
    \end{enumerate}
    for all $a,b\in M$, we have $f\in\Hom(\mathcal{M},\mathbb{Q}.\Gamma)\iff g\in\Hom(\mathcal{M},\mathbb{Q}.\Gamma)$.
\end{lemma}
\begin{proof}
    By Fact~\ref{fact:quantifier_elimination_order}, for every $R\in\tau$ there exists a quantifier-free formula $\phi_R(\bar{x})$ over the language $\sigma\coloneq\{x+m<y\mid m\in\mathbb{Z}\}$ such that $(\mathbb{Q}\times\mathbb{Z},\sigma)\models\phi_R(\bar{a})\iff\bar{a}\in R^{\mathbb{Q}.\Gamma}$.
    Let $n$ be the maximum among all $|m|$ such that an atomic formula of the form $x+m<y$ appears in $\phi_R$ for some $R\in\tau$.

    If $a,b\in M$ with $|f(a)-f(b)|\leq n$ then $f(a)-f(b)=g(a)-g(b)$ and thus $f(a)+m<f(b)\iff g(a)+m<g(b)$ for all $m\in\mathbb{Z}$.
    Suppose that $a,b\in M$ with $|f(a)-f(b)|>n$.
    Then $|g(a)-g(b)|>n$ and for all $m\in\{-n,\dots,n\}$ we have $f(a)+m<f(b)\iff f(a)<f(b)\iff g(a)<g(b)\iff g(a)+m<g(b)$.
    Now, for all $R\in\tau$ (of arity $k$) and $\bar{a}\in M^k$ it follows that $f(\bar{a})\in R^{\mathbb{Q}.\Gamma}\iff(\mathbb{Q}\times\mathbb{Z},\sigma)\models\phi_R(f(\bar{a}))\iff(\mathbb{Q}\times\mathbb{Z},\sigma)\models\phi_R(g(\bar{a}))\iff g(\bar{a})\in R^{\mathbb{Q}.\Gamma}$.
\end{proof}
Using this, we can now capture the whole sampling argument.
\begin{lemma}\label{lem:sampling_for_fo_reducts_of_integers}
    Let $\Gamma$ be a finite-signature first-order reduct of $(\mathbb{Z},<)$.
    If there exists a uniform algorithm that solves $\Pcsp(\Gamma|_{[n]},\Gamma)$ for every $n\in\mathbb{N}$, then $\Csp(\Gamma)\in\Ptime$.
\end{lemma}
\begin{proof}
    By Lemma~\ref{lem:substitution} there exists $n\in\mathbb{N}$ such that a finite structure $\mathcal{A}$ maps homomorphically to $\Gamma$ iff it maps homomorphically to $\Gamma|_{[m]}$ for $m\coloneq(n+1)\cdot|A|$.
    Using the quantifier-free definitions of the relations of $\Gamma$ given by Fact~\ref{fact:quantifier_elimination_order}, we can compute $\Gamma|_{[m]}$ in time polynomial in $|A|$.
    Invoking the uniform algorithm that solves $\Pcsp(\Gamma|_{[m]},\Gamma)$ on $(\mathcal{A},\Gamma|_{[m]})$ now computes whether $\mathcal{A}$ homomorphically maps to $\Gamma$ in time polynomial in the size of $\mathcal{A}$.
\end{proof}

\subsection{Tractability of modular min and max}
We now aim to prove Proposition~\ref{prop:tractability_of_mod_max_min}: if a finite-signature first-order reduct of $(\mathbb{Z},<)$ admits a modular min or modular max polymorphism, then the uniform algorithm singleton-arc-consistency (\Sac{}) \cite[see][]{zhuk2026} solves the CSP of every finite substructure.
Combined with the sampling lemma, this yields tractability of $\Csp(\Gamma)$ for every finite-signature first-order reduct $\Gamma$ of $(\mathbb{Z},<)$ admitting such a polymorphism.
\begin{definition}
    For $d\in\mathbb{N}_{>0}$ let
    \[\Max_d:\mathbb{Z}^2\to\mathbb{Z},(a,b)\mapsto
    \begin{cases}
        \max(a,b) & \text{if }a\equiv b\mod{d}\\
        a & \text{otherwise.}
    \end{cases}\]
    The function $\Min_d$ is defined analogously.
\end{definition}
\textcite{bodirsky2018} already correctly identify preservation by a modular min or max polymorphism as a source of tractability, but for the proof they refer to \cite{bodirsky2016}, asserting that the local-finiteness assumption underlying that paper is not needed in this particular case.
This assertion is incorrect: the proof of \cite[Lemma~34]{bodirsky2016} does rely on local-finiteness, via the claim that every relation $R\subseteq\mathbb{Z}^k$ that is first-order definable in $(\mathbb{Z},\Succ)$ and preserved by $\Min_d$ or $\Max_d$ satisfies $a_i-a_1\equiv b_i-b_1\mod{d}$ for all $\bar{a},\bar{b}\in R$ and $i\in[k]$.
As the following example shows, this claim fails in general.
\begin{example}
    The structure $(\mathbb{Z},\{\bar{a}\in\mathbb{Z}^4\mid a_2-a_1\in\{2,4\},a_4-a_3=1\})$ is preserved by $\Min_2$ and $\Max_2$.
    However, its relation contains tuples $\bar{a}$ and $\bar{b}$ such that $a_1-a_3\not\equiv b_1-b_3\mod{2}$.
\end{example}
What the proof in \cite{bodirsky2016} actually needs is a suitable decomposability property for relations, which we establish in Lemma~\ref{lem:modular_relations_decompose}.
Combined with results from \cite{zhuk2026}, this decomposition already suffices to establish tractability, so we do not need to rely on the techniques of \cite{bodirsky2016} at all.

Towards Lemma~\ref{lem:modular_relations_decompose}, we first investigate what it means for a binary relation to be preserved by a modular min or max polymorphism.
This requires understanding which binary relations are definable in $(\mathbb{Z},<)$; while we are at it, we also classify which of them are already definable in $(\mathbb{Z},\Succ)$, for later use.
\begin{fact}\label{fact:quantifier_elimination_succ}
    The structures $(\mathbb{Z},\Succ)$ and $(\mathbb{Z},(\Diff_n)_{n\in\mathbb{Z}})$ are interdefinable and $\Th(\mathbb{Z},(\Diff_n)_{n\in\mathbb{Z}})$ has quantifier elimination.
\end{fact}
\begin{lemma}\label{lem:fo_definable_binary_relations}
    A binary relation $R\subseteq\mathbb{Z}^2$ is first-order definable over $(\mathbb{Z},<)$ if and only if $R=\Diff_S$ for some $S\subseteq\mathbb{Z}$ such that $\{s\in S\mid s+1\notin S\}$ is finite.
    It is first-order definable over $(\mathbb{Z},\Succ)$ iff $S$ is finite or cofinite.
\end{lemma}
\begin{proof}
    For a set $S\subseteq\mathbb{Z}$ let $S'\coloneq\{s\in S\mid s+1\notin S\}$.

    Suppose $R$ is first-order definable over $(\mathbb{Z},<)$.
    Then by Fact~\ref{fact:quantifier_elimination_order} there exists a quantifier-free definition $\phi(x,y)$ of $R$ in $(\mathbb{Z},(x+n<y)_{n\in\mathbb{Z}})$.
    Since $\lnot(x+n<y)$ is equivalent to $y+(-n-1)<x$ we can even assume that $\phi$ is positive.
    We prove the claim by structural induction on $\phi$.
    If $\phi(x,y)=(x+n<y)$, then $R=\Diff_S$ for $S\coloneq\{n+1,n+2,\dots\}$ and $S'=\emptyset$.
    If $\phi(x,y)=(y+n<x)$, we have $R=\Diff_S$ for $S\coloneq\{-n-1,-n-2,\dots\}$ and $S'=\{-n-1\}$.
    Suppose $\psi(x,y)$ and $\chi(x,y)$ define $\Diff_U$ and $\Diff_V$ respectively with $U'$ and $V'$ finite.
    If $\phi=\psi\land\chi$ or $\phi=\psi\lor\chi$ then $R=\Diff_S$ for $S\coloneq U\cap V$ or $S\coloneq U\cup V$.
    In both cases we have $S'\subseteq U'\cup V'$ and hence $S'$ is finite.

    For the other direction let $S\subseteq\mathbb{Z}$ with $S'$ finite.
    If $S$ is finite, $\Diff_S$ is defined by $\bigvee_{s\in S}\Succ^s(x,y)$.
    If $S$ is cofinite, $\lnot\Diff_{\mathbb{Z}\setminus S}(x,y)$ defines $\Diff_S$.
    Suppose $S$ is neither finite nor cofinite.
    Since $S'$ is finite, there exists $n\in\mathbb{N}$ such that either $-m\notin S$ and $m\in S$ for all $m>n$ or $-m\in S$ and $m\notin S$ for all $m>n$.
    Because both cases are symmetric let us assume that the former holds.
    Then $x+n<y\lor\Diff_{S\cap\{-n,\dots,n\}}(x,y)$ defines $\Diff_S$.

    The proof for definability over $(\mathbb{Z},\Succ)$ is similar using Fact~\ref{fact:quantifier_elimination_succ}.
\end{proof}
\begin{lemma}\label{lem:infinite_binary_modular_relations_are_full}
    Let $R\subseteq\mathbb{Z}^2$ be first-order definable over $(\mathbb{Z},<)$ and preserved by $\Min_d$ or $\Max_d$ for some $d>1$.
    If for some $a\in\mathbb{Z}$ there are infinitely many $b\in\mathbb{Z}$ with $(a,b)\in R$ then $R=\mathbb{Z}^2$.
\end{lemma}
\begin{proof}
    Suppose $R$ is preserved by $\Max_d$.
    The proof for $\Min_d$ is analogous.
    Suppose w.l.o.g.\ that there are infinitely many $b\geq a$ such that $(a,b)\in R$.
    Then by Lemma~\ref{lem:fo_definable_binary_relations} there is some $c\in\mathbb{Z}$ such that $(a,b)\in R$ for all $b\geq c$.
    Assuming that $R\neq\mathbb{Z}^2$ we can pick $c$ minimal with that property.
    Then $(a-d,c-1)\in R$ and thus $(a,c-1)=(\Max_d(a-d,a),\Max_d(c-1,c))\in R$ contradicting the minimality of $c$.
\end{proof}
Given a relation $R\subseteq\mathbb{Z}^k$ and $d\in\mathbb{N}$, for $i,j\in[k]$ we write $i\sim j$ if $a_i-a_j\equiv b_i-b_j\mod{d}$ for all $\bar{a},\bar{b}\in R$.
Both $R$ and $d$ will always be clear from context.
\begin{lemma}\label{lem:shifting_modular_blocks}
    Let $R\subseteq\mathbb{Z}^k$ be first-order definable over $(\mathbb{Z},<)$ and preserved by $\Min_d$ or $\Max_d$.
    Then for all $\bar{a}\in R$ and $n_{[1]_\sim},\dots,n_{[k]_\sim}\in\mathbb{Z}$ we have $(a_i+n_{[i]_\sim})_{i\in[k]}\in R$.
\end{lemma}
\begin{proof}
    If $d=1$ we are done, so assume $d>1$.
    Suppose $R$ is preserved by $\Max_d$.
    The proof for $\Min_d$ is analogous.
    It suffices to show that for every $\sim$-equivalence class $I\subseteq[k]$, $\bar{a}\in R$ and $n\in\mathbb{Z}$ we have $\bar{a}'\in R$ for
    \[a_i'\coloneq
    \begin{cases}
        a_i+n & \text{if }i\in I\\
        a_i & \text{if }i\notin I.
    \end{cases}\]
    Let $\bar{a}_0\coloneq\bar{a}$ and $\bar{a}_1,\dots,\bar{a}_\ell\in R$ such that for all $i\in I$ and $j\in[k]\setminus I$
    there is $r\in[\ell]$ with $a_i\equiv a_r^{(i)}\mod{d}$ and $a_j\not\equiv a_r^{(j)}\mod{d}$.

    Pick $m\in d\mathbb{Z}$ such that $m\geq|a_r^{(i)}-a_s^{(j)}|$ for all $i,j\in[k]$ and $r,s\in\{0,\dots,\ell\}$.
    Extending addition to tuples componentwise, let $\bar{b}_0\coloneq\bar{a}$ and
    \[\bar{b}_r\coloneq\Max_d\big(\Max_d(\bar{b}_{r-1},\bar{a}_r+m_r),\Max_d(\bar{a}_r+m_r,\bar{b}_{r-1}+2m_r)\big)\in R\]
    for all $r\in[\ell]$ where $m_r\coloneq(3^r-1)m/2$.
    Suppose $|a_r^{(i)}-b_{r-1}^{(i)}|\leq m_r$ for some $r\in[\ell]$, which holds for $r=1$ by the choice of $m$.
    Then
    \[b_r^{(i)}=
    \begin{cases}
        b_{r-1}^{(i)}+2m_r & \text{if }a_r^{(i)}\equiv b_{r-1}^{(i)}\mod{d}\\
        b_{r-1}^{(i)} & \text{otherwise}
    \end{cases}\]
    and hence, if $r<\ell$, we have $|a_{r+1}^{(i)}-b_r^{(i)}|\leq|a_{r+1}^{(i)}-a_r^{(i)}|+|a_r^{(i)}-b_{r-1}^{(i)}|+|b_{r-1}^{(i)}-b_r^{(i)}|\leq m+m_r+2m_r=m_{r+1}$.
    Now, let $\bar{b}\coloneq\bar{b}_\ell$.
    Then by induction and since $a_i\equiv b_r^{(i)}\mod{d}$ for all $i\in[k]$ and $r\in[\ell]$, we get
    \begin{align*}
        b_i&=a_i+\sum\{2m_r\mid r\in[\ell],a_i\equiv a_r^{(i)}\mod{d}\}\\
        &=a_i+m\cdot\underbrace{\sum\{3^r-1\mid r\in[\ell],a_i\equiv a_r^{(i)}\mod{d}\}}_{\eqcolon N_i}
    \end{align*}
    for all $i\in[k]$.

    Let $N\coloneq\sum_{r\in[\ell]}(3^r-1)$ and note that $N_i\leq N$ for all $i\in[k]$.
    For $i\in I$ we have $a_i\equiv a_r^{(i)}\mod{d}$ for all $r\in[\ell]$ and thus $N_i=N$.
    Moreover, by our choice of $\bar{a}_1,\dots,\bar{a}_\ell$ we get $N_j<N$ for all $j\in[k]\setminus I$.
    Now, let $\bar{c}\coloneq\Max_d(\bar{a},\bar{b}-m(N-1))\in R$.
    Then
    \[c_i=
    \begin{cases}
        a_i+m & \text{if }i\in I\\
        a_i & \text{if }i\notin I
    \end{cases}\]
    for all $i\in[k]$.
    Hence, $\bar{c}$ witnesses that the binary relation pp-defined by
    \[\exists\bar{z}(R\bar{z}\land\bigwedge_{i\in I}\Succ^{a_i}(x,z_i)\land\bigwedge_{i\in[k]\setminus I}\Succ^{a_i}(y,z_i))\]
    contains $(m,0)$.
    Since we can pick $m$ arbitrarily large, by Lemma~\ref{lem:infinite_binary_modular_relations_are_full} the relation is full, which concludes the proof.
\end{proof}
\begin{lemma}\label{lem:modular_relations_decompose}
    Let $R\subseteq\mathbb{Z}^k$ be first-order definable over $(\mathbb{Z},<)$ and preserved by $\Min_d$ or $\Max_d$.
    Then
    \[R=\bigcap_{i\in[k]}\{(a_1,\dots,a_k)\in\mathbb{Z}^k\mid(a_j)_{j\sim i}\in\Proj_{[i]_\sim}(R)\}.\]
\end{lemma}
\begin{proof}
    Suppose $R$ is preserved by $\Max_d$.
    The proof for $\Min_d$ is analogous.
    The inclusion from left to right is trivial.
    For the other direction we show that for all $\bar{a},\bar{b}\in R$ and $i\in[k]$ the tuple $\bar{a}'$ with
    \[a_j'\coloneq
    \begin{cases}
        a_j & \text{if }i\not\sim j\\
        b_j & \text{if }i\sim j
    \end{cases}\]
    for all $j\in[k]$ is in $R$.
    The claim then follows by induction.
    Let $n_{[1]_\sim},\dots,n_{[k]_\sim}\in\mathbb{Z}$ such that for all $j\in[k]$ we have
    both $a_j\equiv b_j+n_{[j]_\sim}\mod{d}$ and $a_j\leq b_j+n_{[j]_\sim}$ if and only if $j\sim i$.
    By Lemma~\ref{lem:shifting_modular_blocks} $(b_j+n_{[j]_\sim})_{j\in[k]}\in R$ and thus $\bar{c}$ with
    \[c_j\coloneq\Max_d(a_j,b_j+n_{[j]_\sim})=
    \begin{cases}
        a_j & \text{if }j\not\sim i\\
        b_j+n_{[i]_\sim} & \text{if }j\sim i
    \end{cases}\]
    for all $j\in[k]$ is in $R$.
    Applying Lemma~\ref{lem:shifting_modular_blocks} again yields $\bar{a}'\in R$.
\end{proof}
With Lemma~\ref{lem:modular_relations_decompose} in place, we can now use techniques by \textcite{zhuk2026} to prove tractability of the CSP for first-order reducts of $(\mathbb{Z},<)$ admitting a modular min or max polymorphism.
\begin{definition}
    Let $n,k_1,\dots,k_n\in\mathbb{N}_{>0}$.
    A tuple $(a_{1,1},\dots,a_{1,k_1};\dots;a_{n,1},\dots,a_{n,k_n})$ is \emph{$(k_1,\dots,k_n)$-palette} if for all $i\in[n]$ and $j\in[k_i]$ there is some $\ell\in[n]$ such that $a_{i,j}=a_{\ell,1}=\dots=a_{\ell,k_\ell}$.
    If a tuple is $(k_1,\dots,k_n)$-palette for $k_1=\dots=k_n\eqcolon m$, we call it \emph{$(n\times m)$-palette}.
\end{definition}
\begin{definition}
    Let $n,k_1,\dots,k_n\in\mathbb{N}_{>0}$.
    A function $f:A^{k_1+\dots+k_n}\to B$ is \emph{$(k_1,\dots,k_n)$-palette totally symmetric} if there exists a function $g:\Pow(A)^n\to B$ such that
    \[f(a_{1,1},\dots,a_{1,k_1};\dots;a_{n,1},\dots,a_{n,k_n})=g(\{a_{1,1},\dots,a_{1,k_1}\},\dots,\{a_{n,1},\dots,a_{n,k_n}\})\]
    for all $(k_1,\dots,k_n)$-palette tuples from $A^{k_1+\dots+k_n}$.
    If a function is $(k_1,\dots,k_n)$-palette totally symmetric for $k_1=\dots=k_n\eqcolon m$, we call it \emph{$(n\times m)$-palette totally symmetric}.
\end{definition}
\begin{lemma}\label{lem:mod_max_min_palette_totally_symmetric}
    Let $\Gamma$ be a first-order reduct of $(\mathbb{Z},<)$ that admits $\Min_d$ or $\Max_d$ as a polymorphism.
    Then $\Gamma$ admits a conservative $(k_1,\dots,k_n)$-palette totally symmetric polymorphism for all $n,k_1,\dots,k_n\in\mathbb{N}_{>0}$.
\end{lemma}
\begin{proof}
    From any conservative $(n\times m)$-palette totally symmetric polymorphism, by identifying variables, we can construct a conservative $(k_1,\dots,k_n)$-palette totally symmetric polymorphism for arbitrary $k_1,\dots,k_n\in[m]$.
    Hence, it is enough to show that $\Gamma$ admits a conservative $(n\times m)$-palette totally symmetric polymorphism for all $n,m\in\mathbb{N}_{>0}$.
    Suppose $\Max_d$ is a polymorphism of $\Gamma$.
    The proof for $\Min_d$ is analogous.
    Let $f:\mathbb{Z}^{nm}\to\mathbb{Z}$ such that
    \[f(a_{1,1},\dots,a_{1,m};\dots;a_{n,1},\dots,a_{n,m})=\max(a_{i,1},\dots,a_{i,m})\]
    where $i\in[n]$ is minimal with $a_{i,1}\equiv\dots\equiv a_{i,m}\mod{d}$ if such an $i$ exists and
    \[f(a_{1,1},\dots,a_{1,m};\dots;a_{n,1},\dots,a_{n,m})=a_{1,1}\]
    otherwise.
    Clearly, $f$ is conservative, i.e., $f(\bar{a})\in\bar{a}$ for all $\bar{a}\in\mathbb{Z}^{nm}$.

    Let $\bar{a}\coloneq(a_{1,1},\dots,a_{1,m};\dots;a_{n,1},\dots,a_{n,m})\in\mathbb{Z}^{nm}$ be $(n\times m)$-palette.
    Then there exists $i\in[n]$ such that $a_{i,1}=\dots=a_{i,m}$ and hence $f(\bar{a})=\max(a_{j,1},\dots,a_{j,m})$ for some $j\in[n]$.
    Thus, $f$ is $(n\times m)$-palette totally symmetric.

    Next, we show that $f$ is a polymorphism of $\Gamma$.
    Let $R$ be a relation of $\Gamma$ and $\bar{a}_{i,j}\coloneq(a_{i,j}^{(1)},\dots,a_{i,j}^{(k)})\in R$ for all $i\in[n]$ and $j\in[m]$.
    By Lemma~\ref{lem:modular_relations_decompose} it suffices to show that $f$ preserves $\Proj_{[i]_\sim}(R)$ for every $i\in[k]$, so assume for simplicity that ${\sim}=[k]^2$.
    Suppose there are $\ell\in[n]$ and $i\in[k]$ such that $a_{\ell,1}^{(i)}\equiv\dots\equiv a_{\ell,m}^{(i)}\mod{d}$ and let $\ell$ be minimal with this property.
    Since ${\sim}=[k]^2$, for any $j\in[k]$ we get $a_{\ell,1}^{(j)}\equiv\dots\equiv a_{\ell,m}^{(j)}\mod{d}$.
    Hence, applying $f$ coordinatewise yields
    \begin{align*}
        f(\bar{a}_{1,1},\dots,\bar{a}_{1,m};\dots;\bar{a}_{n,1},\dots,\bar{a}_{n,m})&=(f(a_{1,1}^{(i)},\dots,a_{1,m}^{(i)};\dots;a_{n,1}^{(i)},\dots,a_{n,m}^{(i)}))_{i\in[k]}\\
        &=(\max(a_{\ell,1}^{(i)},\dots,a_{\ell,m}^{(i)}))_{i\in[k]}\\
        &=\max(\bar{a}_{\ell,1},\dots,\bar{a}_{\ell,m})\\
        &=\Max_d(\dots(\Max_d(\bar{a}_{\ell,1},\bar{a}_{\ell,2}),\bar{a}_{\ell,3})\dots)\in R.
    \end{align*}
    On the other hand, if there are no $\ell\in[n]$ and $i\in[k]$ with $a_{\ell,1}^{(i)}\equiv\dots\equiv a_{\ell,m}^{(i)}$,
    we get
    \[f(\bar{a}_{1,1},\dots,\bar{a}_{1,m};\dots;\bar{a}_{n,1},\dots,\bar{a}_{n,m})=\bar{a}_{1,1}\in R.\]
    Therefore, $f$ preserves $R$.
\end{proof}
By \cite[Theorem~3.5]{zhuk2026} the next proposition follows immediately.
\begin{proposition}\label{prop:tractability_of_mod_max_min}
    Let $\Gamma$ be a finite-signature first-order reduct of $(\mathbb{Z},<)$ that admits $\Min_d$ or $\Max_d$ as a polymorphism.
    Then for any finite $A\subseteq\mathbb{Z}$ the problem $\Csp(\Gamma|_A)$ is solved by \Sac{}.
\end{proposition}

\subsection{Horn relations}
Another source of tractability for $\Csp(\Gamma)$ arises when $\Gamma$ is a first-order reduct of $(\mathbb{Z},\Succ)$ whose relations can be expressed as a collection of implications stating that fixed differences among some pairs of variables force a fixed difference between two others.
To make this precise, call a first-order formula \emph{Horn} if it is quantifier-free, in conjunctive normal form, and each clause contains at most one positive literal.
A relation that is first-order definable over $(\mathbb{Z},\Succ)$ or $(\mathbb{Q}\times\mathbb{Z},\Succ)$ is then \emph{Horn-definable} if it can be defined by a Horn formula over $(\mathbb{Z},(\Diff_n)_{n\in\mathbb{Z}})$ or $(\mathbb{Q}\times\mathbb{Z},(\Diff_n)_{n\in\mathbb{Z}})$, respectively, and we call a first-order reduct of $(\mathbb{Z},\Succ)$ or $(\mathbb{Q}\times\mathbb{Z},\Succ)$ \emph{Horn} if all of its relations are Horn-definable.
Using techniques by \textcite{zhuk2026} once more, we show that the constraint-singleton-arc-consistency algorithm (\Csac{}) efficiently solves the CSP of every such Horn reduct.
\begin{definition}
    Let $\mathcal{A}$ and $\mathcal{B}$ be relational structures over the same signature $\tau$.
    For $n,k_1,\dots,k_n\in\mathbb{N}_{>0}$ with $m\coloneq k_1+\dots+k_n$, by $(k_1,\dots,k_n)\textnormal{-}\Pol(\mathcal{A},\mathcal{B})$ we denote the set of functions $f:A^m\to B$ such that for every $R\in\tau$ and every $(k_1,\dots,k_n)$-palette tuple $(\bar{a}_1,\dots,\bar{a}_m)\in(R^\mathcal{A})^m$ we have $f(\bar{a}_1,\dots,\bar{a}_m)\in R^\mathcal{B}$.
\end{definition}
\begin{lemma}\label{lem:horn_palette_totally_symmetric}
    Let $\Gamma$ be a finite-signature first-order Horn reduct of $(\mathbb{Z},\Succ)$.
    Then $(k_1,\dots,k_\ell)\textnormal{-}\Pol(\Gamma|_{[n]},\Gamma)$ contains a $(k_1,\dots,k_\ell)$-palette totally symmetric function for all $n\in\mathbb{N}$ and $\ell,k_1,\dots,k_\ell\in\mathbb{N}_{>0}$.
\end{lemma}
\begin{proof}
    For every relation $R$ of $\Gamma$ fix a defining Horn formula and let $t'\in\mathbb{N}$ such that $t'\geq|t|$ for every negative literal $\lnot\Diff_t(x_i,x_j)$ appearing in one of these formulas.
    It is not hard to see that there exists an injective function $g:[n]^\ell\to\mathbb{Z}$ preserving $\Succ$ such that $|g(\bar{a})-g(\bar{b})|>t'$ whenever there are $i,j\in[\ell]$ with $a_i-b_i\neq a_j-b_j$.
    We claim that for $m\coloneq k_1+\dots+k_\ell$
    \[f:[n]^m\to\mathbb{Z},(a_{1,1},\dots,a_{1,k_1};\dots;a_{\ell,1},\dots,a_{\ell,k_\ell})\mapsto g(\min\bar{a}_1,\dots,\min\bar{a}_\ell)\]
    is the desired function.
    Clearly, $f$ is $(k_1,\dots,k_\ell)$-palette totally symmetric, so we just have to show that it is in $(k_1,\dots,k_\ell)\textnormal{-}\Pol(\Gamma|_{[n]},\Gamma)$.
    Let $R$ be a relation of $\Gamma$, $(\bar{a}_{1,1},\dots,\bar{a}_{1,k_1};\dots;\bar{a}_{\ell,1},\dots,\bar{a}_{\ell,k_\ell})\in(R|_{[n]})^m$ a $(k_1,\dots,k_\ell)$-palette tuple and $\bar{b}\coloneq f(\bar{a}_{1,1},\dots,\bar{a}_{1,k_1};\dots;\bar{a}_{\ell,1},\dots,\bar{a}_{\ell,k_\ell})$.
    Since $\Gamma$ is Horn, $R$ is definable by a Horn-formula $\phi(\bar{x})$ in $(\mathbb{Z},(\Diff_n)_{n\in\mathbb{Z}})$.
    Let $\psi$ be a clause of $\phi$.
    If $\psi$ contains a positive literal $\Diff_t(x_r,x_s)$ which is satisfied by $\bar{a}_{i,j}$ for all $i\in[\ell]$ and $j\in[k_i]$, then $\min_{j\in k_i}a_{i,j}^{(s)}-\min_{j\in k_i}a_{i,j}^{(r)}=t$ for all $i\in[\ell]$ and since $g$ preserves $\Succ$ it follows that $\bar{b}$ satisfies $\Diff_t(x_r,x_s)$ and thus $\psi$.
    Otherwise there are $i\in[\ell]$ and $j\in[k_i]$ such that $\bar{a}_{i,j}$ satisfies a negative literal $\lnot\Diff_t(x_r,x_s)$ of $\psi$.
    By the palette property there exists $i'\in[\ell]$ such that $\bar{a}_{i,j}=\bar{a}_{i',1}=\dots=\bar{a}_{i',k_{i'}}$.
    Now, $\min_{j\in[k_{i'}]}a_{i',j}^{(s)}-\min_{j\in[k_{i'}]}a_{i',j}^{(r)}=a_{i,j}^{(s)}-a_{i,j}^{(r)}\neq t$.
    By our choice of $g$, it follows that $\bar{b}$ satisfies $\lnot\Diff_t(x_r,x_s)$ and hence $\psi$.
    Thus, $\bar{b}$ satisfies every clause of $\phi$, which means $\bar{b}\in R$.
\end{proof}
By \cite[Lemma~6.13]{zhuk2026} we get the following proposition.
\begin{proposition}\label{prop:tractable_horn}
    Let $\Gamma$ be a finite-signature first-order reduct of $(\mathbb{Z},\Succ)$.
    If $\Gamma$ is Horn, then for every $n\in\mathbb{N}$ the algorithm \Csac{} solves $\Pcsp(\Gamma|_{[n]},\Gamma)$.
\end{proposition}

\section{Hard cases}\label{sec:hardness}
Having covered all sources of tractability for the CSP of first-order reducts of $(\mathbb{Z},<)$, we now turn to sources of $\Np$-hardness.
All of our hardness proofs proceed by reduction from an $\Np$-hard finite-domain CSP, using the following lemma.
\begin{lemma}\label{lem:use_relation_to_reduce_to_finite_csp}
    Let $\Gamma$ be a finite-signature first-order reduct of $(\mathbb{Z},<)$.
    If $\Gamma$ pp-defines $\Diff_S$ for some $S\subseteq\mathbb{Z}$, then $\Csp(\Gamma|_S)\leq_\Ptime\Csp(\Gamma)$.
\end{lemma}
\begin{proof}
    Let $\phi$ be an instance of $\Csp(\Gamma|_S)$ and $x$ a fresh variable which does not occur in $\phi$.
    Let $\psi(x)$ be the formula $\phi$ with every atomic subformula $Rx_1\dots x_k$ replaced by
    \[Rx_1\dots x_k\land\bigwedge_{i\in[k]}\Diff_S(x,x_i).\]
    Then it is easy to see that $\Gamma|_S\models\phi\iff\Gamma\models\psi(0)\iff\Gamma\models\exists x\psi$.
    Since $\exists x\psi$ can be computed from $\phi$ in polynomial time we get $\Csp(\Gamma|_S)\leq_\Ptime\Csp(\Gamma)$.
\end{proof}

\subsection{Hardness from binary relations}
We begin by identifying cases where binary relations alone already cause $\Np$-hardness.
The main role here is played by relations of the form $\Diff_D$, where, for some $d\in\mathbb{N}$, $D$ is a convex subset of $d\mathbb{Z}$, possibly shifted by a constant.
\begin{definition}
    A set $D\subseteq\mathbb{Z}$ is a \emph{$d$-progression}, for $d\in\mathbb{N}$, if $a-b\in d\mathbb{Z}$ for all $a,b\in D$, and $a+d\in D$ whenever $a\in D$ is not the maximum of $D$.
\end{definition}
Since we always assume $\Succ$ to be present, shifts by constants will not concern us much, as the following observation shows.
\begin{observation}\label{obs:shift_binary_relations}
    Let $S\subseteq\mathbb{Z}$ and $n\in\mathbb{Z}$.
    Then in $(\mathbb{Z},\Succ,\Diff_S)$ the formula $\exists z(\Diff_S(x,z)\land\Succ^n(y,z))$ pp-defines $\Diff_{S-n}$.
\end{observation}
We give a simpler proof of the next lemma than the original one by \textcite{bodirsky2016}.
\begin{lemma}[{\cite[Lemma~43]{bodirsky2016}}]\label{lem:pp_defining_gcd_progressions}
    Let $S\subseteq\mathbb{Z}$ be finite with $|S|>1$ and $d\coloneq\gcd\{s-t\mid s,t\in S\}$.
    Then $(\mathbb{Z},\Succ,\Diff_S)$ pp-defines $\Diff_D$ for every finite $d$-progression $D\subseteq\mathbb{Z}$.
\end{lemma}
\begin{proof}
    Let $S=\{s_1,\dots,s_k\}$.
    By Observation~\ref{obs:shift_binary_relations} we can assume w.l.o.g.\ that $0\in S$.
    Then $d=\gcd(S)$ and by B\'{e}zout's identity there are $n_1,\dots,n_k\in\mathbb{Z}$ such that $n_1s_1+\dots+n_ks_k=d$.
    Hence,
    \[\exists x_2\dots x_k\bigwedge_{i\in[k]}\Diff_S^{n_i}(x_i,x_{i+1})\]
    pp-defines $\Diff_T$ for some finite $T\subseteq d\mathbb{Z}$ with $0,d\in T$.
    Again by Observation~\ref{obs:shift_binary_relations} we can pp-define $\Diff_U$ for $U\coloneq T-\min T$.
    Now, for $i\in\mathbb{N}$ let $U_i\subseteq\mathbb{Z}$ such that $\Diff_{U_i}=\Diff_U^i$, i.e.,
    \[U_i=
    \begin{cases}
        \{0\} & \text{if }i=0\\
        U_{i-1}+U & \text{if }i>0.
    \end{cases}\]
    Since $0\in U\subseteq d\mathbb{N}$ we get $0\in U_i\subseteq d\mathbb{N}$ for all $i\in\mathbb{N}$.
    Therefore, $\max U_i=i\cdot\max U$ and hence
    \begin{equation}\label{eq:pp_defining_gcd_progressions_upper_bound}
        |U_i|\leq\max U_i/d+1=i\cdot\max U/d+1
    \end{equation}
    for all $i\in\mathbb{N}$.
    Now, for all $i\in\mathbb{N}$ let $D_i$ be the largest $d$-progression contained in $U_i$ and note that $|D_1|>1$.
    Since $D_i+D_1\subseteq U_{i+1}$ is a $d$-progression we get $|D_{i+1}|\geq|D_i+D_1|=|D_i|+|D_1|-1\geq|D_i|+1$ for all $i\in\mathbb{N}$.
    Thus, there exists some $n\in\mathbb{N}$ such that for all $i\geq n$ we have $\max D_i-\min D_i\geq\max U-d$, which means that $D_i+\{0,\max U\}\subseteq U_{i+1}$ is a $d$-progression.
    Hence,
    \begin{equation}\label{eq:pp_defining_gcd_progressions_lower_bound}
        |D_i|\geq|D_{i-1}+\{0,\max U\}|=|D_{i-1}|+\max U/d\geq|D_n|+(i-n)\max U/d
    \end{equation}
    for all $i>n$.
    Combining \eqref{eq:pp_defining_gcd_progressions_upper_bound} and \eqref{eq:pp_defining_gcd_progressions_lower_bound} yields
    \[|U_i\setminus D_i|=|U_i|-|D_i|\leq i\cdot\max U/d+1-(|D_n|+(i-n)\max U/d)=\underbrace{n\cdot\max U/d-|D_n|+1}_{\eqcolon r}\]
    for all $i>n$.
    Let $m\in\mathbb{N}$ such that $|D_m|>r+1$ and define
    \[V_i\coloneq
    \begin{cases}
        U_m & \text{if }i=0\\
        V_{i-1}\cap(V_{i-1}+d) & \text{if }i>0
    \end{cases}\]
    for all $i\in\mathbb{N}$.
    Then
    \[V_i=\{\max E\mid E\subseteq U_m\text{ is a $d$-progression with }|E|=i+1\}\]
    for all $i\in\mathbb{N}$.
    By our choice of $m$ there are exactly two $d$-progressions of size $|D_m|-1>r$ in $U_m$, namely $D_m\setminus\{\min D_m\}$ and $D_m\setminus\{\max D_m\}$.
    Thus, $V\coloneq V_{|D_m|-2}=\{\max D_m-d,\max D_m\}$.
    Using Observation~\ref{obs:shift_binary_relations} and the fact that $\Diff_{V_i}(x,y)\land\Diff_{V_i+d}(x,y)$ pp-defines $\Diff_{V_{i+1}}$ for all $i\in\mathbb{N}$, we get that $\Diff_{\{0,d\}}$ is pp-definable.

    Finally, let $D\subseteq\mathbb{Z}$ be an arbitrary non-empty finite $d$-progression.
    Then $\Diff_{\{0,d\}}^{|D|-1}=\Diff_{D-\min D}$ and by Observation~\ref{obs:shift_binary_relations} we can indeed pp-define $\Diff_D$.
\end{proof}
\textcite{bodirsky2016} prove the following lemma using Lemmas~\ref{lem:pp_defining_gcd_progressions} and~\ref{lem:pairwise_different_distances_are_hard}.
Their proof of Lemma~\ref{lem:pairwise_different_distances_are_hard}, however, is incomplete; we therefore first prove the next lemma, using techniques by \textcite{barto2009}, and then use it to give an independent proof of Lemma~\ref{lem:pairwise_different_distances_are_hard}.
\begin{lemma}[{\cite[Lemma~44]{bodirsky2016}}]\label{lem:non_progressions_are_hard}
    Let $S\subseteq\mathbb{Z}$ be finite and not a $d$-progression for any $d\in\mathbb{N}$.
    Then $\Csp(\mathbb{Z},\Succ,\Diff_S)$ is $\mathsf{NP}$-hard.
\end{lemma}
\begin{proof}
    Let $d\coloneq\gcd\{s-t\mid s,t\in S\}$.
    Since $S$ is not a $d$-progression there exists $a\in\{\min S,\min S+d,\dots,\max S\}\setminus S$.
    Let $T\coloneq S-a\subseteq d\mathbb{Z}$ and note that by Observation~\ref{obs:shift_binary_relations}, $\Diff_T$ is pp-definable using $\Succ$ and $\Diff_S$.
    Since $0\notin T$ the graph $\mathcal{G}\coloneq(\mathbb{Z},\Diff_T)$ has no loops.

    First, we claim that $\mathcal{G}$ has algebraic length~1,
    i.e., there are
    \[(\delta_1,b_1,b_2),(\delta_2,b_2,b_3),\dots,(\delta_{k-1},b_{k-1},b_k),(\delta_k,b_k,b_1)\in\{(1,b,c),(-1,c,b)\mid(b,c)\in\Diff_T\}\]
    such that $|\sum_{i\in[k]}\delta_i|=1$.
    Let $T=\{t_1,\dots,t_\ell\}$.
    It is not hard to see, that $\mathcal{G}$ has algebraic length~1 iff
    \begin{align*}
        x_1t_1+\dots+x_\ell t_\ell&=y_1t_1+\dots+y_\ell t_\ell\\
        x_1+\dots+x_\ell&=y_1+\dots+y_\ell+1
    \end{align*}
    has a solution in $\mathbb{N}$ or equivalently in $\mathbb{Z}$.
    Clearly, this is equivalent to
    \begin{align*}
        z_1t_1+\dots+z_\ell t_\ell&=0\\
        z_1+\dots+z_\ell&=1
    \end{align*}
    having a solution in $\mathbb{Z}$.
    Combining both equations we get
    \[z_1t_1+\dots+z_{\ell-1}t_{\ell-1}+(1-z_1-\dots-z_{\ell-1})t_\ell=0\]
    and equivalently
    \[z_1(t_1-t_\ell)+\dots+z_{\ell-1}(t_{\ell-1}-t_\ell)=-t_\ell.\]
    Now, this equation has a solution in $\mathbb{Z}$ iff $\gcd(t_1-t_\ell,\dots,t_{\ell-1}-t_\ell)$ divides $t_\ell$.
    But $\gcd(t_1-t_\ell,\dots,t_{\ell-1}-t_\ell)=d$, which divides $t_\ell$.
    This proves the claim.

    Since $T\subseteq d\mathbb{Z}$ it follows that $\mathcal{G}|_D$ has algebraic length~1 for every finite but sufficiently large $d$-progressions $D\subseteq\mathbb{Z}$.

    Next, we show that for large finite $d$-progressions $D\subseteq\mathbb{Z}$ the induced subgraph $\mathcal{G}|_D$ is smooth,
    i.e., for all $a\in D$ there are $b,c\in D$ such that $(a,b),(c,a)\in\Diff_T$.
    Note, that $\min T<0<\max T$.
    So if $a\geq\min D-\min T$ or $a\leq\max D-\max T$ there exists $b\in D$ with $(a,b)\in\Diff_T$, namely $b\coloneq a+\min T$ or $b\coloneq a+\max T$ respectively.
    Similarly, if $a\geq\min D+\max T$ or $a\leq\max D+\min T$ there exists $c\in D$ with $(c,a)\in\Diff_T$, namely $c\coloneq a-\max T$ or $c\coloneq a-\min T$ respectively.
    Hence, picking $D$ such that $\max D-\min D\geq \max T-\min T$ is sufficient.

    Now, let $D\subseteq\mathbb{Z}$ be a finite $d$-progression which is large enough such that $\mathcal{G}|_D$ is smooth and has algebraic length~1.
    Then $\Csp(\mathcal{G}|_D)$ is $\Np$-hard \cite{barto2009}.
    Since by Lemma~\ref{lem:pp_defining_gcd_progressions} the structure $(\mathbb{Z},\Succ,\Diff_S)$ pp-defines $\Diff_D$, by Lemma~\ref{lem:use_relation_to_reduce_to_finite_csp} we get $\Csp(\mathcal{G}|_D)\leq_\Ptime\Csp(\mathcal{G},\Diff_D)\leq_\Ptime\Csp(\mathbb{Z},\Succ,\Diff_S)$ and thus $\Csp(\mathbb{Z},\Succ,\Diff_S)$ is $\Np$-hard.
\end{proof}
As mentioned before, the proof of the following lemma given in \cite{bodirsky2016} is incomplete: it fails whenever $a,b,c\in\mathbb{Z}$ are pairwise distinct and $(b-a)/(c-a)$ or $(c-a)/(b-a)$ is an even integer.
We give a complete proof below.
\begin{lemma}[{\cite[see][Lemma~40]{bodirsky2016}}]\label{lem:pairwise_different_distances_are_hard}
    Let $a,b,c\in\mathbb{Z}$ be pairwise different with $|a-b|\neq|a-c|$.
    Then $\Csp(\mathbb{Z},\Succ,\Diff_{\{a,b\}},\Diff_{\{a,c\}})$ is $\mathsf{NP}$-hard.
\end{lemma}
\begin{proof}
    We can assume w.l.o.g.\ that $0=a<b<c$.
    If $2b\neq c$ then $\{a,b,c,b+c\}$ is not a $d$-progression for any $d$ and since $\Diff_{\{a,b,c,b+c\}}$ is pp-definable
    we get hardness by Lemma~\ref{lem:non_progressions_are_hard}.
    Suppose $2b=c$. We claim that $\Csp(\mathbb{Z},\Succ^b,\Diff_{\{a,b\}},\Diff_{\{a,c\}})$ is $\mathsf{NP}$-hard.
    Because of homomorphic equivalence we can assume that $b=1$.
    From Lemmas~\ref{lem:pp_defining_gcd_progressions} and \ref{lem:use_relation_to_reduce_to_finite_csp} it follows that $\Csp(\{0,1,2\},\Succ,\Diff_{\{0,1\}},\Diff_{\{0,2\}})$ reduces to this problem.
    We show that this finite structure does not have a cyclic polymorphism of arity greater than~1,
    which means that it is $\mathsf{NP}$-hard by \cite{barto2012}.
    Suppose for contradiction that $f:\{0,1,2\}^k\to\{0,1,2\}$ is such a polymorphism.
    Since it preserves $\Succ$ it has to be idempotent.
    Now, for all $\ell<k-1$, since $f$ preserves $\Diff_{\{0,1\}}$
    \[f(\underbrace{0,\dots,0}_\ell,\textcolor{red}{1},\textcolor{red}{2},2,\dots,2)-f(\underbrace{0,\dots,0}_\ell,\textcolor{red}{0},\textcolor{red}{1},2,\dots,2)\in\{0,1\}\]
    and because it preserves $\Diff_{\{0,2\}}$
    \[f(\underbrace{0,\dots,0}_\ell,1,2,2,\dots,\textcolor{red}{2})-f(\underbrace{0,\dots,0}_\ell,1,2,\dots,2,\textcolor{red}{0})\in\{0,2\}.\]
    Since $f$ is cyclic the right side of both differences are equal and hence it follows that
    \[f(\underbrace{0,\dots,0}_\ell,1,2,\dots,2)=f(\underbrace{0,\dots,0}_{\ell+1},1,2,\dots,2)\]
    for all $\ell<k-1$.
    In particular, $f(0,\dots,0,1)=f(1,2,\dots,2)$.
    Because $f$ is idempotent and preserves $\Diff_{\{0,1\}}$ we derive $f(0,\dots,0,1)=f(1,2,\dots,2)=1$.
    Since $f$ preserves $\Succ$ we get $f(1,\dots,1,2)=2$ and $f(0,1,\dots,1)=0$.
    However, this contradicts $f$ preserving $\Diff_{\{0,1\}}$.
\end{proof}
\begin{proposition}\label{prop:hardness_unless_all_d_progressions}
    Let $\Diff_S$ and $\Diff_T$ be first-order definable over $(\mathbb{Z},<)$ for some $S,T\subsetneq\mathbb{Z}$ with $|S|>1$.
    If $S$ is not cofinite and there is no $d\in\mathbb{N}$ such that both $\Diff_S$ and $\Diff_T$ are $d$-progressions, then $\Csp(\mathbb{Z},\Succ,\Diff_S,\Diff_T)$ is $\Np$-hard.
\end{proposition}
\begin{proof}
    First, we show that we can assume $S$ to be finite and $T$ to be finite or cofinite.
    Suppose that one of them, let us say $S$, is neither finite nor cofinite.
    By Lemma~\ref{lem:fo_definable_binary_relations} there exists $n\in\mathbb{N}$ such that either $-m\notin S$ and $m\in S$ for all $m\geq n$ or $-m\in S$ and $m\notin S$ for all $m\geq n$.
    Since both cases are symmetric we suppose that the former applies.
    Note, that $\Diff_S(y,x)$ pp-defines $\Diff_{-S}$ and hence $\Diff_U$ for $U\coloneq-S+2n+1$ is pp-definable by Observation~\ref{obs:shift_binary_relations}.
    Now, $S\cap U$ is finite and contains $n$ and $n+1$.
    Since $\Diff_S(x,y)\land\Diff_U(x,y)$ pp-defines $\Diff_{S\cap U}$ it follows by Lemma~\ref{lem:pp_defining_gcd_progressions} that $\Diff_D$ is pp-definable for every finite $1$-progression $D\subseteq\mathbb{Z}$.
    Suppose $S$ is not a $d$-progression for any $d\in\mathbb{N}$.
    Then there exists $m\in\mathbb{N}$ such that $S\cap\{-m,\dots,m\}$ is not a $d$-progression for any $d\in\mathbb{N}$.
    Since $\Diff_S(x,y)\land\Diff_{\{-m,\dots,m\}}(x,y)$ pp-defines $\Diff_{S\cap\{-m,\dots,m\}}$, $\Np$-hardness of $\Csp(\mathbb{Z},\Succ,\Diff_S,\Diff_T)$ follows by Lemma~\ref{lem:non_progressions_are_hard}.
    So suppose $S$ is a $d$-progression for some $d\in\mathbb{N}$.
    Since $S$ is infinite, by Lemma~\ref{lem:fo_definable_binary_relations} it can only be a $1$-progression.
    However, we already showed that we can pp-define $\Diff_D$ for every finite $1$-progression $D\subseteq\mathbb{Z}$ so we may as well use $\{0,1\}$ instead of $S$.

    Since $S$ is finite, if it is not a $d$-progression for any $d\in\mathbb{N}$ we get hardness by Lemma~\ref{lem:non_progressions_are_hard}, so assume $S$ is a $d$-progression.

    Suppose $T$ is finite.
    If $T$ is not a $d'$-progression for any $d'\in\mathbb{N}$, again we get hardness from Lemma~\ref{lem:non_progressions_are_hard}.
    On the other hand, if $T$ is a $d'$-progression, it contains at least two elements because otherwise it is a $d$-progressions.
    Hence, by Lemma~\ref{lem:pp_defining_gcd_progressions} we can pp-define $\Diff_{\{0,d\}}$ and $\Diff_{\{0,d'\}}$, which yields hardness by Lemma~\ref{lem:pairwise_different_distances_are_hard}.

    Suppose now that $T$ is cofinite.
    Let $m\coloneq\min(\mathbb{Z}\setminus T)$ and $n\coloneq\max(\mathbb{Z}\setminus T)$.
    Pick $k\in\mathbb{N}$ minimal such that $m+kd>n$ and let $U\coloneq\{m-2d,m-d,m,m+d,\dots,m+kd\}$.
    Then $m-2d,m-d,m+kd\in T\cap U$ with $k>0$ and $m\notin T\cap U$, which means that $T\cap U$ is not a $d$-progrssion.
    Lemma~\ref{lem:pp_defining_gcd_progressions} yields a pp-definition of $\Diff_U$ and therefore $\Diff_{T\cap U}(x,y)$ is pp-defined by $\Diff_T(x,y)\land\Diff_U(x,y)$.
    Hence, we get hardness by Lemma~\ref{lem:non_progressions_are_hard}.
\end{proof}

\subsection{Hardness from higher arity relations}
Of course, hardness of a CSP can also originate from relations of arity greater than two.
Our next goal is to show that, provided $\Succ$ and a relation of the form $\Diff_{\{0,d\}}$ are present, the absence of $d$-modular max and min polymorphisms already yields hardness for CSPs of first-order reducts of $(\mathbb{Z},\Succ)$.

\textcite{bodirsky2016} invest some work into showing that a relation satisfying the assumptions of the next lemma pp-defines an even nicer relation, from which they then derive hardness.
We instead give a direct and simple proof of $\Np$-hardness.
\begin{definition}
    A function $f:A^k\to A$ with $k\geq 2$ is called a \emph{weak near-unanimity (WNU)} if $f(b,a,\dots,a)=f(a,b,a,\dots,a)=\dots=f(a,\dots,a,b)$ for all $a,b\in A$.
\end{definition}
\begin{lemma}\label{lem:wnu_is_max_or_min}
    Let $\Gamma$ be a first-order reduct of $(\mathbb{Z},<)$ containing $\Succ$, $\Diff_{\{0,1\}}$ and a relation $R\subseteq\mathbb{Z}^k\setminus\{(0,\dots,0)\}$ with $k\geq 3$ such that for all $i\in[k]$ we have
    \[(\underbrace{0,\dots,0}_{i-1},1,0,\dots,0)\in R\text{ or }(\underbrace{0,\dots,0}_{i-1},-1,0,\dots,0)\in R.\]
    Then $\Gamma$ admits $\max$ or $\min$ as a polymorphism or $\Csp(\Gamma)$ is $\mathsf{NP}$-hard.
\end{lemma}
\begin{proof}
    If $(\underbrace{0,\dots,0}_{i-1},1,0,\dots,0)\in R$ for all $i\in[k]$ or $(\underbrace{0,\dots,0}_{i-1},-1,0,\dots,0)\in R$ for all $i\in[k]$, then $\min$, $\max$ and the ternary majority operation do not preserve $R|_{\{0,1\}}$.
    Since $\Diff_{\{0,1\}}$ restricted to $\{0,1\}$ is not preserved by the ternary minority operation, it follows that $\Csp(\Gamma|_{\{0,1\}})$ is $\Np$-hard by Schaefer's Theorem \cite{schaefer1978} and using Lemma~\ref{lem:use_relation_to_reduce_to_finite_csp} we are done.
    Hence, suppose w.l.o.g.\ that
    \[(\delta,0,\dots,0),(0,\delta,0,\dots,0),(0,0,-\delta,0,\dots,0)\in R\]
    for some $\delta\in\{-1,1\}$.

    Assuming that $\Csp(\Gamma)$ is not $\Np$-hard, we show that $\Gamma$, depending on $\delta$, admits $\min$ or $\max$ as a polymorphism.
    Let $n\in\mathbb{N}_{>0}$.
    Then by Lemmas~\ref{lem:use_relation_to_reduce_to_finite_csp} and~\ref{lem:pp_defining_gcd_progressions}, the CSP of $\Gamma|_{\{-2n,\dots,2n\}}$ is not $\Np$-hard.
    Note, that this structure is a finite core.
    Hence, by combining \cite{bulatov2005} and \cite{maroti2008} (see also \cite[Corollary~42]{barto2017}), it admits an $\ell$-ary idempotent WNU polymorphism $f$.

    We assume for the sake of contradiction that $f(\delta,0,\dots,0)=0$ and prove by induction on $i\in[\ell]$ that $f(\underbrace{\delta,\dots,\delta}_i,0,\dots,0)=0$, which holds for $i=0$ by idempotency of $f$.
    For the induction step we use that
    \[\begin{array}{c@{\,}ccccccc@{\,}c@{\;}c@{\;}c}
        f( & \delta & \dots & \delta & 0 & 0 & \dots & 0 & ) & = & 0\\
        f( & 0 & \dots & 0 & \delta & 0 & \dots & 0 & ) & = & 0\\
        f( & 0 & \dots & 0 & 0 & -\delta & \dots & -\delta & ) & \in & \{0,-\delta\}\\
        f( & 0 & \dots & 0 & 0 & 0 & \dots & 0 & ) & = & 0\\
         & \vdots & & \vdots & \vdots & \vdots & & \vdots & & & \vdots\\
        f( & 0 & \dots & 0 & 0 & 0 & \dots & 0 & ) & = & 0
    \end{array}\]
    where the first row contains $i$ non-zero entries,
    implies $f(\underbrace{0,\dots,0}_{i+1},-\delta,\dots,-\delta)=-\delta$ since $(0,\dots,0)\notin R$ and $f$ preserves $R$.
    By preservation of $\Succ$ we get $f(\underbrace{\delta,\dots,\delta}_{i+1},0,\dots,0)=0$.
    Hence, by induction it follows that $f(\delta,\dots,\delta)=0$, which contradicts $f$ being idempotent.
    Thus, $f(\delta,0,\dots,0)\neq 0$ and by idempotency and preservation of $\Diff_{\{0,1\}}$ we get $f(\delta,0,\dots,0)=\delta$.
    With preservation of $\Succ$ it also follows that $f(2\delta,\delta,\dots,\delta)=2\delta$
    and since $\Diff_{\{0,1\}}$ is preserved we get $f(0,\delta,\dots,\delta)=f(\delta,\dots,\delta,0)=\delta$.

    We now prove inductively for all $a,b\in\{-n,\dots,n\}$ that $f(a,b,\dots,b)=\max(a,b)$ if $\delta=1$ and that $f(a,b,\dots,b)=\min(a,b)$ if $\delta=-1$.
    Suppose this holds for $|a-b|\leq t$.
    Then $f(\delta t,0,\dots,0)=\delta t=f(0,\delta t,\dots,\delta t)$ and by preservation of $\Diff_{\{0,1\}}$ we have $f(\delta t+\delta,0,\dots,0),f(0,\delta t+\delta,\dots,\delta t+\delta)\in\{\delta t,\delta t+\delta\}$.
    Now,
    \[\begin{array}{c@{\,}cccc@{\,}c@{\;}c@{\;}c}
        f( & \delta t+\delta & 0 & \dots & 0 & ) & \in & \{\delta t,\delta t+\delta\}\\
        f( & \delta t & \delta & \dots & \delta & ) & = & \delta t\\
        f( & \delta t & 0 & \dots & 0 & ) & = & \delta t\\
         & \vdots & \vdots & & \vdots & & & \vdots\\
        f( & \delta t & 0 & \dots & 0 & ) & = & \delta t
    \end{array}\]
    and
    \[\begin{array}{c@{\,}cccc@{\,}c@{\;}c@{\;}c}
        f( & 0 & \delta t+\delta & \dots & \delta t+\delta & ) & \in & \{\delta t,\delta t+\delta\}\\
        f( & \delta  & \delta t & \dots & \delta t & ) & = & \delta t\\
        f( & 0 & \delta t & \dots & \delta t & ) & = & \delta t\\
         & \vdots & \vdots & & \vdots & & & \vdots\\
        f( & 0 & \delta t & \dots & \delta t & ) & = & \delta t
    \end{array}\]
    imply $f(\delta t+\delta,0,\dots,0)=\delta t+\delta=f(0,\delta t+\delta,\dots,\delta t+\delta)$.
    Since $f$ is defined on $\{-2n,\dots,2n\}$ and preserves $\Succ$, this concludes the proof.
\end{proof}
\begin{lemma}[{\cite[see][Proposition~7.11]{bodirsky2018}}]\label{lem:one_progression_and_no_min_max_is_hard}
    Let $\Gamma$ be a first-order reduct of $(\mathbb{Z},<)$ containing $\Succ$ and $\Diff_{\{0,1\}}$.
    If $\Gamma$ admits neither $\max$ nor $\min$ as a polymorphism then $\Csp(\Gamma)$ is $\Np$-hard.
\end{lemma}
\begin{proof}
    Let $R\subseteq\mathbb{Z}^k$ be a relation of minimal arity that is pp-definable in $\Gamma$ and not preserved by $\max$.
    By Proposition~\ref{prop:hardness_unless_all_d_progressions}, if $\Gamma$ pp-defines $\Diff_T$ for some $T\subseteq\mathbb{Z}$ which is not a $1$-progression, then $\Csp(\Gamma)$ is $\Np$-hard.
    Hence, using Lemma~\ref{lem:fo_definable_binary_relations} we assume that all pp-definable binary relations are of the form $\Diff_T$ for a $1$-progression $T\subseteq\mathbb{Z}$.
    Since all such relations are preserved by $\max$ we get that $k>2$.
    Note, that $S\coloneq\{\bar{a}\in\mathbb{Z}^k\mid(a_i)_{i\in I}\in\Proj_I(R),I\subseteq[k],|I|=k-1\}$ is pp-definable in $\Gamma$.
    Thus, $S$ is preserved by $\max$ iff $\Proj_I(R)$ is preserved by $\max$ for all $I\subseteq[k]$ with $|I|=k-1$.
    Because $R$ has minimal arity it follows that $R\subsetneq S$.
    Hence, there exists $\bar{a}\in S\setminus R$, which means that there are $b_1,\dots,b_k\in\mathbb{Z}$ such that $(a_1,\dots,a_{i-1},b_i,a_{i+1},\dots,a_k)\in R$ for all $i\in[k]$.
    Let $Q\subseteq\mathbb{Z}^k$ be the relation pp-defined by $\exists\bar{y}(R\bar{y}\land\bigwedge_{i\in[k]}\Succ^{a_i}(x_i,y_i))$.
    Then $(0,\dots,0)\notin Q$ but for each $i\in[k]$ we have
    \[(\underbrace{0,\dots,0}_{i-1},b_i-a_i,0,\dots,0)\in Q.\]
    We set $Q_0\coloneq Q$ and for $i\in[k]$ pp-define $Q_i\subseteq\mathbb{Z}^k$ inductively by
    \[\exists y_i(Q_{i-1}x_1\dots x_{i-1}y_ix_{i+1}\dots x_k\land\Diff_{\{-n_i,\dots,n_i\}}(x_i,y_i))\]
    where $n_i\in\mathbb{N}$ is maximal such that
    \[(\underbrace{0,\dots,0}_{i-1},c,0,\dots,0)\notin Q_{i-1}\text{ for all }c\in\{-n_i,\dots,n_i\}.\]
    Note, that $\Diff_{\{-n_i,\dots,n_i\}}$ is pp-definable from $\Succ$ and $\Diff_{\{0,1\}}$ by Lemma~\ref{lem:pp_defining_gcd_progressions}.
    Now, for each $i\in[k]$ we have $(0,\dots,0)\notin Q_i$ and there are $c_1,\dots,c_i\in\{-1,1\}$ such that
    \[(\underbrace{0,\dots,0}_{j-1},c_j,0,\dots,0)\in Q_i\text{ for all }j\in[i].\]
    We get that $Q_k$ satisfies the requirements of Lemma~\ref{lem:wnu_is_max_or_min} and since $\Gamma$ does not admit $\min$ or $\max$ as a polymorphism it follows that $\Csp(\Gamma)$ is $\Np$-hard.
\end{proof}
\begin{lemma}[{\cite[see][Lemma~7.4]{bodirsky2018}}]\label{lem:positive_relations}
    Let $R\subseteq\mathbb{Z}^k$ be first-order definable in $(\mathbb{Z},\Succ)$.
    Then $R$ has a positive quantifier-free definition in $(\mathbb{Z},(\Diff_n)_{n\in\mathbb{Z}})$ or the structure $(\mathbb{Z},\Succ,R)$ pp-defines $\Diff_S$ for some cofinite $S\subsetneq\mathbb{Z}$.
\end{lemma}
\begin{proof}
    By Fact~\ref{fact:quantifier_elimination_succ} the relation $R$ has a quantifier-free definition in $(\mathbb{Z},(\Diff_n)_{n\in\mathbb{Z}})$.
    Let $m\in\mathbb{N}$ be maximal so that every such definition which is in disjunctive normal form contains a clause with at least $m$ negative literals.
    Among the quantifier-free definitions of $R$ that are in disjunctive normal form and contain no more than $m$ negative literals in each clause, we choose $\phi(\bar{x})$ to be one that minimizes the number of clauses containing $m$ negated literals.
    If $m=0$, then $\phi$ is a positive quantifier-free definition of $R$ in $(\mathbb{Z},(\Diff_n)_{n\in\mathbb{Z}})$, so suppose that $m>0$.
    Let $\psi(\bar{x})$ be a clause of $\phi$ containing $m$ negative literals and let $\lnot\Diff_\ell(x_r,x_s)$ be such a literal.
    Further, let $\phi'(\bar{x})$ be the formula $\phi$ without the clause $\psi$ and let $\psi'(\bar{x})$ be the clause $\psi$ without the literal $\lnot\Diff_\ell(x_r,x_s)$.
    If in $\psi'$ the variables $x_r$ and $x_s$ are connected by a chain of positive literals, then $\psi'$ implies $\Diff_n(x_r,x_s)$ for some $n\in\mathbb{Z}$, which means that $\phi'\lor\psi'$ already defines $R$ contradicting the minimality assumptions.
    Hence, there exists $I\subseteq[k]$ containing $r$ but not $s$ such that for every positive literal $\Diff_n(x_i,x_j)$ in $\psi'$ we have $i\in I\iff j\in I$.
    Moreover, there is some $\bar{a}\in\mathbb{Z}^k$ satisfying $\psi'\land\Diff_\ell(x_r,x_s)$ but not $\phi'$ because otherwise $\phi'\lor\psi'$ would again define $R$.

    We now claim that existentially quantifying all variables except for $x_r$ and $x_s$ in
    \[\chi(\bar{x})\coloneq R\bar{x}\land\bigwedge_{i,j\in I}\Succ^{a_j-a_i}(x_i,x_j)\land\bigwedge_{i,j\in[k]\setminus I}\Succ^{a_j-a_i}(x_i,x_j)\]
    yields a pp-definition of $\Diff_S$ for some cofinite $S\subsetneq\mathbb{Z}$.
    By Lemma~\ref{lem:fo_definable_binary_relations} the relation is of the form $\Diff_S$ for some finite or cofinite $S\subseteq\mathbb{Z}$.
    Since $\bar{a}\notin R$ we get that $a_s-a_r\notin S$ and hence $S\subsetneq\mathbb{Z}$.
    To show that $S$ is infinite, fix $n\in\mathbb{N}$ with $n>\max\{|a_i-a_j|+|n'|\mid\lnot\Diff_{n'}(x_i,x_j)\text{ appears in }\psi\}$ and let
    \[b_i\coloneq
    \begin{cases}
        a_i & \text{if }i\in I\\
        a_i+n & \text{if }i\notin I
    \end{cases}\]
    for all $i\in[k]$.
    Since $\bar{a}$ satisfies $\psi'$, $\bar{b}$ satisfies every literal of $\psi$ containing only variables $x_i$ and $x_j$ with $i\in I\iff j\in I$.
    By our choice of $I$, every other literal of $\psi$ is of the form $\lnot\Diff_{n'}(x_i,x_j)$.
    Because we picked $n$ such that $|n'|<n-|a_i-a_j|\leq|b_i-b_j|\in\{|a_i-a_j\pm n|\}$, $\bar{b}$ satisfies all these literals and hence $\psi$.
    It follows that $\bar{b}$ satisfies $\chi$ and because we have infinitely many choices for $n$ the set $S$ must be infinite.
\end{proof}
\textcite[Proposition~7.12]{bodirsky2018} establish hardness of the CSP for positive first-order expansions of $(\mathbb{Z},\Succ)$ that admit neither $\Min_d$ nor $\Max_d$ as a polymorphism.
Their proof, however, assumes that every binary projection has the form $\Diff_S$ for some finite $S\subseteq\mathbb{Z}$, which does not hold in general, since a binary projection may also be the full relation $\mathbb{Z}^2$.
We resolve this issue and thereby close the gap.
\begin{proposition}\label{prop:d_progression_and_no_mod_min_max_is_hard}
    Let $\Gamma$ be a first-order reduct of $(\mathbb{Z},\Succ)$ containing $\Succ$ and $\Diff_{\{0,d\}}$ for some $d\in\mathbb{N}_{>0}$.
    If $\Gamma$ admits neither $\Max_d$ nor $\Min_d$ as a polymorphism then $\Csp(\Gamma)$ is $\Np$-hard.
\end{proposition}
\begin{proof}
    If $d=1$ we are done by Lemma~\ref{lem:one_progression_and_no_min_max_is_hard}, so suppose that $d>1$.
    We can assume that all relations of $\Gamma$ have a positive quantifier-free definition in $(\mathbb{Z},(\Diff_n)_{n\in\mathbb{Z}})$ because otherwise by Lemma~\ref{lem:positive_relations} for some cofinite $S\subsetneq\mathbb{Z}$ the relation $\Diff_S$ is pp-definable in $\Gamma$, which means that $\Csp(\Gamma)$ is $\Np$-hard by Proposition~\ref{prop:hardness_unless_all_d_progressions}.
    It can easily be verified that $\Diff_S$ is preserved by $\Max_d$ and $\Min_d$ for every $d$-progression $S\subseteq\mathbb{Z}$ and for $S=\mathbb{Z}$.
    Thus, if $\Gamma$ pp-defines a binary relation that is not preserved by $\Min_d$ or $\Max_d$, then by Lemma~\ref{lem:fo_definable_binary_relations} this relation is of the form $\Diff_S$ where $S\subsetneq\mathbb{Z}$ is not a $d$-progression and again $\Csp(\Gamma)$ is $\Np$-hard by Proposition~\ref{prop:hardness_unless_all_d_progressions}.
    Hence, let us assume that all pp-definable binary relations are preserved by $\Min_d$ and $\Max_d$.

    Let $Q\subseteq\mathbb{Z}^k$ be of minimal arity such that it is pp-definable in $\Gamma$ and not preserved by $\Max_d$.
    We claim that there are $\bar{a},\bar{b}\in Q$ with $\Max_d(\bar{a},\bar{b})\notin Q$ such that $a_i\equiv b_i\mod{d}$ for all $i\in[k]$.
    For $i,j\in[k]$ let $i\sim j$ if $a_i-a_j\equiv b_i-b_j\mod{d}$ for all $\bar{a},\bar{b}\in Q$.
    Suppose ${\sim}=[k]^2$ and let $\bar{a},\bar{b}\in Q$ with $\Max_d(\bar{a},\bar{b})\notin Q$.
    Then $a_i\equiv b_i\mod{d}$ for all $i\in[k]$ because otherwise $a_i\not\equiv b_i\mod{d}$ for all $i\in[k]$ and $\Max_d(\bar{a},\bar{b})=\bar{a}\in Q$.
    So assume that ${\sim}\neq[k]^2$.
    Since $Q$ is not preserved by $\Max_d$, there is some $\bar{c}\in\Max_d(Q\times Q)\setminus Q$.
    Because $Q$ has minimal arity among the pp-definable relations that are not preserved by $\Max_d$, for every $\ell\in[k]$ the relation $\Proj_{[k]\setminus\{\ell\}}(Q)$ is preserved by $\Max_d$.
    With Lemma~\ref{lem:shifting_modular_blocks} it follows that for all $\ell\in[k]$ and $n_{[1]_\sim},\dots,n_{[k]_\sim}\in\mathbb{Z}$ there is some $c'\in\mathbb{Z}$ such that
    \[(c_1+n_{[1]_\sim},\dots,c_{\ell-1}+n_{[\ell-1]_\sim},c',c_{\ell+1}+n_{[\ell+1]_\sim},\dots,c_k+n_{[k]_\sim})\in Q.\]
    Let $\phi(\bar{x})$ be a positive quantifier-free definition of $Q$ in $(\mathbb{Z},(\Diff_n)_{n\in\mathbb{Z}})$ which is in disjunctive normal form and note that if $i\sim j$ then $x_i$ and $x_j$ are connected by a chain of (positive) literals in every clause of $\phi$.
    It follows that if $\ell\in[k]$ with $[\ell]_\sim\neq\{\ell\}$, then $\phi$ contains a clause that is equivalent to
    \[\psi_\ell(\bar{x})\coloneq\bigwedge_{i\sim j\not\sim\ell}\Diff_{c_j-c_i}(x_i,x_j)\land\bigwedge_{\substack{i\sim\ell\\i\neq\ell}}\Diff_{c_i-c_\ell'}(x_\ell,x_i)\]
    for some $c_\ell'\in\mathbb{Z}$.
    We prove the claim seperately for the following three cases.
    \begin{description}
        \item[Case~1:] There are $r,s\in[k]$ with $[r]_\sim\neq\{r\}$, $[s]_\sim\neq\{s\}$ and $r\not\sim s$.

            Take any $n_{[r]_\sim}\in d\mathbb{Z}_{\leq 0}$ with $n_{[r]_\sim}\leq c_r-c_r'$ and set $n_{[i]_\sim}\coloneq 0$ for all $i\not\sim r$.
            Similarly, let $m_{[s]_\sim}\in d\mathbb{Z}_{\leq 0}$ with $m_{[s]_\sim}\leq c_s-c_s'$ and $m_{[i]_\sim}\coloneq 0$ for all $i\not\sim s$.
            Then
            \[\bar{a}\coloneq(c_1+n_{[1]_\sim},\dots,c_{r-1}+n_{[r-1]_\sim},c_r'+n_{[r]_\sim},c_{r+1}+n_{[r+1]_\sim},\dots,c_k+n_{[k]_\sim})\]
            satisfies $\psi_r$ and
            \[\bar{b}\coloneq(c_1+m_{[1]_\sim},\dots,c_{s-1}+m_{[s-1]_\sim},c_s'+m_{[s]_\sim},c_{s+1}+m_{[s+1]_\sim},\dots,c_k+m_{[k]_\sim})\]
            satisfies $\psi_s$, which means that $\bar{a},\bar{b}\in Q$.
            Note, that for all $i\notin[r]_\sim\cup[s]_\sim$ we have $a_i=b_i=c_i$.
            For $i\in[r]_\sim\setminus\{r\}$ we get $a_i-b_i=(c_i+n_{[i]_\sim})-c_i=n_{[r]_\sim}\in d\mathbb{Z}$ and thus $a_i\equiv b_i\mod{d}$.
            Similarly, $a_i\equiv b_i\mod{d}$ for all $i\in[s]_\sim\setminus\{s\}$.
            Since $[r]_\sim\neq\{r\}$ there exists some $i\in[r]_\sim\setminus\{r\}$.
            Hence, $a_r-a_i\equiv b_r-b_i\mod{d}$ and with $a_i\equiv b_i\mod{d}$ it follows that $a_r\equiv b_r\mod{d}$.
            The same argument shows that $a_s\equiv b_s\mod{d}$ and therefore $a_i\equiv b_i\mod{d}$ for all $i\in[k]$.
            To prove the claim, it now suffices to verify that $\max(a_i,b_i)=c_i$ for all $i\in[k]$.
            Trivially, this holds for all $i\notin[r]_\sim\cup[s]_\sim$.
            For $i\in[r]_\sim\setminus\{r\}$ we have $\max(a_i,b_i)=\max(c_i+n_{[r]_\sim},c_i)=c_i$ because $n_{[r]_\sim}\leq 0$.
            Similarly, we get $\max(a_i,b_i)=c_i$ for all $i\in[s]_\sim\setminus\{s\}$.
            Finally, $\max(a_r,b_r)=\max(c_r'+n_{[r]_\sim},c_r)=c_r$ and $\max(a_s,b_s)=\max(c_s,c_s'+m_{[s]_\sim})=c_s$ by our choice of $n_{[r]_\sim}$ and $m_{[s]_\sim}$.
        \item[Case~2:] There are $r,s\in[k]$ with $[r]_\sim\neq\{r\}$ and $[s]_\sim=\{s\}$.

            Let $c_s'\in\mathbb{Z}$ such that $\bar{a}\coloneq(c_1,\dots,c_{s-1},c_s',c_{s+1},\dots,c_k)\in Q$.
            Fix $n_{[r]_\sim}\in d\mathbb{Z}_{\leq 0}$ with $n_{[r]_\sim}\leq c_r-c_r'$ and let $n_{[i]_\sim}\coloneq 0$ for all $i\notin[r]_\sim\cup[s]_\sim$.
            Then for any $n_{[s]_\sim}\in\mathbb{Z}$ the tuple
            \[\bar{b}\coloneq(c_1+n_{[1]_\sim},\dots,c_{r-1}+n_{[r-1]_\sim},c_r'+n_{[r]_\sim},c_{r+1}+n_{[r+1]_\sim},\dots,c_k+n_{[k]_\sim})\]
            satisfies $\psi_r$ and hence $\bar{b}\in Q$.
            Note, that for all $i\notin[r]_\sim\cup[s]_\sim$ we have $a_i=b_i=c_i$.
            For all $i\in[r]_\sim\setminus\{r\}$ we get $b_i-a_i=(c_i+n_{[i]_\sim})-c_i=n_{[r]_\sim}\in d\mathbb{Z}$ and thus $a_i\equiv b_i\mod{d}$.
            Since there exists such an $i\in[r]_\sim\setminus\{r\}$ we have $a_r-a_i\equiv b_r-b_i\mod{d}$ and it follows that $a_r\equiv b_r\mod{d}$.
            Now, choosing $n_{[s]_\sim}$ such that $n_{[s]_\sim}\equiv c_s'-c_s\mod{d}$ yields $a_i\equiv b_i\mod{d}$ for all $i\in[k]$.
            If additionally $n_{[s]_\sim}\geq c_s'-c_s$ then $\Max_d(\bar{a},\bar{b})=\max(\bar{a},\bar{b})=(c_1,\dots,c_{s-1},c_s+n_{[s]_\sim},c_{s+1},\dots,c_k)$.
            Suppose $(c_1,\dots,c_{s-1},c_s+n_{[s]_\sim},c_{s+1},\dots,c_k)\in Q$ for all $n_{[s]_\sim}\geq c_s'-c_s$ with $n_{[s]_\sim}\equiv c_s'-c_s\mod{d}$.
            Then, because $\phi$ is positive, we have $(c_1,\dots,c_{s-1},c',c_{s+1},\dots,c_k)\in Q$ for all $c'\in\mathbb{Z}$, which contradicts $\bar{c}\notin Q$.
            Hence, we can choose $n_{[s]_\sim}$ in a way such that $a_i\equiv b_i\mod{d}$ for all $i\in[k]$ and $\max(\bar{a},\bar{b})\notin Q$.
        \item[Case~3:] $[r]_\sim=\{r\}$ for all $r\in[k]$.

            We have $\Proj_{[k]\setminus\{r\}}(Q)=\mathbb{Z}^{k-1}$.
            Hence, $\phi$ contains a clause which has no chain of literals connecting $x_i$ and $x_j$ for any distinct $i,j\in[k]\setminus\{r\}$.
            Since $Q\neq\mathbb{Z}^k$, this clause must be equivalent to $\Diff_{n_r}(x_r,x_{r'})$ for some $n_r\in\mathbb{Z}$ and $r'\in[k]\setminus\{r\}$.
            We also get such a clause equivalent to $\Diff_{n_s}(x_s,x_{s'})$ for $s\in[k]\setminus\{r,r'\}$.
            If $r=s'$ we swap $r$ and $r'$, so that out of $r,r',s,s'$ only $r'$ and $s'$ can be equal.
            Let $m\in\mathbb{N}$ with $m>d$ such that for every literal $\Diff_{m'}(x_i,x_j)$ in $\phi$ we have $|m'|<m$.
            Writing $\{i_1,\dots,i_{k-2}\}\coloneq[k]\setminus\{r,s\}$, let $R\subseteq\mathbb{Z}^2$ be the relation pp-defined by existentially quantifying all variables but $x_r$ and $x_{i_1}$ in
            $Q\bar{x}\land\Succ^{km}(x_r,x_s)\land\bigwedge_{j\in\{1,\dots,k-3\}}\Succ^m(x_{i_j},x_{i_{j+1}})$.
            \begin{center}
                \begin{tikzcd}[sep=huge]
                    x_r \arrow[rrrrr, bend right=15, "km"] & x_{i_1} \arrow[r, "m"] & \exists x_{i_2} \arrow[r, "m"] & \cdots \arrow[r, "m"] & \exists x_{i_{k-2}} & \exists x_s
                \end{tikzcd}\
            \end{center}
            By Lemma~\ref{lem:fo_definable_binary_relations}, $R=\Diff_S$ for some $S\subseteq\mathbb{Z}$.
            Since $\phi$ contains a clause equivalent to $\Diff_{n_r}(x_r,x_{r'})$, there is some $t\in S$ with $t<m$.
            On the other hand, because $\phi$ contains a clause equivalent to $\Diff_{n_s}(x_s,x_{s'})$, there exists $t\in S$ with $t>km-(k-2)m=2m$.
            Suppose $S\cap\{m,m+1\dots,2m\}\neq\emptyset$.
            Then $Q$ contains a tuple with $|a_i-a_j|\geq m$ for all distinct $i,j\in[k]$.
            By our choice of $m$, this tuple cannot satisfy a literal in $\phi$, which means that $Q=\mathbb{Z}^k$, a contradiction.
            Thus, $S\cap\{m,m+1,\dots,2m\}=\emptyset$ and because $m>d$ the set $S$ cannot be a $d$-progression.
            It is now easy to verify that $\Diff_S$ is not preserved by $\Max_d$ contradicting our assumption that all pp-definable binary relations are preserved.
    \end{description}

    We have now proved the claim that $Q$ contains tuples $\bar{a}$ and $\bar{b}$ with $\Max_d(\bar{a},\bar{b})\notin Q$ such that $a_i\equiv b_i\mod{d}$ for all $i\in[k]$.
    Let $D\subseteq d\mathbb{Z}$ be a finite $d$-progression containing $b_i-a_i-b_1+a_1$ for all $i\in[k]$.
    Note, that $\Diff_D$ is pp-definable in $\Gamma$ by Lemma~\ref{lem:pp_defining_gcd_progressions} and let $R\subseteq\mathbb{Z}^k$ be the relation pp-defined by
    \[\exists\bar{y}(Q\bar{y}\land\bigwedge_{i\in[k]}\Succ^{a_i}(x_i,y_i))\land\bigwedge_{i\in\{2,\dots,k\}}\Diff_D(x_1,x_i).\]
    Then $c_i\equiv c_j\mod{d}$ for all $\bar{c}\in R$ and $i,j\in[k]$.
    Moreover, $R|_{d\mathbb{Z}}$ contains $\bar{a}'\coloneq(0,\dots,0)$ and $\bar{b}'\coloneq(b_i-a_i)_{i\in[k]}$ but not $\max(\bar{a}',\bar{b}')$.
    An analogous argument shows that $\Gamma$ pp-defines a relation $T\subseteq\mathbb{Z}^\ell$ with $c_i\equiv c_j\mod{d}$ for all $\bar{c}\in T$ and $i,j\in[\ell]$ such that $T|_{d\mathbb{Z}}$ is not preserved by $\min$.
    It is now easy to verify that $\Delta\coloneq(\mathbb{Z},\Succ^d,\Diff_{\{0,d\}},R,T)$ is homomorphically equivalent to $\Delta|_{d\mathbb{Z}}$, which by Lemma~\ref{lem:substructure_isomorphic_to_fo_reduct} is isomorphic to a first-order reduct $(\mathbb{Z},\Succ,\Diff_{\{0,1\}},R',T')$ of $(\mathbb{Z},<)$.
    Hence, $\Csp(\mathbb{Z},\Succ,\Diff_{\{0,1\}},R',T')$ reduces to $\Csp(\Gamma)$ and is $\Np$-hard by Lemma~\ref{lem:one_progression_and_no_min_max_is_hard}.
\end{proof}

\subsection{Non-Horn-definable relations}
Earlier, we established tractability for CSPs of first-order reducts of $(\mathbb{Z},<)$ that admit a modular max or min polymorphism, or whose relations are Horn-definable in $(\mathbb{Z},\Succ)$.
We now complete the hardness side of the classification for first-order expansions of $(\mathbb{Z},\Succ)$: using the results established so far in this section, we show that the CSP of any such expansion is $\Np$-hard unless it is tractable for one of these two reasons.
\begin{lemma}\label{lem:quantifier_elimination_rank}
    Let $R\subseteq(\mathbb{Q}\times\mathbb{Z})^k$ be first-order definable over $(\mathbb{Q}\times\mathbb{Z},\Succ)$.
    Then there exists a finite $S\subseteq\mathbb{Z}$ such that for all $\bar{a},\bar{b}\in(\mathbb{Q}\times\mathbb{Z})^k$ that satisfy
    \[a_i-a_j=n\iff b_i-b_j=n\]
    for all $n\in S$, we have $\bar{a}\in R\iff\bar{b}\in R$.
\end{lemma}
\begin{proof}
    By Fact~\ref{fact:quantifier_elimination_succ} there exists a quantifier-free first-order formula over the language $(\Diff_n)_{n\in\mathbb{Z}}$ that defines $R$.
    Let $S$ be the set of constants $n$ appearing in this formula.
    Then the claim follows immediately.
\end{proof}
\begin{lemma}[{\cite[Proposition~6.1]{bodirsky2018}}]\label{lem:characterization_of_horn_relations}
    Let $R\subseteq(\mathbb{Q}\times\mathbb{Z})^k$ be first-order definable over $(\mathbb{Q}\times\mathbb{Z},\Succ)$.
    Then $R$ is Horn-definable iff it is preserved by any (equivalently all) binary injective polymorphism(s) of $(\mathbb{Q}\times\mathbb{Z},\Succ)$.
\end{lemma}
\begin{proof}
    Let $\mathcal{M}\coloneq(\mathbb{Q}\times\mathbb{Z},(\Diff_n)_{n\in\mathbb{Z}})$.
    By Fact~\ref{fact:quantifier_elimination_succ} the relation $R$ has a quantifier-free first-order definition over $\mathcal{M}$.
    It is not hard to see that there exists some embedding $\mathcal{M}^2\to\mathcal{M}$.
    By \cite[Proposition~5.9]{bodirsky2012} it follows that $R$ is Horn-definable over $\mathcal{M}$ iff it is preserved by any (equivalently all) embedding(s) $\mathcal{M}^2\to\mathcal{M}$.
    Therefore, it suffices to show that the emeddings $\mathcal{M}^2\to\mathcal{M}$ are exactly the binary injective polymorphisms of $(\mathbb{Q}\times\mathbb{Z},\Succ)$.
    Since the other direction is trivial, we prove that every binary injective polymorphism $f$ of $(\mathbb{Q}\times\mathbb{Z},\Succ)$ is an embedding $\mathcal{M}^2\to\mathcal{M}$.
    First, we observe that $f$ is a homomorphism $\mathcal{M}^2\to\mathcal{M}$ because for every $n\in\mathbb{Z}$ the relation $\Diff_n$ is pp-definable by $\Succ^n(x,y)$ over $(\mathbb{Q}\times\mathbb{Z},\Succ)$.
    Now, let $(a,b),(c,d)\in(\mathbb{Q}\times\mathbb{Z})^2$ and $n\in\mathbb{Z}$ such that $((a,b),(c,d))\notin\Diff_n^{\mathcal{M}^2}$.
    Then $a+n\neq c$ or $b+n\neq d$ and since $f$ is an injective homomorphism $\mathcal{M}^2\to\mathcal{M}$ preserving $\Succ$, we get $f(a,b)+n=f(a+n,b+n)\neq f(c,d)$.
    Hence, $f$ is an embedding.
\end{proof}
\begin{lemma}\label{lem:characterization_of_non_horn_relations}
    Let $R\subseteq(\mathbb{Q}\times\mathbb{Z})^k$ be first-order definable over $(\mathbb{Q}\times\mathbb{Z},\Succ)$.
    Then $R$ is not Horn-definable if and only if $(\mathbb{Q}\times\mathbb{Z},\Succ,R)$ pp-defines a non-empty relation $Q\subseteq(\mathbb{Q}\times\mathbb{Z})^\ell$ with $\ell\leq k$ such that
    \begin{enumerate}
        \item\label{item:characterization_of_non_horn_relations_finite_distance} there exists $n\in\mathbb{N}$ so that for all $\bar{a}\in Q$ there are distinct $i,j\in[\ell]$ with $|a_i-a_j|\leq n$, and
        \item\label{item:characterization_of_non_horn_relations_different_distances} for all distinct $i,j\in[\ell]$ there are $\bar{a},\bar{b}\in Q$ with $a_i-a_j\neq b_i-b_j$ or $|a_i-a_j|=|b_i-b_j|=\infty$.
    \end{enumerate}
\end{lemma}
\begin{proof}
    Suppose for contradiction that $R$ is Horn-definable and satisfies Items~\ref{item:characterization_of_non_horn_relations_finite_distance} and \ref{item:characterization_of_non_horn_relations_different_distances}.
    Then by Lemma~\ref{lem:characterization_of_horn_relations} the structure $(\mathbb{Q}\times\mathbb{Z},\Succ,R)$ admits a binary injective polymorphism $f$, which preserves $Q$.
    Note, that for all $\bar{a},\bar{b}\in(\mathbb{Q}\times\mathbb{Z})^\ell$ and $i,j\in[\ell]$ we have
    \[a_i-a_j=b_i-b_j\in\mathbb{Z}\iff f(a_i,b_i)-f(a_j,b_j)\in\mathbb{Z}.\]
    Let $\bar{a}\in Q$ and $I\coloneq\{(i,j)\in[\ell]^2\mid|a_i-a_j|=\infty\}$.
    Suppose there are distinct $i,j\in[\ell]$ with $(i,j)\notin I$.
    By Item~\ref{item:characterization_of_non_horn_relations_different_distances} there exists $\bar{b}\in Q$ such that $a_i-a_j\neq b_i-b_j$.
    Now, for $\bar{c}\coloneq f(\bar{a},\bar{b})\in Q$ we have $|c_i-c_j|=\infty$ and $|c_{i'}-c_{j'}|=\infty$ for all $(i',j')\in I$.
    By induction it follows that there exists $\bar{d}\in Q$ with $|d_i-d_j|=\infty$ for all distinct $i,j\in[\ell]$.
    However, this contradicts Item~\ref{item:characterization_of_non_horn_relations_finite_distance}.

    For the other direction, fix a binary injective polymorphism $f$ of $(\mathbb{Q}\times\mathbb{Z},\Succ)$ and suppose that $R$ is not Horn-definable.
    Then $f$ does not preserve $R$ by Lemma~\ref{lem:characterization_of_horn_relations}.
    Suppose there are $\bar{a},\bar{b}\in R$ and distinct $i,j\in[k]$ such that $m\coloneq a_i-a_j=b_i-b_j\in\mathbb{Z}$ and $\bar{c}\coloneq f(\bar{a},\bar{b})\notin R$.
    Since $f$ preserves $\Succ$ we have $c_i-c_j=m$ and thus the relation of arity $k-1$ that is pp-defined by $\exists x_j(R\bar{x}\land\Succ^m(x_j,x_i))$ is not preserved by $f$.
    By induction we get a pp-definable relation $Q$ of arity $\ell\leq k$ that contains $\bar{a},\bar{b}$ with $\bar{c}\coloneq f(\bar{a},\bar{b})\notin Q$ and such that for all distinct $i,j\in[\ell]$ we have $a_i-a_j\neq b_i-b_j$ or $|a_i-a_j|=|b_i-b_j|=\infty$.
    In particular, $Q$ satisfies Item~\ref{item:characterization_of_non_horn_relations_different_distances}.
    Since $|c_i-c_j|=\infty$ for all distinct $i,j\in[\ell]$ and $\bar{c}\notin Q$, by Lemma~\ref{lem:quantifier_elimination_rank} it follows that $Q$ satisfies Item~\ref{item:characterization_of_non_horn_relations_finite_distance} as well.
\end{proof}
\textcite{bodirsky2018} give a syntactic proof of the following lemma.
We take a different route, using the characterization of non-Horn-definable relations from Lemma~\ref{lem:characterization_of_non_horn_relations}.
\begin{lemma}[{\cite[Proposition~7.7]{bodirsky2018}}]\label{lem:pp_defining_binary_non_horn_relation}
    Let $R\subseteq\mathbb{Z}^k$ be first-order definable over $(\mathbb{Z},\Succ)$.
    If $R$ is not Horn-definable then $(\mathbb{Z},\Succ,R)$ pp-defines a binary relation that is not Horn-definable.
\end{lemma}
\begin{proof}
    Let $Q\subseteq(\mathbb{Q}\times\mathbb{Z})^\ell$ be a relation that is pp-definable over $\mathbb{Q}.(\mathbb{Z},\Succ,R)$ and satisfies Items~\ref{item:characterization_of_non_horn_relations_finite_distance} and~\ref{item:characterization_of_non_horn_relations_different_distances} from Lemma~\ref{lem:characterization_of_non_horn_relations}.
    Suppose $Q$ has minimal arity among all such relations.
    If $\ell=2$ we are done so suppose $\ell>2$.

    Let $m\in\mathbb{N}$ with $m>n$, where $n$ is from Item~\ref{item:characterization_of_non_horn_relations_finite_distance} of Lemma~\ref{lem:characterization_of_non_horn_relations}, and pick distinct $r,s\in[\ell]$.
    We define $Q_{r,s}^m\subseteq(\mathbb{Q}\times\mathbb{Z})^\ell$ to be the relation pp-defined by $Q\bar{x}\land\Succ^m(x_r,x_s)$.
    Suppose $Q_{r,s}^m$ is non-empty and let $I\subseteq[\ell]$ be maximal such that $\Proj_I(Q_{r,s}^m)$ satisfies Item~\ref{item:characterization_of_non_horn_relations_different_distances} from Lemma~\ref{lem:characterization_of_non_horn_relations}.
    Clearly, $\emptyset\subsetneq I\subsetneq[\ell]$.
    Note, that for every $i\in[\ell]$ there is a unique $i'\in I$ such that for all $\bar{a},\bar{b}\in Q_{r,s}^m$ we have $a_i-a_{i'}=b_i-b_{i'}\in\mathbb{Z}$.
    We claim that there are $\bar{a}\in Q_{r,s}^m$ and distinct $i,j\in[\ell]$ with $|a_i-a_j|\leq n$ and $i'=j'$.
    Suppose otherwise.
    Let $t\coloneq\max\{|a_i-a_{i'}|\mid i\in[\ell],\bar{a}\in Q_{r,s}^m\}$.
    Since $Q\supseteq Q_{r,s}^m$ satisfies Item~\ref{item:characterization_of_non_horn_relations_finite_distance} from Lemma~\ref{lem:characterization_of_non_horn_relations}, for any $\bar{a}\in Q_{r,s}^m$ there are distinct $i,j\in[\ell]$ with $|a_i-a_j|\leq n$.
    Now, $|a_{i'}-a_{j'}|\leq|a_{i'}-a_i|+|a_i-a_j|+|a_j-a_{j'}|\leq t+n+t=n+2t$ and with $i',j'\in I$ being distinct it follows that $\Proj_I(Q_{r,s}^m)$ also satisfies Item~\ref{item:characterization_of_non_horn_relations_finite_distance} of Lemma~\ref{lem:characterization_of_non_horn_relations}.
    This contradicts the minimality assumption on the arity of $Q$ and thus proves the claim.
    We now have some $\bar{a}\in Q_{r,s}^m$ and distinct $i,j\in[\ell]$ with $|a_i-a_j|\leq n$ and $i'=j'$.
    By the definition of $i'$ and $j'$, it follows that $a_i-a_j=a_i-a_{i'}+a_{j'}-a_j=b_i-b_{i'}+b_{j'}-b_j=b_i-b_j$ for all $\bar{b}\in Q_{r,s}^m$.
    Hence, for all $m>n$ and distinct $r,s\in[\ell]$ there are distinct $i,j\in[\ell]$ such that
    \[a_s-a_r=m\implies|a_i-a_j|\leq n\]
    for all $\bar{a}\in Q$.

    We claim that for all distinct $r,s\in[\ell]$ there are distinct $i,j\in[\ell]$ such that
    \[|a_r-a_s|=\infty\implies|a_i-a_j|\leq n\]
    for all $\bar{a}\in Q$.
    Otherwise, there are distinct $r,s\in[\ell]$ such that for all distinct $i,j\in[\ell]$ there exists $\bar{a}\in Q$ with $|a_r-a_s|=\infty$ and $|a_i-a_j|>n$.
    By Lemma~\ref{lem:quantifier_elimination_rank} there exists $m>n$ such that for all distinct $i,j\in[\ell]$ there is some $\bar{a}\in Q$ with $a_s-a_r=m$ and $|a_i-a_j|>n$.
    However, this contradicts what we showed before and thus the claim holds.

    By the minimality assumption on the arity of $Q$, the relation $\Proj_{[\ell-1]}(Q)$ violates Item~\ref{item:characterization_of_non_horn_relations_finite_distance} or~\ref{item:characterization_of_non_horn_relations_different_distances} from Lemma~\ref{lem:characterization_of_non_horn_relations}.
    Since it cannot violate Item~\ref{item:characterization_of_non_horn_relations_different_distances}, for every $m\in\mathbb{N}$ there is some $\bar{a}\in Q$ with $|a_r-a_s|>m$ for all distinct $r,s\in[\ell-1]$.
    By Lemma~\ref{lem:quantifier_elimination_rank} (or a simple compactness argument) there exists $\bar{a}\in Q$ such that $|a_r-a_s|=\infty$ for all distinct $r,s\in[\ell-1]$.
    Since $Q$ satisfies Item~\ref{item:characterization_of_non_horn_relations_finite_distance} of Lemma~\ref{lem:characterization_of_non_horn_relations}, there is some $r\in[\ell-1]$ with $|a_r-a_\ell|\leq n$.
    Let $s\in[\ell-1]\setminus\{r\}$.
    Then $|a_i-a_j|=\infty>n$ for all distinct $i,j\in[\ell]$ with $\{i,j\}\neq\{r,\ell\}$.
    Hence, for all $\bar{b}\in Q$ we have $|b_r-b_s|=\infty\implies|b_r-b_\ell|\leq n$ and $|b_s-b_\ell|=\infty\implies|b_r-b_\ell|\leq n$.
    Now, by the same argument as before there is some $\bar{b}\in Q$ such that $|b_i-b_j|=\infty$ for all distinct $i,j\in[\ell]\setminus\{s\}$.
    In particular, $|b_r-b_\ell|=\infty$ and thus $|b_r-b_s|=\infty$ or $|b_s-b_\ell|=\infty$, contradicting at least one of the implications we just proved.
    Therefore, $Q$ cannot have arity greater than 2, which concludes the proof.
\end{proof}
\begin{proposition}\label{prop:non_horn_relations}
    Let $\Gamma$ be a first-order reduct of $(\mathbb{Z},\Succ)$ containing $\Succ$.
    If $\Gamma$ contains a non-Horn-definable relation and for all $d\in\mathbb{N}_{>0}$ neither admits $\max_d$ nor $\min_d$ as a polymorphism, then $\Csp(\Gamma)$ is $\Np$-hard.
\end{proposition}
\begin{proof}
    Since $\Gamma$ contains a non-Horn-definable relation, by Lemma~\ref{lem:pp_defining_binary_non_horn_relation} it pp-defines a binary relation that is not Horn-definable.
    By Lemma~\ref{lem:fo_definable_binary_relations} this relation is of the form $\Diff_S$ for some finite or cofinite $S\subseteq\mathbb{Z}$.
    Since for cofinite $S$ the relation $\Diff_S$ is defined by the Horn-formula $\bigwedge_{n\in\mathbb{Z}\setminus S}\lnot\Diff_n(x,y)$, $S$ must be finite.
    If $S=\{n\}$ for some $n\in\mathbb{Z}$ then the Horn-formula $\Diff_n(x,y)$ defines $\Diff_S$.
    Hence, $S$ is finite with $|S|>1$.
    Let $d\coloneq\gcd\{s-t\mid s,t\in S\}$.
    Then by Lemma~\ref{lem:pp_defining_gcd_progressions} we can pp-define $\Diff_{\{0,d\}}$.
    With Proposition~\ref{prop:d_progression_and_no_mod_min_max_is_hard} it follows that $\Csp(\Gamma)$ is $\Np$-hard.
\end{proof}

\section{PP-Defining the order}\label{sec:defining_order}
Proposition~\ref{prop:non_horn_relations} settles the last remaining case of the complexity classification of the CSPs of first-order reducts of $(\mathbb{Z},\Succ)$.
Of course, not every first-order reduct of $(\mathbb{Z},<)$ is a first-order reduct of $(\mathbb{Z},\Succ)$.
\textcite[Lemma~5.37]{bodirsky2018} claim, and give a proof, that every first-order reduct of $(\mathbb{Z},<)$ that pp-defines $\Succ$ but is not a first-order reduct of $(\mathbb{Z},\Succ)$ also pp-defines $<$, and use this to obtain a full classification.
As the following example shows, this is false.
\begin{example}\label{exmp:cyclic_order}
    Let $\Gamma\coloneq(\mathbb{Z},\Succ,\Cycl)$ where ${\Cycl}\subseteq\mathbb{Z}^3$ is defined by $x<y<z\lor y<z<x\lor z<x<y$.
    It is easy to find an automorphism of $(\mathbb{Q}\times\mathbb{Z},\Succ)$ which is not an automorphism of $\mathbb{Q}.\Gamma$ and hence $\Gamma$ is not a first-order reduct of $(\mathbb{Z},\Succ)$.
    However, for any $r\in\mathbb{R}\setminus\mathbb{Q}$ and order preserving maps $f:(-\infty,r)\to(r,\infty)$ and $g:(r,\infty)\to(-\infty,r)$ the mapping
    \[\mathbb{Q}\times\mathbb{Z}\to\mathbb{Q}\times\mathbb{Z},(p,a)\mapsto
    \begin{cases}
        (f(p),a) & \text{if }p\in(-\infty,r)\\
        (g(p),a) & \text{if }p\in(r,\infty)
    \end{cases}\]
    is an endomorphism of $\mathbb{Q}.\Gamma$, which does not preserve $<$.
    Therefore, $<$ is not pp-definable in $\Gamma$.
\end{example}
To give a correct proof of the full complexity classification, we spend this section investigating the \enquote{global} behavior of endomorphisms of first-order reducts of $(\mathbb{Q}\times\mathbb{Z},<)$, and use a result by \textcite{bodirsky2010} to show that every counterexample to the claim of \textcite{bodirsky2018} has an $\Np$-hard CSP.
\begin{remark}
    Lemma~5.37 in \cite{bodirsky2018} is false, and as a consequence the proof of Theorem~5.1 in that paper is incorrect.
    Nevertheless, the statement of Theorem~5.1 can still be true: this holds provided that every first-order reduct of $(\mathbb{Z},<)$ that contains $\Succ$, but is neither a first-order reduct of $(\mathbb{Z},\Succ)$ nor pp-defines $<$, instead pp-defines $\Dist_n$ for every $n\in\mathbb{N}$.
\end{remark}
\begin{definition}
    Let $\hat{\cdot}:\mathbb{Q}\to\mathbb{Q}\times\mathbb{Z},p\mapsto(p,0)$.
    Given a first-order reduct $\Gamma$ of $(\mathbb{Z},<)$ containing $\Succ$, for $e\in\End(\mathbb{Q}.\Gamma)$ we define $\tilde{e}\coloneq\pi_1\circ e\circ\hat{\cdot}$.
    We write $\widetilde{\End}(\mathbb{Q}.\Gamma)\coloneq\{\tilde{e}\mid e\in\End(\mathbb{Q}.\Gamma)\}$.
\end{definition}
\begin{lemma}\label{lem:global_endomorphism_to_local}
    Let $\Gamma$ be a first-order reduct of $(\mathbb{Z},<)$ containing $\Succ$.
    For all $e\in\End(\mathbb{Q}.\Gamma)$ and $\bar{b}\in\mathbb{Z}^\mathbb{Q}$ the map $\mathbb{Q}\times\mathbb{Z}\to\mathbb{Q}\times\mathbb{Z},(p,a)\mapsto(\tilde{e}(p),a+b_p)$ is an endomorphism of $\mathbb{Q}.\Gamma$.
\end{lemma}
\begin{proof}
    Since $e$ preserves $\Succ$, for every $p\in\mathbb{Q}$ there is some $c_p\in\mathbb{Z}$ such that $e(p,c_p)=(\tilde{e}(p),b_p)$.
    Let $f\in\Aut(\mathbb{Q}\times\mathbb{Z},<)$ such that $f(p,0)=(p,c_p)$ for all $p\in\mathbb{Q}$.
    Then $g\coloneq e\circ f$ is an endomorphism of $\mathbb{Q}.\Gamma$ with $g(p,0)=e(p,c_p)=(\tilde{e}(p),b_p)$ for all $p\in\mathbb{Q}$.
    Since $g$ preserves $\Succ$ we get $g(p,a)=(\tilde{e}(p),a+b_p)$ as desired.
\end{proof}
\begin{fact}\label{fact:closed_transformation_monoid}
    Let $A$ be a set.
    Then $E\subseteq A^A$ is a closed submonoid of $A^A$ with the topology of pointwise convergence iff there exists a relational structure $\mathcal{A}$ with domain $A$ and $E=\End(\mathcal{A})$.
\end{fact}
\begin{lemma}\label{lem:fo_reduct_of_q_capturing_global_endomorphisms}
    Let $\Gamma$ be a first-order reduct of $(\mathbb{Z},<)$ containing $\Succ$.
    Then there exists a first-order reduct $\widetilde{\Gamma}$ of $(\mathbb{Q},<)$ such that $\widetilde{\End}(\mathbb{Q}.\Gamma)=\End(\widetilde{\Gamma})$.
\end{lemma}
\begin{proof}
    First, we show that $\widetilde{\End}(\mathbb{Q}.\Gamma)$ is closed in $\mathbb{Q}^\mathbb{Q}$ equipped with the topology of pointwise convergence.
    Let $\Sigma$ be a directed set and $(e_\sigma)_{\sigma\in\Sigma}\in\widetilde{\End}(\mathbb{Q}.\Gamma)^\Sigma$ with a limit point $f\in\mathbb{Q}^\mathbb{Q}$.
    Then for every $p\in\mathbb{Q}$, $f(p)$ is a limit point of $(e_\sigma(p))_{\sigma\in\Sigma}$ \cite[see][Proposition~2.6.5]{engelking1989}.
    Hence, $(f(p),a)$ is a limit point of $((e_\sigma(p),a))_{\sigma\in\Sigma}$ for every $(p,a)\in\mathbb{Q}\times\mathbb{Z}$.
    It follows that for $g_\sigma:\mathbb{Q}\times\mathbb{Z}\to\mathbb{Q}\times\mathbb{Z},(p,a)\mapsto(e_\sigma(p),a)$ the function $h:\mathbb{Q}\times\mathbb{Z}\to\mathbb{Q}\times\mathbb{Z},(p,a)\mapsto(f(p),a)$ is a limit point of $(g_\sigma)_{\sigma\in\Sigma}$ in the topology of pointwise convergence on $(\mathbb{Q}\times\mathbb{Z})^{\mathbb{Q}\times\mathbb{Z}}$ \cite[see][Proposition~2.6.5]{engelking1989}.
    By Lemma~\ref{lem:global_endomorphism_to_local}, we have $g_\sigma\in\End(\mathbb{Q}.\Gamma)$ for every $\sigma\in\Sigma$.
    Since $\End(\mathbb{Q}.\Gamma)$ is closed (Fact~\ref{fact:closed_transformation_monoid}) we get that $h\in\End(\mathbb{Q}.\Gamma)$ \cite[see][Corollary~1.6.4]{engelking1989}.
    Thus, $f=\tilde{h}\in\widetilde{\End}(\mathbb{Q}.\Gamma)$ and hence $\widetilde{\End}(\mathbb{Q}.\Gamma)$ is closed \cite[see][Corollary~1.6.4]{engelking1989}.

    By Fact~\ref{fact:closed_transformation_monoid}, $\widetilde{\End}(\mathbb{Q}.\Gamma)$ being closed implies that it is the endomorphism monoid of a relational structure $\widetilde{\Gamma}$ with domain $\mathbb{Q}$.
    It is not hard to see that $\Aut(\mathbb{Q},<)\subseteq\widetilde{\End}(\mathbb{Q}.\Gamma)$ and thus every relation of $\widetilde{\Gamma}$ is preserved by $\Aut(\mathbb{Q},<)$.
    By \cite[Corollary~3]{bodirsky2010} it follows that $\widetilde{\Gamma}$ is a first-order reduct of $(\mathbb{Q},<)$.
\end{proof}
\begin{lemma}\label{lem:all_permutations_of_q_implies_fo_reduct_of_succ}
    Let $\Gamma$ be a first-order reduct of $(\mathbb{Z},<)$ containing $\Succ$.
    If $\widetilde{\End}(\mathbb{Q}.\Gamma)$ contains all permutations of $\mathbb{Q}$, then $\Gamma$ is a first-order reduct of $(\mathbb{Z},\Succ)$.
\end{lemma}
\begin{proof}
    It is easy to see that if $\widetilde{\End}(\mathbb{Q}.\Gamma)$ contains all permutations of $\mathbb{Q}$ then $\Aut(\mathbb{Q}\times\mathbb{Z},\Succ)\subseteq\Aut(\mathbb{Q}.\Gamma)$.
    The result now immediately follows from \cite[Theorem~3.4]{bodirsky2018}.
\end{proof}
\begin{lemma}\label{lem:constant_operation_implies_fo_reduct_of_succ}
    Let $\Gamma$ be a first-order reduct of $(\mathbb{Z},<)$ containing $\Succ$.
    If $\widetilde{\End}(\mathbb{Q}.\Gamma)$ contains a constant operation, then $\Gamma$ is a first-order reduct of $(\mathbb{Z},\Succ)$.
\end{lemma}
\begin{proof}
    Let $\prec$ be any linear order on $\mathbb{Q}$.
    We claim that there exists $e\in\End(\mathbb{Q}.\Gamma)$ such that $\tilde{e}(p)<\tilde{e}(q)$ for all $p,q\in\mathbb{Q}$ with $p\prec q$.
    Suppose otherwise.
    Then by Proposition~\ref{prop:compactness}, setting
    \[\Phi\coloneq\{\Succ^n(x_{\hat{p}},x_{\hat{q}})\mid p,q\in\mathbb{Q},p\neq q,n\in\mathbb{Z}\}\cup\{\lnot(x_{\hat{p}}<x_{\hat{q}})\mid p,q\in\mathbb{Q},p\prec q\},\]
    there are a finite $\{p_1,\dots,p_k\}\subseteq\mathbb{Q}$ with $p_1\prec\dots\prec p_k$ and $n\in\mathbb{N}$ such that for all $e\in\End(\mathbb{Q}.\Gamma)$ there are distinct $i,j\in[k]$ such that $|e(\hat{p_i})-e(\hat{p_j})|<n$ or both $i<j$ and $e(\hat{p_i})\geq e(\hat{p_j})$.
    However, by assumption there exists $f\in\End(\mathbb{Q}.\Gamma)$ such that $\tilde{f}(\mathbb{Q})=\{r\}$ for some $r\in\mathbb{Q}$ and with Lemma~\ref{lem:global_endomorphism_to_local} it follows that there exists $e\in\End(\mathbb{Q}.\Gamma)$ with $e(\hat{p_i})=(r,in)$ for all $i\in[k]$.
    This is a contradiction since $e(\hat{p_j})-e(\hat{p_i})=(j-i)n\geq n$ for all $i,j\in[k]$ with $i<j$.

    By Lemma~\ref{lem:fo_reduct_of_q_capturing_global_endomorphisms} there is a first-order reduct $\widetilde{\Gamma}$ of $(\mathbb{Q},<)$ with $\End(\widetilde{\Gamma})=\widetilde{\End}(\mathbb{Q}.\Gamma)$.
    Noting that the $\tilde{e}$ from our claim are injective, it follows from \cite[Proposition~19]{bodirsky2010} that $\End(\widetilde{\Gamma},\neq)$ is the set of all injective unary operations on $\mathbb{Q}$ and thus $\widetilde{\End}(\mathbb{Q}.\Gamma)$ contains all permutations of $\mathbb{Q}$.
    By Lemma~\ref{lem:all_permutations_of_q_implies_fo_reduct_of_succ} we get that $\Gamma$ is a first-order reduct of $(\mathbb{Z},\Succ)$.
\end{proof}
\begin{lemma}\label{lem:minus_implies_fo_reduct_of_succ}
    Let $\Gamma$ be a first-order reduct of $(\mathbb{Z},<)$ containing $\Succ$.
    If $\widetilde{\End}(\mathbb{Q}.\Gamma)$ contains $\mathbb{Q}\to\mathbb{Q},p\mapsto-p$, then $\Gamma$ is a first-order reduct of $(\mathbb{Z},\Succ)$.
\end{lemma}
\begin{proof}
    If $\widetilde{\End}(\mathbb{Q}.\Gamma)$ contains a constant operation or all permutations of $\mathbb{Q}$, we are done by Lemmas~\ref{lem:all_permutations_of_q_implies_fo_reduct_of_succ} and~\ref{lem:constant_operation_implies_fo_reduct_of_succ}.
    Suppose for contradiction that this is not the case.
    By Lemma~\ref{lem:fo_reduct_of_q_capturing_global_endomorphisms} there exists a first-order reduct $\widetilde{\Gamma}$ of $(\mathbb{Q},<)$ with $\widetilde{\End}(\mathbb{Q}.\Gamma)=\End(\widetilde{\Gamma})$ and by \cite[Proposition~19]{bodirsky2010} $\End(\widetilde{\Gamma})$ preserves the relation defined by the formula
    \begin{align*}
        &&&\textcolor{red}{x_1}<\textcolor{blue}{x_2}<\textcolor{red}{x_3}<\textcolor{blue}{x_4}&&\lor&&\textcolor{red}{x_1}<\textcolor{blue}{x_4}<\textcolor{red}{x_3}<\textcolor{blue}{x_2}&&\lor&&\textcolor{red}{x_3}<\textcolor{blue}{x_2}<\textcolor{red}{x_1}<\textcolor{blue}{x_4}&&\lor&&\textcolor{red}{x_3}<\textcolor{blue}{x_4}<\textcolor{red}{x_1}<\textcolor{blue}{x_2}\\
        &\lor&&\textcolor{blue}{x_2}<\textcolor{red}{x_1}<\textcolor{blue}{x_4}<\textcolor{red}{x_3}&&\lor&&\textcolor{blue}{x_2}<\textcolor{red}{x_3}<\textcolor{blue}{x_4}<\textcolor{red}{x_1}&&\lor&&\textcolor{blue}{x_4}<\textcolor{red}{x_1}<\textcolor{blue}{x_2}<\textcolor{red}{x_3}&&\lor&&\textcolor{blue}{x_4}<\textcolor{red}{x_3}<\textcolor{blue}{x_2}<\textcolor{red}{x_1},
    \end{align*}
    which we call $\phi_{\Sep}(\bar{x})$.
    Let $p_1,\dots,p_4\in\mathbb{Q}$ with $p_1<\dots<p_4$.
    Then by Proposition~\ref{prop:compactness}, setting $\Phi\coloneq\{\phi_{\Sep}(x_{\hat{p_1}},\dots,x_{\hat{p_4}})\}$, there are a finite $A\subseteq\mathbb{Q}\times\mathbb{Z}$ containing $\hat{p_1},\dots,\hat{p_4}$ and a finite subset $\rho$ of the signature of $\Gamma$ such that for all $h\in\Hom(\mathcal{A},\mathbb{Q}.\Gamma\restriction\rho)$ where $\mathcal{A}\coloneq(\mathbb{Q}.\Gamma|_A)\restriction\rho$, we have $(\mathbb{Q}\times\mathbb{Z},<)\models\phi_{\Sep}(h(\hat{p_1}),\dots,h(\hat{p_4}))$.
    Since $A$ is finite, we can choose $n\in\mathbb{N}$ sufficiently large such that by Lemma~\ref{lem:substitution}
    \[f:\mathcal{A}\to\mathbb{Q}.\Gamma\restriction\rho,(q,a)\mapsto
    \begin{cases}
        (q,a) & \text{if }q\leq p_3\\
        (p_3,a+\lceil(q-p_3)n\rceil) & \text{if }q>p_3
    \end{cases}\]
    is a homomorphism.
    Let $g\in\Aut(\mathbb{Q}\times\mathbb{Z},<)$ with $g(\hat{p_2})=(0,0)$.
    By assumption and Lemma~\ref{lem:global_endomorphism_to_local}, $e:\mathbb{Q}.\Gamma\to\mathbb{Q}.\Gamma,(q,a)\mapsto(-q,a)$ is an endomorphism and thus $h\coloneq e\circ g\circ f$ is a homomorphism $\mathcal{A}\to\mathbb{Q}.\Gamma\restriction\rho$.
    \begin{center}
        \begin{tikzpicture}
            \foreach \y in {0,-1,-2} {
                \draw[<->] (-5,\y) -- (-1,\y);
                \foreach \x in {-4.5,-4,-3.5,-3,-2.5,-2,-1.5} {
                    \draw (\x,\y+0.1) -- (\x,\y-0.1);
                }
                \draw[<->] (1,\y) -- (5,\y);
                \foreach \x in {1.5,2,2.5,3,3.5,4,4.5} {
                    \draw (\x,\y+0.1) -- (\x,\y-0.1);
                }
            }
            \node (f1) at (-4,-1.6) {$f(\hat{p_1})$};
            \node (f2) at (-4,-0.6) {$f(\hat{p_2})$};
            \node (f3) at (-4,0.4) {$f(\hat{p_3})$};
            \node (f4) at (-2,0.4) {$f(\hat{p_4})$};

            \node (h3) at (2,-1.6) {$h(\hat{p_3})$};
            \node (h4) at (4,-1.6) {$h(\hat{p_4})$};
            \node (h2) at (2,-0.6) {$h(\hat{p_2})$};
            \node (h1) at (2,0.4) {$h(\hat{p_1})$};
            \draw[thick,->] (-0.5,-1) to node[above] {$e\circ g$} (0.5,-1);
        \end{tikzpicture}
    \end{center}
    We now have $h(\hat{p_3})<h(\hat{p_4})<h(\hat{p_2})<h(\hat{p_1})$, which contradicts $(\mathbb{Q}\times\mathbb{Z},<)\models\phi_{\Sep}(h(\hat{p_1}),\dots,h(\hat{p_4}))$.
\end{proof}
\begin{lemma}\label{lem:preserving_order_implies_pp_defining_order}
    Let $\Gamma$ be a first-order reduct of $(\mathbb{Z},<)$ containing $\Succ$.
    If $\widetilde{\End}(\mathbb{Q}.\Gamma)$ preserves $<$, then $\mathbb{Q}.\Gamma$ pp-defines $<$.
\end{lemma}
\begin{proof}
    Let $p_1,p_2\in\mathbb{Q}$ with $p_1<p_2$.
    By Proposition~\ref{prop:compactness}, setting $\Phi\coloneq\{x_{\hat{p_1}}<x_{\hat{p_2}}\}$, there are a finite $A\subseteq\mathbb{Q}\times\mathbb{Z}$ containing $\hat{p_1},\hat{p_2}$ and a finite subset $\rho$ of the signature of $\Gamma$ such that for all $h\in\Hom(\mathcal{A},\mathbb{Q}.\Gamma\restriction\rho)$, where $\mathcal{A}\coloneq(\mathbb{Q}.\Gamma|_A)\restriction\rho$, we have $h(\hat{p_1})<h(\hat{p_2})$.
    From Lemma~\ref{lem:substitution} it follows that there exists $n\in\mathbb{N}$ such that for all $a,b\in\mathbb{Q}\times\mathbb{Z}$ with $b-a\geq n$ there exists a homomorphism $h:\mathcal{A}\to\mathbb{Q}.\Gamma\restriction\rho$ with $h(\hat{p_1})=a$ and $h(\hat{p_2})=b$.
    Let $\phi(x_{\hat{p_1}},x_{\hat{p_2}})$ be the pp-formula obtained by existentially quantifying all variables but $x_{\hat{p_1}}$ and $x_{\hat{p_2}}$ in $\bigwedge\{Rx_{a_1}\dots x_{a_k}\mid R\in\rho,(a_1,\dots,a_k)\in R^\mathcal{A}\}$.
    Observe that $\mathbb{Q}.\Gamma\models\phi(a,b)$ iff there exists a homomorphism $h:\mathcal{A}\to\mathbb{Q}.\Gamma\restriction\rho$ with $h(\hat{p_1})=a$ and $h(\hat{p_2})=b$.
    Hence, $\phi$ pp-defines a relation $\Diff_S$ in $\mathbb{Q}.\Gamma$ where $S\subseteq\mathbb{N}_{>0}$ with $m\in S$ for all $m\geq n$.

    We fix $T$ minimizing $n'\in\mathbb{N}_{>0}$ among all $T\subseteq\mathbb{N}_{>0}$ such that $\mathbb{Q}.\Gamma$ pp-defines $\Diff_T$ and $m\in T$ for all $m\geq n'$.
    If $n'=1$ then $\Diff_T={<}$ and we are done, so suppose $n'>1$.
    Let $U\coloneq(T-1)\cap(T-n'+1)$ and note that by Observation~\ref{obs:shift_binary_relations} $\mathbb{Q}.\Gamma$ pp-defines $\Diff_U$.
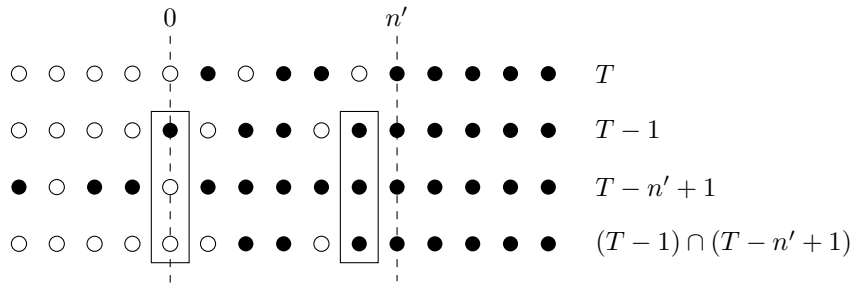
    \begin{center}
        \begin{tikzpicture}
            \node[anchor=west] at (5.5,0) {$T$};
            \node[anchor=west] at (5.5,-0.75) {$T-1$};
            \node[anchor=west] at (5.5,-1.5) {$T-n'+1$};
            \node[anchor=west] at (5.5,-2.25) {$(T-1)\cap(T-n'+1)$};
            \node[anchor=south] (0) at (0,0.5) {$0$};
            \node[anchor=south] (n) at (3,0.5) {$n'$};
            \draw[dashed] (0) -- +(0,-3.5);
            \draw[dashed] (n) -- +(0,-3.5);
            \foreach \x in {0.5,1.5,2,3,3.5,4,4.5,5} {
                \fill (\x,0) circle (1mm);
            }
            \foreach \x in {-2,-1.5,-1,-0.5,0,1,2.5} {
                \draw[black,fill=white] (\x,0) circle (1mm);
            }
            \foreach \x in {0.5,1.5,2,3,3.5,4,4.5,5,5.5} {
                \fill (\x-0.5,-0.75) circle (1mm);
            }
            \foreach \x in {-1.5,-1,-0.5,0,1,2.5} {
                \draw[black,fill=white] (\x-0.5,-0.75) circle (1mm);
            }
            \foreach \x in {0.5,1.5,2,3,3.5,4,4.5,5,5.5,6,6.5,7,7.5} {
                \fill (\x-2.5,-1.5) circle (1mm);
            }
            \foreach \x in {1,2.5} {
                \draw[black,fill=white] (\x-2.5,-1.5) circle (1mm);
            }
            \foreach \x in {1,1.5,2.5,3,3.5,4,4.5,5} {
                \fill (\x,-2.25) circle (1mm);
            }
            \foreach \x in {-2,-1.5,-1,-0.5,0,0.5,2} {
                \draw[black,fill=white] (\x,-2.25) circle (1mm);
            }
            \draw (-0.25,-0.5) rectangle (0.25,-2.5);
            \draw (2.25,-0.5) rectangle (2.75,-2.5);
        \end{tikzpicture}
    \end{center}
    Now, $U\subseteq T-1\subseteq\mathbb{N}$ and because $n'-1\notin T$ we have $0\notin T-n'+1\supseteq U$ and thus $U\subseteq\mathbb{N}_{>0}$.
    Moreover, for every $m\geq n'-1$ we have $m\in U$, which contradicts our choice of $T$.
\end{proof}
\begin{proposition}\label{prop:defining_order}
    Let $\Gamma$ be a first-order reduct of $(\mathbb{Z},<)$ that contains $\Succ$.
    Then one of the following holds:
    \begin{enumerate}
        \item $\Gamma$ is a first-order reduct of $(\mathbb{Z},\Succ)$;
        \item $\Gamma$ pp-defines $<$;
        \item $\Csp(\Gamma)$ is $\Np$-hard.
    \end{enumerate}
\end{proposition}
\begin{proof}
    If $\widetilde{\End}(\mathbb{Q}.\Gamma)$ contains all permutations of $\mathbb{Q}$, a constant operation or the map $\mathbb{Q}\to\mathbb{Q},p\mapsto-p$, then $\Gamma$ is a first-order reduct of $(\mathbb{Z},\Succ)$ by Lemmas~\ref{lem:all_permutations_of_q_implies_fo_reduct_of_succ},~\ref{lem:constant_operation_implies_fo_reduct_of_succ} and~\ref{lem:minus_implies_fo_reduct_of_succ} respectively.
    In the case that $\widetilde{\End}(\mathbb{Q}.\Gamma)$ preserves $<$ we get from Lemma~\ref{lem:preserving_order_implies_pp_defining_order} that $\Gamma$ pp-defines $<$.
    Since by Lemma~\ref{lem:fo_reduct_of_q_capturing_global_endomorphisms} there exists a first-order reduct $\widetilde{\Gamma}$ of $(\mathbb{Q},<)$ with $\widetilde{\End}(\mathbb{Q}.\Gamma)=\End(\widetilde{\Gamma})$, by \cite[Proposition~19]{bodirsky2010} the only remaining case is that
    $\Cycl=\bar{p}^{\widetilde{\End}(\mathbb{Q}.\Gamma)}$ for any $\bar{p}\in\Cycl\coloneq\{(p_1,p_2,p_3)\in\mathbb{Q}\mid p_1<p_2<p_3\text{ or }p_2<p_3<p_1\text{ or }p_3<p_1<p_2\}$.
    Fix $\bar{p}\in\Cycl$ and let $\phi_{\Cycl}(x_1,x_2,x_3)\coloneq x_1<x_2<x_3\lor x_2<x_3<x_1\lor x_3<x_1<x_2$.
    Then by Proposition~\ref{prop:compactness}, setting $\Phi\coloneq\{\phi_{\Cycl}(x_{\hat{p_1}},x_{\hat{p_2}},x_{\hat{p_3}})\}$, there are a finite $A\subseteq\mathbb{Q}\times\mathbb{Z}$ containing $\hat{p_1},\hat{p_2},\hat{p_3}$ and a finite subset $\rho$ of the signature of $\Gamma$ such that for all $h\in\Hom(\mathcal{A},\mathbb{Q}.\Gamma\restriction\rho)$ where $\mathcal{A}\coloneq(\mathbb{Q}.\Gamma|_A)\restriction\rho$ we have $(\mathbb{Q}\times\mathbb{Z},<)\models\phi_{\Cycl}(h(\hat{p_1}),h(\hat{p_2}),h(\hat{p_3}))$.
    By Lemma~\ref{lem:substitution} there is some $n\in\mathbb{N}$ such that for all $\bar{a}\in(\mathbb{Q}\times\mathbb{Z})^3$ with $(\mathbb{Q}\times\mathbb{Z},<)\models\phi_{\Cycl}(\bar{a})$ and $|a_i-a_j|\geq n$ for all distinct $i,j\in[3]$ there exists a homomorphism $h:\mathcal{A}\to\mathbb{Q}.\Gamma\restriction\rho$ with $h(\hat{p_i})=a_i$ for all $i\in[3]$.
    Let $\psi(x_{\hat{p_1}},x_{\hat{p_2}},x_{\hat{p_3}})$ be the pp-formula obtained by existentially quantifying all variables but $x_{\hat{p_1}},x_{\hat{p_2}},x_{\hat{p_3}}$ in $\bigwedge\{Rx_{a_1}\dots x_{a_k}\mid R\in\rho,(a_1,\dots,a_k)\in R^\mathcal{A}\}$.
    We now have $\mathbb{Q}.\Gamma\models\psi(\bar{a})$ iff there is a homomorphism $h:\mathcal{A}\to\mathbb{Q}.\Gamma\restriction\rho$ with $h(\hat{p_i})=a_i$ for all $i\in[3]$.
    Thus,
    \[\Th(\Gamma,<)\vdash\psi(\bar{x})\rightarrow\phi_{\Cycl}(\bar{x})\]
    and
    \[\Th(\Gamma,<)\cup\{\lnot\Succ^m(x_i,x_j)\mid m\in\mathbb{N}_{<n},i\neq j\}\vdash\psi(\bar{x})\leftrightarrow\phi_{\Cycl}(\bar{x}).\]
    It follows that in $\Gamma$ the pp-formula $\exists z(\psi(x,y,z)\land\Succ^{2n+3}(x,z))$ defines $\Diff_S$ for some finite $S\subseteq\mathbb{Z}$ containing $n+1$ and $n+2$.
    Moreover, $\exists y(\psi(x,y,z)\land\Succ^{n+1}(x,y))$ pp-defines $\Diff_T$ for a cofinite $T\subseteq\mathbb{Z}$ with $0\notin T$.
    Now, by Proposition~\ref{prop:hardness_unless_all_d_progressions}, $\Csp(\mathbb{Z},\Succ,\Diff_S,\Diff_T)$ and thus also $\Csp(\Gamma)$ are $\Np$-hard.
\end{proof}

\section{Classification}\label{sec:classification}
We now combine the results from Sections~\ref{sec:defining_successor} to~\ref{sec:defining_order} to prove the full complexity classification.
\begin{lemma}\label{lem:csp_in_np}
    Let $\Gamma$ be a finite-signature first-order reduct of $(\mathbb{Z},<)$.
    Then $\Csp(\Gamma)\in\Np$.
\end{lemma}
\begin{proof}
    By Lemma~\ref{lem:substitution} there exists $n\in\mathbb{N}$ such that a finite structure $\mathcal{A}$ maps homomorphically to $\Gamma$ iff it maps homomorphically to $\Gamma|_{[m]}$ for $m\coloneq(n+1)\cdot|A|$.
    Hence, we can use a homomorphism $h:\mathcal{A}\to\Gamma|_{[m]}$ as a certificate for the existance of a homomorphism $\mathcal{A}\to\Gamma$.
    Using the quantifier-free definitions of the relations of $\Gamma$, which exist by Fact~\ref{fact:quantifier_elimination_order}, we can efficiently verify that $h$ is a homomorphism $\mathcal{A}\to\Gamma|_{[m]}$.
\end{proof}
\begin{theorem}[{\cite[see][Theorem~4.2]{bodirsky2018}}]\label{thm:classification}
    Let $\Gamma$ be a finite-signature first-order reduct of $(\mathbb{Z},<)$.
    Then there exists $\Delta$ with $\Csp(\Gamma)=\Csp(\Delta)$ such that one of the following holds.
    \begin{enumerate}
        \item\label{item:classification_integers_finite} $\Delta$ is finite; then $\Csp(\Delta)\in\Ptime\cup\Npc$ by \cite{bulatov2017} and \cite{zhuk2017,zhuk2020}.
        \item\label{item:classification_integers_rationals} $\Delta$ is a first-order reduct of $(\mathbb{Q},<)$; then $\Csp(\Delta)\in\Ptime\cup\Npc$ by \cite{bodirsky2010}.
        \item\label{item:classification_integers_min_max} $\Delta$ is a first-order reduct of $(\mathbb{Z},<)$ preserved by $\Min_d$ or $\Max_d$ for some $d$; then $\Csp(\Delta)\in\Ptime$.
        \item\label{item:classification_integers_horn} $\Delta$ is a first-order Horn reduct of $(\mathbb{Z},\Succ)$; then $\Csp(\Delta)\in\Ptime$.
        \item\label{item:classification_integers_hard} $\Csp(\Delta)\in\Npc$.
    \end{enumerate}
\end{theorem}
\begin{proof}
    We have $\Csp(\Gamma)\in\Np$ by Lemma~\ref{lem:csp_in_np}.
    By Theorem~\ref{thm:defining_successor}, there exists $\Delta$ with $\Csp(\Gamma)=\Csp(\Delta)$, and one of the following four cases applies.
    If $\Delta$ is finite, we are in Case~\ref{item:classification_integers_finite}.
    If $\Delta$ is a first-order reduct of $(\mathbb{Q},<)$, we are in Case~\ref{item:classification_integers_rationals}.
    If $\Delta$ is a first-order reduct of $(\mathbb{Z},<)$ that pp-defines $\Dist_n$ for all $n\in\mathbb{N}$, then by Lemma~\ref{lem:hardness_from_distances} we are in Case~\ref{item:classification_integers_hard}.
    Otherwise, $\Delta$ is a first-order reduct of $(\mathbb{Z},<)$ that pp-defines $\Succ$, and we assume from now on that $\Delta$ contains the relation $\Succ$.
    Now, by Proposition~\ref{prop:defining_order} we are in Case~\ref{item:classification_integers_hard} or either $\Delta$ is a first-order reduct of $(\mathbb{Z},\Succ)$ or pp-defines $<$.
    In any case, if there is $d\in\mathbb{N}_{>0}$ such that $\Delta$ admits $\Min_d$ or $\Max_d$ as a polymorphism, then by Lemma~\ref{lem:sampling_for_fo_reducts_of_integers} and Proposition~\ref{prop:tractability_of_mod_max_min} we have $\Csp(\Delta)\in\Ptime$ and are in Case~\ref{item:classification_integers_min_max}.
    So suppose that $\Delta$ does not admit $\Min_d$ or $\Max_d$ for any $d\in\mathbb{N}_{>0}$.
    If $\Delta$ is a first-order reduct of $(\mathbb{Z},\Succ)$, then either it is Horn and we are in Case~\ref{item:classification_integers_horn} by Lemma~\ref{lem:sampling_for_fo_reducts_of_integers} and Proposition~\ref{prop:tractable_horn}, or it is not Horn and we are in Case~\ref{item:classification_integers_hard} by Proposition~\ref{prop:non_horn_relations}.
    On the other hand, if $\Delta$ pp-defines $<$, then $\exists uv(\Succ(u,x)\land\Succ^2(x,v)\land u<y<v)$ pp-defines $\Diff_{\{0,1\}}$ and hence we are in Case~\ref{item:classification_integers_hard} by Lemma~\ref{lem:one_progression_and_no_min_max_is_hard}.
\end{proof}
\begin{corollary}
    Let $\Gamma$ be a finite-signature first-order reduct of $(\mathbb{Z},<)$.
    Then $\Csp(\Gamma)\in\Ptime\cup\Npc$.
\end{corollary}
To conclude our study of first-order reducts of $(\mathbb{Z},<)$, we show that a single uniform algorithm suffices to solve all tractable CSPs of first-order reducts of $(\mathbb{Z},<)$ that are not the CSP of a finite structure.
\begin{definition}
    A \emph{height~1 identity} is an expression of the form $f(x_{i_1},\dots,x_{i_m})\approx g(x_{j_1},\dots,x_{j_n})$, where $f$ and $g$ are function symbols of arity $m$ and $n$, respectively, $i_1,\dots,i_m,j_1,\dots,j_n\in\mathbb{N}$, and $x_0,x_1,\dots$ are variables.
    An interpretation $\mathcal{I}$, sending each $k$-ary function symbol to some function in $B^{A^k}$, \emph{satisfies} such an identity if $\mathcal{I}(f)(a_{i_1},\dots,a_{i_m})=\mathcal{I}(g)(a_{j_1},\dots,a_{j_n})$ for all $\bar{a}\in A^\mathbb{N}$.
    We call a finite set $\Sigma$ of height~1 identities a \emph{height~1 condition}.
    A set of functions $F\subseteq\bigcup_{k\in\mathbb{N}}B^{A^k}$ \emph{satisfies} $\Sigma$ if there exists an interpretation $\mathcal{I}$, sending every function symbol occurring in $\Sigma$ to a function in $F$ of matching arity, such that $\mathcal{I}$ satisfies every identity in $\Sigma$.
\end{definition}
\begin{lemma}\label{lem:height_one_conditions}
    Let $\mathcal{M}$ and $\mathcal{N}$ be relational structures with finite signatures and $\Csp(\mathcal{M})=\Csp(\mathcal{N})$.
    If a height~1 condition is satisfied by $\Pol(\mathcal{M}|_A,\mathcal{M})$ for every finite $A\subseteq M$, then it is satisfied by $\Pol(\mathcal{N}|_B,\mathcal{N})$ for every finite $B\subseteq N$.

    The same holds for $(k_1,\dots,k_n)\textnormal{-}\Pol(\mathcal{M}|_A,\mathcal{M})$ and $(k_1,\dots,k_n)\textnormal{-}\Pol(\mathcal{N}|_B,\mathcal{N})$.
\end{lemma}
\begin{proof}
    Suppose $\Sigma$ is a height~1 condition satisfied by $\Pol(\mathcal{M}|_A,\mathcal{M})$ for all finite $A\subseteq M$ and let $B\subseteq N$ be finite.
    Since $\Csp(\mathcal{M})=\Csp(\mathcal{N})$, there exists a homomorphism $h:\mathcal{N}|_B\to\mathcal{M}$.
    Let $f_1,\dots,f_m\in\Pol(\mathcal{M}|_{h(B)},\mathcal{M})$ of arities $\ell_1,\dots,\ell_m\in\mathbb{N}$ satisfy $\Sigma$.
    Because $\Csp(\mathcal{M})=\Csp(\mathcal{N})$, for $A\coloneq\bigcup_{i\in[m]}f_i(h(B)^{\ell_i})$ there is a homomorphism $g:\mathcal{M}|_A\to\mathcal{N}$.
    It is now easy to verify that $f_1',\dots,f_m'$ with $f_i':B^{\ell_i}\to N,(a_1,\dots,a_{k_i})\mapsto(g\circ f_i)(h(a_1),\dots,h(a_{k_i}))$ for all $i\in[m]$ satisfy $\Sigma$ in $\Pol(\mathcal{N}|_B,\mathcal{N})$.

    The proof for $(k_1,\dots,k_n)\textnormal{-}\Pol(\mathcal{M}|_A,\mathcal{M})$ and $(k_1,\dots,k_n)\textnormal{-}\Pol(\mathcal{N}|_B,\mathcal{N})$ is analogous.
    We only have to check that $f_1',\dots,f_m'\in(k_1,\dots,k_n)\textnormal{-}\Pol(\mathcal{N}|_B,\mathcal{N})$.
    This follows immediately from $g$ and $h$ being homomorphisms and the fact that applying a unary function componentwise preserves the property of being a $(k_1,\dots,k_n)$-palette tuple.
\end{proof}
A useful feature of uniform algorithms is that they can be combined: for instance, running two of them independently on the same input and accepting only if both accept yields another uniform algorithm, one that solves every CSP solved by either of the two.
The algorithm \Csac~$\land$~\Saip{} is exactly this combination of \Csac{} and \Saip{} \cite[see][]{zhuk2026}.
By the next theorem, this also means that a uniform algorithm for finite-domain CSPs, should one be found, would immediately yield a uniform algorithm solving all tractable CSPs of first-order reducts of $(\mathbb{Z},<)$.
\begin{theorem}
    Let $\Gamma$ be a finite-signature first-order reduct of $(\mathbb{Z},<)$.
    Then one of the following holds.
    \begin{enumerate}
        \item $\Csp(\Gamma)=\Csp(\Delta)$ for a finite structure $\Delta$.
        \item $\Pcsp(\Gamma|_{[n]},\Gamma)$ is solved by \Csac~$\land$~\Saip{} for every $n\in\mathbb{N}$.
        \item $\Csp(\Gamma)$ is $\Np$-complete.
    \end{enumerate}
\end{theorem}
\begin{proof}
    Suppose $\Csp(\Gamma)$ is neither the CSP of a finite structure nor $\Np$-hard.
    Then, by Theorem~\ref{thm:classification}, $\Csp(\Gamma)=\Csp(\Delta)$ for some $\Delta$ that is a first-order reduct of $(\mathbb{Q},<)$, a first-order reduct of $(\mathbb{Z},<)$ admitting a modular min or max polymorphism, or a first-order Horn reduct of $(\mathbb{Z},\Succ)$.
    By \cite{zhuk2026}, to show that $\Pcsp(\Gamma|_{[n]},\Gamma)$ is solved by \Csac~$\land$~\Saip{}, it suffices to establish one of the following.
    \begin{itemize}
        \item $\Pol(\Gamma|_{[n]},\Gamma)$ contains a so-called $(2k_1+1,\dots,2k_\ell+1)$-palette alternating function for all $\ell,k_1,\dots,k_\ell\in\mathbb{N}_{>0}$.
            Equivalently, $\Pol(\Gamma|_{[n]},\Gamma)$ satisfies a certain height~1 condition $\Sigma_\text{\Saip}^{k_1,\dots,k_\ell}$ for all $\ell,k_1,\dots,k_\ell\in\mathbb{N}_{>0}$.
        \item For all $\ell,k_1,\dots,k_\ell\in\mathbb{N}_{>0}$ the set $(k_1,\dots,k_\ell)\textnormal{-}\Pol(\Gamma|_{[n]},\Gamma)$ contains a $(k_1,\dots,k_\ell)$-palette totally symmetric function.
            Equivalently, $(k_1,\dots,k_\ell)\textnormal{-}\Pol(\Gamma|_{[n]},\Gamma)$ satisfies a certain height~1 condition $\Sigma_\text{\Csac}^{k_1,\dots,k_\ell}$ for all $\ell,k_1,\dots,k_\ell\in\mathbb{N}_{>0}$.
    \end{itemize}
    Since every finite substructure of $\Delta$ embeds into $\Delta|_{[m]}$ for some $m\in\mathbb{N}$, Lemma~\ref{lem:height_one_conditions} lets us establish the above for $\Delta$ instead of $\Gamma$.
    For $\Delta$ a first-order reduct of $(\mathbb{Q},<)$, this was done by \textcite[Section~6]{zhuk2026}.
    For $\Delta$ a first-order reduct of $(\mathbb{Z},<)$ admitting a modular min or max polymorphism, the statement follows from Lemma~\ref{lem:mod_max_min_palette_totally_symmetric}.
    For $\Delta$ a first-order Horn reduct of $(\mathbb{Z},\Succ)$, it follows from Lemma~\ref{lem:horn_palette_totally_symmetric}.
\end{proof}

\section{The triangular lattice graph}\label{sec:grids}
In this section, we give a first, simple application of the techniques from Section~\ref{sec:model_theory} beyond $(\mathbb{Z},<)$: we classify the complexity of the CSP for first-order reducts of the triangular lattice graph, and establish a $\Ptime$ versus $\Np$-complete dichotomy for this class of structures.
Throughout this section, let $\mathcal{G}\coloneq(V,E)$ denote the triangular lattice graph, illustrated in Figure~\ref{fig:grid}.
\begin{figure}[ht]
    \centering
    \begin{tikzpicture}
        \clip (0:2.5) -- (60:2.5) -- (120:2.5) -- (180:2.5) -- (240:2.5) -- (300:2.5) -- cycle;
        \foreach \x in {-5,...,5}
            \foreach \y in {-5,...,5}
                \foreach \d in {(1,0),(.5,{sqrt(3)/2}),(-.5,{sqrt(3)/2})}
                    \draw[black!50] ({\x+\y/2},{sqrt(3)*\y/2}) -- ++\d;
    \end{tikzpicture}
    \caption{A finite part of the triangular lattice graph.}
    \label{fig:grid}
\end{figure}
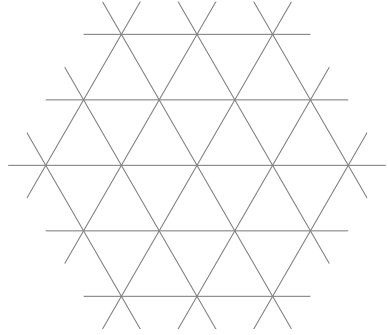
\begin{lemma}\label{lem:grids_csp_in_np}
    Let $\Gamma$ be a finite-signature first-order reduct of $\mathcal{G}$.
    Then $\Csp(\Gamma)\in\Np$.
\end{lemma}
\begin{proof}
    The structure $\mathcal{G}$ is interpretable in $(\mathbb{Z},\Succ)$.
    Indeed, it is not hard to see that $(\mathbb{Z}^2,E)$ with $E\subseteq(\mathbb{Z}^2)^2$ defined by
    \begin{align*}
        \phi(x_1,y_1,x_2,y_2)\coloneq&(x_1=x_2\land(\Succ(y_1,y_2)\lor\Succ(y_2,y_1)))\\
        \lor&(y_1=y_2\land(\Succ(x_1,x_2)\lor\Succ(x_2,x_1)))\\
        \lor&(\Succ(x_1,x_2)\land\Succ(y_1,y_2))\\
        \lor&(\Succ(x_2,x_1)\land\Succ(y_2,y_1))
    \end{align*}
    is isomorphic to $\mathcal{G}$.
    \begin{center}
        \begin{tikzpicture}[baseline,scale=0.5]
            \clip (-2.5,-2.5) -- (-2.5,0) -- (0,2.5) --(2.5,2.5) --(2.5,0) -- (0,-2.5) -- cycle;
            \foreach \x in {-5,...,5}
                \foreach \y in {-5,...,5}
                    \foreach \d in {(0,1),(1,0),(1,1)}
                        \draw[black!50] (\x,\y) -- ++\d;
        \end{tikzpicture}
        \hspace{0.5cm}
        \begin{tikzpicture}[baseline]
            \draw[<->] (-1,0) to node[above] {\small isomorphic} (1,0);
        \end{tikzpicture}
        \hspace{0.5cm}
        \begin{tikzpicture}[baseline,scale=0.6]
            \clip (0:2.5) -- (60:2.5) -- (120:2.5) -- (180:2.5) -- (240:2.5) -- (300:2.5) -- cycle;
            \foreach \x in {-5,...,5}
                \foreach \y in {-5,...,5}
                    \foreach \d in {(1,0),(.5,{sqrt(3)/2}),(-.5,{sqrt(3)/2})}
                        \draw[black!50] ({\x+\y/2},{sqrt(3)*\y/2}) -- ++\d;
        \end{tikzpicture}
    \end{center}
    It follows that $\Csp(\Gamma)$ reduces in polynomial time to $\Csp(\Delta)$ for a finite-signature first-order reduct $\Delta$ of $(\mathbb{Z},<)$.
    By Lemma~\ref{lem:csp_in_np}, $\Csp(\Delta)$ and hence $\Csp(\Gamma)$ are in $\Np$.
\end{proof}
As before, our approach relies on the saturated counterparts of first-order reducts of $\mathcal{G}$.
We construct such a counterpart by letting $\omega.\mathcal{G}\coloneq(\omega\times V,E)$, where $E^{\omega.\mathcal{G}}\coloneq\{((n,a),(n,b))\mid n\in\omega,\ (a,b)\in E^{\mathcal{G}}\}$, so that $\omega.\mathcal{G}$ consists of countably many disjoint copies of $\mathcal{G}$, and for a first-order reduct $\Gamma$ of $\mathcal{G}$, write $\omega.\Gamma$ for the first-order reduct of $\omega.\mathcal{G}$ over the same signature whose relations are defined by the same formulas.
As one would expect from this construction, $\omega.\mathcal{G}$ is nice, with $a\asymp b$ if and only if $a$ and $b$ lie in the same copy of $\mathcal{G}$.
For simplicity, we write $\sim_\Gamma$ instead of $\sim_{\omega.\Gamma}$ (see Definition~\ref{def:rank_relation}).
\begin{lemma}\label{lem:grids_with_rank_zero}
    Let $\Gamma$ be a first-order reduct of $\mathcal{G}$ with no finite-range endomorphism.
    If $\sim_\Gamma$ is trivial, then $\omega.\Gamma$ is homomorphically equivalent to a first-order reduct of $(\mathbb{N},=)$.
\end{lemma}
\begin{proof}
    Since $\Gamma$ does not have a finite-range endomorphism, by Lemma~\ref{lem:equivalence_of_finite_range_endomorphisms} neither does $\omega.\Gamma$.
    Hence, by Proposition~\ref{prop:endomorphism_with_image_in_relation} there exists $e\in\End(\omega.\Gamma)$ such that for all $a,b\in\omega\times V$ either $e(a)\sim_\Gamma e(b)$ or $e(a)\not\asymp e(b)$.
    Fix some $a\in V$.
    Since $\sim_\Gamma$ is trivial, we can compose $e$ with an appropriate automorphism of $\omega.\mathcal{G}$ and get an endomorphism of $\omega.\Gamma$ whose image is contained in $\omega\times\{a\}$.
    Let us assume that already $e(\omega\times V)\subseteq\omega\times\{a\}$.
    Then $\omega.\Gamma$ is homomorphically equivalent to $\omega.\Gamma|_{\omega\times\{a\}}$.
    Since every permutation of $\omega\times\{a\}$ preserves the relations of $\omega.\Gamma|_{\omega\times\{a\}}$, this structure is isomorphic to a first-order reduct of $(\mathbb{N},=)$ \cite[see][]{bodirsky2008}.
\end{proof}
\begin{proposition}\label{prop:grids_hardness}
    Let $\Gamma$ be a first-order reduct of $\mathcal{G}$ with no finite-range endomorphism.
    If $\sim_\Gamma$ is non-trivial, then $\Csp(\Gamma)$ is $\Np$-hard.
\end{proposition}
\begin{proof}
    By Lemmas~\ref{lem:equivalence_of_finite_range_endomorphisms} and~\ref{lem:avoiding_orbitals} there are $a,b\in\omega\times V$ such that $a\sim_\Gamma b$ and $e(a)\neq e(b)$ for all $e\in\End(\omega.\Gamma)$.
    It follows from Lemma~\ref{lem:pp_define_image_under_endomorphisms} that $R\coloneq(a,b)^{\End(\omega.\Gamma)}$ is a finite union of orbitals of $\Aut(\omega.\mathcal{G})$ and pp-definable in $\omega.\Gamma$.
    For each $d\in\omega\times V$, let $f_d$ be an automorphism of $\omega.\mathcal{G}$ that rotates the copy of $\mathcal{G}$ containing $d$ by $60^\circ$ clockwise about $d$.
    Moreover, let $c\coloneq f_a(b)$ and note that $c=f_b^5(a)$.
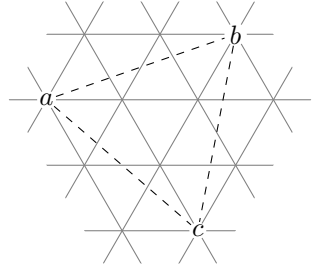
    \begin{center}
        \begin{tikzpicture}[
                node/.style={
                    circle,
                    fill=white,
                    inner sep=0.1pt,
                },
            ]
            \clip (0.75,{3*sqrt(3)/4}) -- (-1.75,{3*sqrt(3)/4}) -- (-2.5,0) -- (-1.25,-{5*sqrt(3)/4}) -- (0.25,-{5*sqrt(3)/4}) -- (1.5,0) -- cycle;
            \foreach \x in {-5,...,5}
                \foreach \y in {-5,...,5}
                    \foreach \d in {(1,0),(.5,{sqrt(3)/2}),(-.5,{sqrt(3)/2})}
                        \draw[black!50] ({\x+\y/2},{sqrt(3)*\y/2}) -- ++\d;
            \node[node] (a) at (-2,0) {$a$};
            \node[node] (b) at (0.5,{sqrt(3)/2}) {$b$};
            \node[node] (c) at (0,{-sqrt(3)}) {$c$};
            \draw[dashed] (a) to (b);
            \draw[dashed] (b) to (c);
            \draw[dashed] (a) to (c);
        \end{tikzpicture}
    \end{center}
    Then $(a,c)=(f_a(a),f_a(b))\in R$ and $(c,b)=(f_b^5(a),f_b^5(b))\in R$.
    Since $R$ is symmetric we get $(a,a),(a,b),(a,c)\in R^2$ and hence $A\coloneq\{d\mid (a,d)\in R^2\}$ contains $a$, $b$ and $c$.
    Moreover, $A$ is finite because $R\subseteq{\sim_\Gamma}\subseteq{\asymp}$ is a finite union of orbitals of $\Aut(\omega.\mathcal{G})$.
    Thus, $(A,R|_A)$ is a finite, undirected graph without loops that contains the triangle $a,b,c$ and by \cite{hell1990} its CSP is $\Np$-hard.
    The same argument as in the proof of Lemma~\ref{lem:use_relation_to_reduce_to_finite_csp} shows that $\Csp(A,R|_A)$ reduces in polynomial time to $\Csp(\Gamma)$, and hence that $\Csp(\Gamma)$ is $\Np$-hard as well.
\end{proof}
\begin{theorem}
    Let $\Gamma$ be a finite-signature first-order reduct of $\mathcal{G}$.
    Then one of the following holds.
    \begin{enumerate}
        \item\label{item:grids_finite} $\Gamma$ is homomorphically equivalent to a finite structure and $\Csp(\Gamma)\in\Ptime\cup\Npc$.
        \item\label{item:grids_trivial} $\omega.\Gamma$ is homomorphically equivalent to a first-order reduct of $(\mathbb{N},=)$ and $\Csp(\Gamma)\in\Ptime\cup\Npc$.
        \item\label{item:grids_hard} $\Csp(\Gamma)\in\Npc$.
    \end{enumerate}
\end{theorem}
\begin{proof}
    If $\Gamma$ has a finite-range endomorphism, we are in Case~\ref{item:grids_finite} by \cite{bulatov2017} and \cite{zhuk2017,zhuk2020}, so suppose otherwise.
    Then if $\sim_\Gamma$ is trivial, by Lemma~\ref{lem:grids_with_rank_zero} and \cite{bodirsky2008} we are in Case~\ref{item:grids_trivial}.
    On the other hand, if $\sim_\Gamma$ is non-trivial we are in Case~\ref{item:grids_hard} by Lemma~\ref{lem:grids_csp_in_np} and Proposition~\ref{prop:grids_hardness}.
\end{proof}

\section{Infinite regular trees}\label{sec:trees}
We conclude the paper with a second application of the techniques from Section~\ref{sec:model_theory}: we classify the complexity of the CSP for first-order reducts of infinite regular trees, and establish a $\Ptime$ versus $\Np$-complete dichotomy for this class of structures.
While the overall approach mirrors the one for the triangular lattice graph in Section~\ref{sec:grids}, the argument here is marginally more involved.
Throughout this section, fix some $n\in\mathbb{N}_{>1}$, and let $\mathcal{T}_n\coloneq(V,E)$ denote the infinite $n$-regular tree, i.e., the (up to isomorphism unique) infinite, undirected, connected, and acyclic graph in which every vertex has degree $n$.

As in Section~\ref{sec:grids}, the saturated model $\omega.\mathcal{T}_n\coloneq(\omega\times V,E)$ consists of countably many disjoint copies of $\mathcal{T}_n$; for a first-order reduct $\Gamma$ of $\mathcal{T}_n$, we again write $\omega.\Gamma$ for the first-order reduct of $\omega.\mathcal{T}_n$ over the same signature whose relations are defined by the same formulas.
This construction yields a nice model $\omega.\mathcal{T}_n$ of the theory of $\mathcal{T}_n$.

For vertices $u,v\in V$, write $|u-v|$ for their distance in $\mathcal{T}_n$, and extend this to $a,b\in\omega\times V$ by setting $|a-b|\coloneq\infty$ whenever $a$ and $b$ have different first coordinates.
It then holds that $a\asymp b\iff|a-b|<\infty$, and $(a,b)^{\Aut(\omega.\mathcal{T}_n)}=\{(c,d)\in(\omega\times V)^2\mid|a-b|=|c-d|\}$ for all $a,b\in\omega\times V$.
Finally, given $S\subseteq\mathbb{N}$, let $\Dist_S\coloneq\{(u,v)\in V^2\mid |u-v|\in S\}$, using the same notation for the corresponding relation on $\omega.\mathcal{T}_n$.
For $m\in\mathbb{N}$, we write $\Dist_m\coloneq\Dist_{\{m\}}$.
\begin{lemma}\label{lem:trees_quantifier_elimination}
    The structures $\mathcal{T}_n$ and $(V,(\Dist_m)_{m\in\mathbb{N}})$ are interdefinable and $\Th(V,(\Dist_m)_{m\in\mathbb{N}})$ has quantifier elimination.
\end{lemma}
\begin{proof}
    We have $E=\Dist_1$ and $\Dist_m$ is defined in $\mathcal{T}_n$ by $E^mxy\land\bigwedge_{\ell\in\{0,\dots,m-1\}}\lnot E^\ell xy$.
    For the quantifier elimination see \cite{hedman2011}.
\end{proof}
\begin{lemma}\label{lem:trees_csp_in_np}
    Let $\Gamma$ be a finite-signature first-order reduct of $\mathcal{T}_n$.
    Then $\Csp(\Gamma)\in\Np$.
\end{lemma}
\begin{proof}
    By Lemma~\ref{lem:trees_quantifier_elimination}, every relation $R$ of $\Gamma$ is definable in $(V,(\Dist_m)_{m\in\mathbb{N}})$ by a quantifier-free formula.
    Let $m\in\mathbb{N}$ be maximal such that $\Dist_m$ appears in one of those formulas.
    It is now easy to see that a finite structure $\mathcal{A}$ maps homomorphically to $\Gamma$ if and only if there exists a subtree $(V',E')$ of $\mathcal{T}_n$ of size $(m+1)\cdot|A|$ such that $\mathcal{A}$ maps homomorphically to $\Gamma|_{V'}$.
    Hence, we can use $(V',E')$ together with a homomorphism $h:\mathcal{A}\to\Gamma|_{V'}$ as a certificate for a homomorphism $\mathcal{A}\to\Gamma$.
    Clearly, its size is polynomial in $|A|$.
    Moreover, we can efficiently verify that $(V',E')$ is a subtree of $\mathcal{T}_n$, use the quantifier-free definitions of the relations of $\Gamma$ to compute $\Gamma|_{V'}$ from $(V',E')$ and verify that $h$ is a homomorphism.
\end{proof}
As before, we write $\sim_\Gamma$ instead of $\sim_{\omega.\Gamma}$ for simplicity.
\begin{lemma}\label{lem:trees_possible_ranks}
    Let $\Gamma$ be a first-order reduct of $\mathcal{T}_n$.
    Then there exists $r\in\mathbb{N}$ such that $a\sim_\Gamma b\iff|a-b|\in r\mathbb{N}$ for all $a,b\in\omega\times V$.
    If $n>2$, we have $r\in\{0,1,2\}$.
\end{lemma}
\begin{proof}
    The same idea as in the proof of Lemma~\ref{lem:equivalence_relations_are_periodic} shows that there exists $r\in\mathbb{N}$ with $a\sim_\Gamma b\iff|a-b|\in r\mathbb{N}$ for all $a,b\in\omega\times V$.
    Suppose $n>2$ and $r\neq 0$.
    Then there are distinct $a,b\in\omega\times V$ with $a\sim_\Gamma b$.
    Since $n>2$ there exists $c\in\omega\times V$ such that $|a-b|=|a-c|$ and $|b-c|=2$.
    
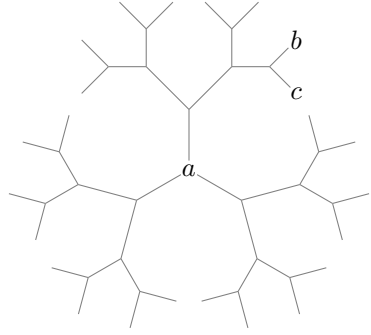
\begin{center}
        \begin{tikzpicture}[
                rotate=90,
                edge/.style={
                    draw=black!50,
                },
                node/.style={
                    circle,
                    fill=white,
                    inner sep=0.1pt,
                },
            ]
            \coordinate (0) at (0,0);
            \foreach \a in {0,1,2} {
                \begin{scope}[rotate={\a*120},shift={(0.8,0)}]
                    \coordinate (0\a) at (0,0);
                    \draw[edge] (0) -- (0\a);
                    \foreach \b in {0,1} {
                        \begin{scope}[rotate={\b*90-45},shift={(0.8,0)}]
                            \coordinate (0\a\b) at (0,0);
                            \draw[edge] (0\a) -- (0\a\b);
                            \foreach \c in {0,1} {
                                \begin{scope}[rotate={\c*90-45},shift={(0.5,0)}]
                                    \coordinate (0\a\b\c) at (0,0);
                                    \draw[edge] (0\a\b) -- (0\a\b\c);
                                    \foreach \d in {0,1} {
                                        \begin{scope}[rotate={\d*90-45},shift={(0.5,0)}]
                                            \coordinate (0\a\b\c\d) at (0,0);
                                            \draw[edge] (0\a\b\c) -- (0\a\b\c\d);
                                        \end{scope}
                                    }
                                \end{scope}
                            }
                        \end{scope}
                    }
                \end{scope}
            }
            \node[node] (a) at (0) {$a$};
            \node[node] (b) at (00001) {$b$};
            \node[node] (c) at (00000) {$c$};
        \end{tikzpicture}
    \end{center}
    Hence, there exists $f\in\Aut(\omega\times V,E)$ with $f(a)=a$ and $f(b)=c$.
    Because $\sim_\Gamma$ is preserved by $\Aut(\omega\times V,E)$ we get $a\sim_\Gamma c$ and thus $b\sim_\Gamma c$.
    It follows that $r$ divides $2$ and hence $r\in\{1,2\}$.
\end{proof}
\begin{lemma}\label{lem:trees_with_rank_zero}
    Let $\Gamma$ be a first-order reduct of $\mathcal{T}_n$ with no finite-range endomorphism.
    If $\sim_\Gamma$ is trivial, then $\omega.\Gamma$ is homomorphically equivalent to a first-order reduct of $(\mathbb{N},=)$.
\end{lemma}
\begin{proof}
    The proof is completely analogous to that of Lemma~\ref{lem:grids_with_rank_zero}; we include it here for the reader's convenience.
    Since $\Gamma$ does not have a finite-range endomorphism, by Lemma~\ref{lem:equivalence_of_finite_range_endomorphisms} neither does $\omega.\Gamma$.
    Hence, by Proposition~\ref{prop:endomorphism_with_image_in_relation} there exists $e\in\End(\omega.\Gamma)$ such that for all $a,b\in\omega\times V$ either $e(a)\sim_\Gamma e(b)$ or $|e(a)-e(b)|=\infty$.
    Fix some $a\in V$.
    Since $\sim_\Gamma$ is trivial, we can compose $e$ with an appropriate automorphism of $(\omega\times V,E)$ and get an endomorphism of $\omega.\Gamma$ whose image is contained in $\omega\times\{a\}$.
    Let us assume that already $e(\omega\times V)\subseteq\omega\times\{a\}$.
    Then $\omega.\Gamma$ is homomorphically equivalent to $\omega.\Gamma|_{\omega\times\{a\}}$.
    Since every permutation of $\omega\times\{a\}$ preserves the relations of $\omega.\Gamma|_{\omega\times\{a\}}$, this structure is isomorphic to a first-order reduct of $(\mathbb{N},=)$ \cite[see][]{bodirsky2008}.
\end{proof}
\begin{proposition}\label{prop:trees_hardness}
    Let $\Gamma$ be a first-order reduct of $\mathcal{T}_n$ with no finite-range endomorphism.
    If $\sim_\Gamma$ is non-trivial, then $\Csp(\Gamma)$ is $\Np$-hard.
\end{proposition}
\begin{proof}
    If $n=2$, then $\Gamma$ is isomorphic to a first-order reduct of $(\mathbb{Z},<)$.
    Let us assume $\Gamma$ itself is a first-order reduct of $(\mathbb{Z},<)$.
    Since $\Gamma$ does not have a finite-range endomorphism and $\sim_\Gamma$ is non-trivial, inspection of the proof of Theorem~\ref{thm:defining_successor} shows that $\mathbb{Q}.\Gamma$ is homomorphically equivalent to a first-order reduct of $(\mathbb{Q}\times\mathbb{Z},<)$ that pp-defines $\Succ$ or $\Dist_m$ for all $m\in\mathbb{N}$.
    Because $\mathbb{Q}.\Gamma$ has $\mathbb{Q}\times\mathbb{Z}\to\mathbb{Q}\times\mathbb{Z},(p,a)\mapsto(p,-a)$ as an automorphism, it cannot be homomorphically equivalent to a first-order reduct of $(\mathbb{Q}\times\mathbb{Z},<)$ that pp-defines $\Succ$.
    It follows by Lemma~\ref{lem:hardness_from_distances} that $\Csp(\Gamma)$ is $\Np$-hard.

    Suppose $n>2$.
    By Lemma~\ref{lem:equivalence_of_finite_range_endomorphisms}, $\omega.\Gamma$ does not have a finite-range endomorphism.
    Hence, by Lemmas~\ref{lem:avoiding_orbitals} and~\ref{lem:trees_possible_ranks} there are $a,b\in\omega\times V$ with $|a-b|\in 4\mathbb{N}_{>0}$ and $e(a)\neq e(b)$ for all $e\in\End(\omega.\Gamma)$.
    Let $c\in\omega\times V$ with $|a-c|=|b-c|=|a-b|/2\eqcolon m$.
    Then by Lemma~\ref{lem:trees_possible_ranks} we have $a\sim_\Gamma b\sim_\Gamma c$.
    It follows from Lemma~\ref{lem:pp_define_image_under_endomorphisms} that $\omega.\Gamma$ pp-defines $(a,b)^{\End(\omega.\Gamma)}$ and $(a,c)^{\End(\omega.\Gamma)}$ and that these relations are of the form $\Dist_S$ and $\Dist_T$, respectively, for finite $S,T\subseteq\mathbb{N}$.
    Note, that $0\notin S$, $m\in T$ and $2m\in S$.
    An argument similar to the proof of Lemma~\ref{lem:use_relation_to_reduce_to_finite_csp} shows that $\Csp(\omega.\Gamma|_A)\leq_\Ptime\Csp(\Gamma)$ where $A\coloneq\{d\in\omega\times V\mid|c-d|\in T\}$.
    Observe, that $(A,{\Dist_S}|_A)$ is an undirected graph with no loops.
    Moreover, since $n>2$ there exists $d\in A$ with $|a-d|=|b-d|=2m$ and thus $a$, $b$ and $d$ form a triangle in $(A,{\Dist_S}|_A)$ (Figure~\ref{fig:tree_with_triangle}).
    By \cite{hell1990} it follows that $\Csp(A,{\Dist_S}|_A)$ and thus $\Csp(\Gamma)$ are $\Np$-hard.
    \begin{figure}[ht]
        \centering
        \begin{tikzpicture}[
                rotate=90,
                edge/.style={
                    draw=black!50,
                },
                tri/.style={
                    dashed,
                },
                dist/.style={
                    dotted,
                },
                label/.style={
                },
                node/.style={
                    circle,
                    fill=white,
                    inner sep=0.1pt,
                },
            ]
            \coordinate (0) at (0,0);
            \foreach \a in {0,1,2} {
                \begin{scope}[rotate={\a*120},shift={(1,0)}]
                    \coordinate (0\a) at (0,0);
                    \draw[edge] (0) -- (0\a);
                    \foreach \b in {0,1} {
                        \begin{scope}[rotate={\b*90-45},shift={(1,0)}]
                            \coordinate (0\a\b) at (0,0);
                            \draw[edge] (0\a) -- (0\a\b);
                            \foreach \c in {0,1} {
                                \begin{scope}[rotate={\c*90-45},shift={(0.6,0)}]
                                    \coordinate (0\a\b\c) at (0,0);
                                    \draw[edge] (0\a\b) -- (0\a\b\c);
                                    \foreach \d in {0,1} {
                                        \begin{scope}[rotate={\d*90-45},shift={(0.6,0)}]
                                            \coordinate (0\a\b\c\d) at (0,0);
                                            \draw[edge] (0\a\b\c) -- (0\a\b\c\d);
                                        \end{scope}
                                    }
                                \end{scope}
                            }
                        \end{scope}
                    }
                \end{scope}
            }
            \draw[tri] (00001) to node[label,right] {$2m$} (02111);
            \draw[tri] (00001) to node[label,above] {$2m$} (01000);
            \draw[tri] (02111) to node[label,below] {$2m$} (01000);
            \draw[dist] (0) to node[label,above] {$m$} (00001);
            \draw[dist] (0) to node[label,below] {$m$} (01000);
            \draw[dist] (0) to node[label,below] {$m$} (02111);
            \node[node] (a) at (00001) {$a$};
            \node[node] (b) at (02111) {$b$};
            \node[node] (c) at (0) {$c$};
            \node[node] (d) at (01000) {$d$};
        \end{tikzpicture}
        \caption{
            The dashed triangle causes $\mathsf{NP}$-hardness.
        }
        \label{fig:tree_with_triangle}
    \end{figure}
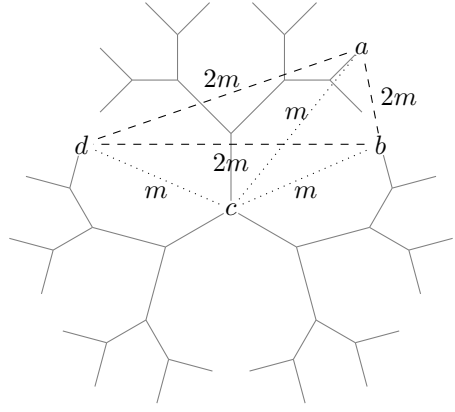
\end{proof}
\begin{theorem}
    Let $\Gamma$ be a finite-signature first-order reduct of $\mathcal{T}_n$.
    Then one of the following holds.
    \begin{enumerate}
        \item\label{item:trees_finite} $\Gamma$ is homomorphically equivalent to a finite structure and $\Csp(\Gamma)\in\Ptime\cup\Npc$.
        \item\label{item:trees_trivial} $\omega.\Gamma$ is homomorphically equivalent to a first-order reduct of $(\mathbb{N},=)$ and $\Csp(\Gamma)\in\Ptime\cup\Npc$.
        \item\label{item:trees_hard} $\Csp(\Gamma)\in\Npc$.
    \end{enumerate}
\end{theorem}
\begin{proof}
    If $\Gamma$ has a finite-range endomorphism, we are in Case~\ref{item:trees_finite} by \cite{bulatov2017} and \cite{zhuk2017,zhuk2020}, so suppose otherwise.
    Then if $\sim_\Gamma$ is trivial, by Lemma~\ref{lem:trees_with_rank_zero} and \cite{bodirsky2008} we are in Case~\ref{item:trees_trivial}.
    On the other hand, if $\sim_\Gamma$ is non-trivial we are in Case~\ref{item:trees_hard} by Lemma~\ref{lem:trees_csp_in_np} and Proposition~\ref{prop:trees_hardness}.
\end{proof}

\printbibliography
\end{document}